\documentclass[english]{iopart}

\usepackage{graphicx}
\usepackage{url}
\usepackage{iopams} 

\expandafter\let\csname equation*\endcsname\relax
\expandafter\let\csname endequation*\endcsname\relax

\usepackage{amssymb}
\usepackage[latin1]{inputenc}
\usepackage{amsmath}
\usepackage{epsfig}
\newcounter{fig}
\begin{document}

\title[Algebraic Geometry approach of Diagonals ]
{\Large  Diagonals of rational functions: from differential algebra to effective algebraic geometry (unabridged version)}

\vskip .3cm 

\author{Y. Abdelaziz$^\dag$  S. Boukraa$^\pounds$, C. Koutschan$^\P$,
J-M. Maillard$^\dag$}

\address{$^\dag$ LPTMC, UMR 7600 CNRS, 
Universit\'e Pierre et Marie Curie, Sorbonne Universit\'e, 
Tour 23, 5\`eme \'etage, case 121, 
 4 Place Jussieu, 75252 Paris Cedex 05, France} 

\address{$^\pounds$  LSA, IAESB,
 Universit\'e de Blida 1, Algeria}

\address{$^\P$  Johann Radon Institute for Computational and Applied Mathematics, RICAM, Altenberger Strasse 69, 
A-4040 Linz,  Austria}

\vskip .2cm 

\begin{abstract}

We show that the results we had obtained on diagonals of
nine and ten parameters families of rational functions 
using creative telescoping, yielding modular forms expressed as
pullbacked $\, _2F_1$ hypergeometric functions, can be obtained,
much more efficiently, calculating the $\, j$-invariant of
an elliptic curve canonically associated with the denominator of
the rational functions.  In the case where creative telescoping yields
pullbacked $\, _2F_1$ hypergeometric functions,
we generalise this result to other families of rational functions
in three, and even more than three, variables.
We also generalise this result to rational functions in more than three
variables when the denominator can be associated to an algebraic
variety corresponding to products of elliptic curves, foliation
in  elliptic curves. We also
extend these results to rational functions in three variables when
the denominator  is associated with a {\em genus-two curve such that its
Jacobian is a split Jacobian} corresponding to the product
of two elliptic curves. We sketch the situation where
the denominator of the rational function is associated with
algebraic varieties that are not of the general type,
having  an infinite set of birational automorphisms. We
finally provide some examples of rational functions
in more than three variables, where the telescopers have pullbacked $\, _2F_1$
hypergeometric solutions, the denominator corresponding to an algebraic variety
being not simply foliated in elliptic curves, but having  a selected elliptic curve
in the variety explaining the pullbacked $\, _2F_1$
hypergeometric solution.

\end{abstract}

\vskip .1cm



\vskip .4cm


\noindent {\bf AMS Classification scheme numbers}: 34M55, 
47E05, 81Qxx, 32G34, 34Lxx, 34Mxx, 14Kxx 

\vskip .2cm


\vskip .1cm

\section{Introduction}
\label{Introduction}

In a previous paper~\cite{DiagJPA,unabri}, using creative  telescoping~\cite{Koutschan}, we have
obtained  {\em diagonals}\footnote[2]{For the introduction of the concept of diagonals of
rational functions see~\cite{Christol84,Christol85,Christol369,Christol111,Purdue,Lipshitz,Denef,DiagSelected}.}
of nine and ten parameters families of rational functions, given by  (classical) {\em modular forms}
expressed as pullbacked $\, _2F_1$ hypergeometric functions\footnote[9]{The lattice Green functions
are the simplest examples of such diagonal of rational
functions~\cite{GlasserGuttmann,GoodGuttmann,LGF,LGF2,HeunJPA,malala}.}. The natural emergence of diagonal
of rational functions in lattice statistical mechanics is explained in~\cite{Short,Big}.  This
can be seen as the reason of the frequent occurrence of modular forms, Calabi-Yau operators
in lattice statistical
mechanics~\cite{CalabiYauIsing,DiffAlgGreen,High,Khi6,bo-ha-ma-ze-07b,perimeter}. In
another previous paper~\cite{HeunJPA,malala},
dedicated to Heun functions that are diagonals of  simple rational
functions, or only solutions of {\em telescopers}~\cite{Telescopers,Telescopers2} of simple rational
functions of three variables, but most of the time four variables, we have obtained
many solutions of order-three telescopers having squares of Heun functions as solutions
that turn out to be squares of pullbacked $\, _2F_1$ hypergeometric solutions
corresponding to {\em classical modular forms} and even {\em Shimura automorphic forms}~\cite{Takeuchi,Voight}, 
strongly reminiscent of periods of  {\em extremal rational surfaces}~\cite{Malmendier4,HigherMalmendier}, 
and other foliation of K3 surfaces in elliptic curves. In other words one finds experimentally that
the $\, _2F_1$ hypergeometric functions emerging in the calculation of diagonal of rational functions,
or of solutions of the telescopers of rational functions, seem to be only special $\, _2F_1([a,b], [c],x)$
hypergeometric functions with a selected set of parameters $\, [a,b], [c]$
(see the list (B.1) in Appendix B of~\cite{HeunJPA}, corresponding to classical modular
forms\footnote[1]{See for instance Felix Klein's connection of  the
  $\, _2F_1([1/12,5/12], [1],x)$ Gauss hypergeometric function
  with modular forms, for instance in the very pedagogical and heuristic
  paper~\cite{Maier1}.}, together with
a finite set of parameters, like $\, [7/24,11/24],[5/4]$,  corresponding to
Shimura automorphic forms~\cite{Takeuchi,Voight}),
pullbacked by selected  pullbacks. This last paper~\cite{HeunJPA}
also underlined the difference between the diagonal of a rational function ${\it Diag}(R)$, and the solutions of the
telescoper of the same rational function.
These results strongly suggested to find an algebraic geometry interpretation for all
these results, and, more generally, suggested to provide an {\em alternative algebraic geometry
  approach of the results emerging from the creative
  telescoping}\footnote[5]{The reader may refer to~\cite{Chyzak} for an extensive survey of ``creative telescoping'' approaches.}.
This is the purpose
of the present paper. In particular, we are going to show that most of these pullbacked $\, _2F_1$ hypergeometric
functions can be obtained efficiently though algebraic geometry
calculations, thus providing a more intrinsic algebraic geometry interpretation
of the  creative telescoping calculations which are typically
{\em differential algebra calculations}~\cite{Telescopers,Telescopers2,Chyzak,Lairez}.

Creative telescoping~\cite{Telescopers,Telescopers2,Chyzak,Zeilberger} is a methodology to deal with parametrized symbolic sums
and integrals that yields differential/recurrence equations for such expressions.
This methodology  became popular in computer algebra in the past twenty five years.
By ``telescoper''
of a rational function, say $\, R(x,y,z)$, we here refer to the output of the creative telescoping
program~\cite{Koutschan}, applied to the {\em transformed} rational function $\, \tilde{R} = R(x/y,y/z,z)/(yz)$.
Such a telescoper is a linear differential operator~$T$ in $ \, x$ and $ \, {{\partial}\over {\partial x}}$, such that
$\,T \, + {{\partial U}\over {\partial y}} \, +  {{\partial V}\over {\partial z}} $  annihilates
$\tilde{R}$, where $U, \,V$ are rational functions in $x,y,z$. In other words, the telescoper~$T$ represents a
linear ODE that is satisfied by ${\it Diag}(R)$. 

The paper is essentially dedicated
to {\em solutions of telescopers} of rational functions which
are {\em not necessarily diagonals} of rational functions. These solutions
correspond to periods~\cite{KontZagier} of algebraic varieties   
over some cycles which are not necessarily vanishing~\cite{Igusa}
cycles\footnote[5]{In french ``cycles \'evanescents''~\cite{Deligne,ChristolPicard}.} 
like in the case of diagonals of rational functions.
The reader interested in the connection between the process of taking diagonals, calculating telescopers,
and the notion of Periods, deRham cohomology (i.e. differential forms)
and other Picard-Fuchs equations can read in detail the thesis of Pierre Lairez~\cite{Lairez}
(see also~\cite{Lairez2}).
We just sketch some of these ideas in \ref{Griffiths}.

The purpose of this paper is not to give an introduction on creative telescoping~\cite{Telescopers,Telescopers2},
but to provide many pedagogical (non-trivial) examples of telescopers using\footnote[9]{One can obtain these telescopers
using Chyzak's algorithm~\cite{Chyzak1} or Koutschan's semi-algorithm~\cite{Koutschan,Koutschan3} (the termination is not proven). For
the examples displayed in this paper, Koutschan's package~\cite{Koutschan} is more efficient.}
extensively the  ``{\em HolonomicFunctions}'' 
Mathematica package~\cite{Koutschan}. 

Let us first recall the exact results of~\cite{DiagJPA,unabri}.

\vskip .1cm   

\section{Classical modular forms and diagonals of nine and ten parameters  family of rational functions }
\label{Introduction}

 \vskip .1cm  

 In a previous paper~\cite{DiagJPA,unabri}, using creative telescoping~\cite{Koutschan},
 we have obtained  diagonals of
nine and ten parameters families of rational functions, given by  (classical) modular forms expressed as
pullbacked $\, _2F_1$ hypergeometric functions. Let us recall these results. 
 
\subsection{Nine-parameters rational functions giving pullbacked $\, _2F_1$ hypergeometric functions for their diagonals}
\label{nine}

Let us recall the {\em nine-parameters} rational function in three variables $\, x$, $\, y$ and $\, z$:
\begin{eqnarray}
\label{Ratfoncplusplus}
\hspace{-0.98in}&&
\quad  \quad 
 {{1} \over {
a \,\, \,+ \, b_1 \, x \, + \, b_2 \, y \,  + \, b_3 \, z \,\,\,
 + \, c_1 \, y\, z \, + \, c_2 \, x \, z \,  + \, c_3 \, x\, y 
\, \,\, + \, \, d \, y^2 \, z \, \, \, + \, \, e \, z \, x^2 }}. 
\end{eqnarray}
Calculating\footnote[1]{Using the ``{\em HolonomicFunctions}'' 
  Mathematica package~\cite{Koutschan}.} the telescoper
of  this rational function (\ref{Ratfoncplusplus}),
one gets an {\em order-two} linear differential  operator annihilating 
the diagonal of the rational function (\ref{Ratfoncplusplus}).
The  diagonal of the rational function (\ref{Ratfoncplusplus})
can be written~\cite{DiagJPA,unabri} as a pullbacked hypergeometric function
\begin{eqnarray}
\label{2F15HypformAplusplus}
\hspace{-0.7in}&&\quad \quad \quad \quad \quad 
{{1} \over { P_4(x)^{1/4}}} \cdot \, 
 _2F_1\Bigl([{{1} \over {12}}, \, {{5} \over {12}}], \, [1],
 \, \, 1 \, - \, {{P_6(x)^2 } \over {P_4(x)^3}}\Bigr),
\end{eqnarray}
where $\,\, P_4(x)\,$ and $\, P_6(x)\,$ are two polynomials of degree 
four and six  in $\, x$ respectively. The Hauptmodul pullback
in (\ref{2F15HypformAplusplus}) has the form
\begin{eqnarray}
\label{P4P3oftheform}
\hspace{-0.98in}&& \quad \quad \quad \quad \quad \quad \quad \, \,\,\,
 {\cal H}  \,  \, = \,    \,  \,  {{1728} \over {j}}
 \,  \, = \,    \,  \,  \,     1 \,\,   - \, {{P_6(x)^2 } \over {P_4(x)^3}} 
 \,  \,  \,  \, = \,   \,  \,  \,  \,  
 {{ 1728 \cdot \, x^3 \, P_8(x)} \over { P_4(x)^3}}, 
\end{eqnarray}
where $\, P_8(x)$ is a polynomial of degree eight in $\, x$. 
Such  a pullbacked $\, _2F_1$ hypergeometric function (\ref{2F15HypformAplusplus}) corresponds 
to a {\em classical modular forms}~\cite{DiagJPA,unabri}.

\vskip .2cm

\subsection{Ten-parameters rational functions giving pullbacked $\, _2F_1$ hypergeometric functions for their diagonals.}
\label{ten}

Let us recall the {\em ten-parameters} rational function  in three variables $\, x$, $\, y$ and $\, z$:
\begin{eqnarray}
\label{Ratfoncplusplusplus}
\hspace{-0.98in}&&
 R(x, \, y, \, z)  \, \, \, = \, \,  
 \\
\hspace{-0.98in}&&
 {{1} \over {
a \, \, \,  + \, b_1 \, x \,  + \, b_2 \, y \,   + \, b_3 \, z \,\,\, 
 + \, c_1 \, y\, z \, + \, c_2 \, x \, z \, + \, c_3 \, x\, y
 \, \, \, + \, \, d_1 \, x^2 \, y \, \,   + \, \, d_2 \, y^2 \, z \,   + \, \, d_3 \, z^2 \, x  }}. 
\nonumber 
\end{eqnarray}
Calculating\footnote[2]{Using the ``{\em HolonomicFunctions}'' 
program~\cite{Koutschan}.} the telescoper  of  this rational function (\ref{Ratfoncplusplusplus}),
one gets an order-two linear differential  operator annihilating 
the diagonal of the rational function (\ref{Ratfoncplusplusplus}).
The  diagonal of the rational function (\ref{Ratfoncplusplusplus})
can be written~\cite{DiagJPA,unabri} as a pullbacked hypergeometric function
\begin{eqnarray}
\label{2F15HypformAplusplusplus}
\hspace{-0.7in}&&\quad \quad \quad \quad \quad 
{{1} \over { P_3(x)^{1/4}}} \cdot \, 
 _2F_1\Bigl([{{1} \over {12}}, \, {{5} \over {12}}], \, [1],
 \, \, 1 \, - \, {{P_6(x)^2 } \over {P_3(x)^3}}\Bigr),
\end{eqnarray}
where $\, P_3(x) \, $ and $\, P_6(x) \, $ are two polynomials of degree 
three and six  in $\, x$ respectively.
Furthermore, the Hauptmodul pullback in (\ref{2F15HypformAplusplusplus}) is seen to be of the form: 
\begin{eqnarray}
\label{oftheform}
  \hspace{-0.98in}&& \quad  \quad \, \quad  \quad \quad \quad \quad
  {\cal H}  \,  \, = \,    \,  \,  {{1728} \over {j}}
 \,  \, = \,    \,  \,\,    
 1 \, \,  - \, {{P_6(x)^2 } \over {P_3(x)^3}}  \, \, \, \, = \, \, \, \, \,
{{1728 \cdot \, x^3 \cdot \, P_9(x)} \over { P_3(x)^3 }}.  
\end{eqnarray}
where $\, P_9(x) \, $ is a polynomial of degree nine in $\, x$.
Again, (\ref{2F15HypformAplusplusplus}) corresponds to a {\em classical modular form}~\cite{DiagJPA,unabri}.

 \vskip .1cm  
 \vskip .2cm 

\section{Deducing creative telescoping results from effective algebraic geometry}
\label{deducing}

Obtaining the previous pullbacked hypergeometric results (\ref{2F15HypformAplusplus}) and (\ref{2F15HypformAplusplusplus})
required~\cite{DiagJPA,unabri} an accumulation of creative telescoping calculations, and a lot of ``guessing'' using all the symmetries
of the diagonals of these rational functions (\ref{Ratfoncplusplus}) and  (\ref{Ratfoncplusplusplus}).
We are looking for a more efficient and intrinsic way of obtaining these exact results.
These two pullbacked hypergeometric results  (\ref{2F15HypformAplusplus}) and (\ref{2F15HypformAplusplusplus}),
are essentially ``encoded'' by their {\em Hauptmodul} pullbacks  (\ref{P4P3oftheform})  and  (\ref{oftheform}),
or, equivalently, their corresponding {\em $\, j$-invariants}. 
The interesting question, which will be addressed in this paper, is whether it is possible to
canonically associate an elliptic curve with precisely $\, j$-invariants corresponding to these
Hauptmoduls  $\, {\cal H}  \,  \, = \,    \,  \,  {{1728} \over {j}}$.

\vskip .1cm

\subsection{Revisiting the pullbacked hypergeometric results in an algebraic geometry perspective.}
\label{revisiting}

One expects such an elliptic curve to correspond to the singular part of the rational function,
namely the {\em denominator} of the rational function. Let us recall that the diagonal of a rational function
is obtained through the multi-Taylor expansion of the rational function~\cite{Short,Big}
\begin{eqnarray}
\label{multiTaylor}
\hspace{-0.98in}&& \quad  \quad \, \quad  \quad \quad  \quad  \quad  \, 
R(x, \, y, \, z) \, \, \,= \, \, \,\,
  \sum_m  \sum_n  \sum_l \,\, a_{m, \, n, \, l} \cdot x^m \, y^n \, z^l, 
\end{eqnarray}
by extracting the "diagonal" terms, i.e. the  powers of the product $\, p \, = \, \, x y z$:
\begin{eqnarray}
\label{diagmultiTaylor}
\hspace{-0.98in}&& \quad  \quad \, \quad  \quad \quad  \quad \quad  \, 
Diag\Bigl(R(x, \, y, \, z)\Bigr) \, \,\, = \, \, \,
 \sum_m \,\,  a_{m, \, m, \, m} \cdot x^n. 
\end{eqnarray}
Consequently, it is natural to consider the algebraic curve corresponding to the
intersection of the surface corresponding to the vanishing condition  $\, D(x, \, y, \, z) \, = \, \, 0\, $
of the denominator $\, D(x, \, y, \, z)$ 
of these rational functions (\ref{Ratfoncplusplus}) and  (\ref{Ratfoncplusplusplus}), with the hyperbola
$\, p \, = \, \, x\, y\, z$ (where $\, p$ is seen, here, as a constant). This amounts, for instance, to eliminating
the variable $\, z$, substituting $\,\, z\, = \, \, p/x/y \, \, $ in  $\, D(x, \, y, \, z) \, = \, \, 0$.

\subsubsection{Nine-parameters case:}

In the case of the rational functions (\ref{Ratfoncplusplus}) this corresponds to the (planar) algebraic curve
\begin{eqnarray}
\label{alg1}
\hspace{-0.98in}&& \quad \quad  \quad  \quad  \quad
 a \,\, \,+  b_1 \, x \,\,  +  b_2 \, y \,\,   +  b_3 \, {{p} \over {x\, y}} \,\,\,\, 
+  c_1 \, y\, {{p} \over {x\, y}} \,\,  +  c_2 \, x \, {{p} \over {x\, y}} \,\,   + c_3 \, x\, y
 \nonumber \\
 \hspace{-0.98in}&& \quad \quad \quad \quad \quad \quad \quad \quad
+  d \, y^2 \, {{p} \over {x\, y}} \, \, \, \, +  e \,  {{p} \over {x\, y}} \, x^2 
\, \ \, = \, \, \, 0,                   
\end{eqnarray}
which can be rewritten as a (general, nine-parameters) {\em biquadratic}:
\begin{eqnarray}
\label{alg1bis}
\hspace{-0.98in}&& \quad \quad  \quad  \quad \quad 
 a \,x \, y \, \,\, + \, b_1 \, x^2 \, y \, + \, b_2 \, x\, y^2 \,  + \, b_3 \, p \,  \,
+ \, c_1 \, p \, \, y \,  \, + \, c_2 \, p \, x \, \, \,  + \, c_3 \, x^2\, y^2
\nonumber \\
\hspace{-0.98in}&& \quad \quad \quad \quad \quad \quad \quad  \quad  \quad  \quad 
\, \,\, +  d \, p \, \, y^2 \, \, \, + e \, p \, x^2 \, \, \ \, = \, \, \, 0.                   
\end{eqnarray}
Using formal calculations\footnote[1]{Namely using with(algcurves) in Maple, and, in particular, the command
  j$\_$invariant.} one can easily calculate the genus of the planar algebraic curve (\ref{alg1bis}),
and find that this  planar algebraic curve {\em is actually an elliptic curve} (genus-one). Furthermore, one can 
(almost instantaneously) find the exact expression of the $\, j$-invariant of this elliptic curve
as a rational function of the nine parameters $\, a, \, b_1, \, b_2, \, \cdots, \, e$
in (\ref{Ratfoncplusplus}). One actually finds that this $\, j$-invariant {\em is precisely the} $\, j$-invariant $\, j$
such that the Hauptmodul $\, {\cal H}  \,  \, = \,    \,  \,  {{1728} \over {j}}$ {\em is the exact expression}
(\ref{P4P3oftheform}). In other words, the classical modular form result (\ref{2F15HypformAplusplus})
could have been obtained, almost instantaneously, calculating the $\, j$-invariant of an elliptic curve
canonically associated with the denominator of the rational function (\ref{Ratfoncplusplus}).
The algebraic planar curve (\ref{alg1bis}) {\em corresponds to the most  general biquadratic of two variables},
which depends on nine homogeneous parameters. Such general biquadratic is well-known to be an elliptic curve
for {\em generic values} of the {\em nine parameters}\footnote[2]{So many results in integrable models correspond
to this most general biquadratic: the Bethe ansatz of the Baxter model~\cite{JphysSym,Quasi},
the elliptic curve foliating the sixteen-vertex model~\cite{Quasi}, 
so many QRT birational maps~\cite{Determinantal}, etc ...}.  

Thus, the nine-parameters exact result (\ref{2F15HypformAplusplus}) {\em can be seen as a simple
consequence of the fact that the most  general nine-parameters biquadratic is an elliptic curve}.

\subsubsection{Ten-parameters case:}

In the case of the  rational function (\ref{Ratfoncplusplusplus}), one must consider the (planar)
algebraic curve
\begin{eqnarray}
\label{alg2}
\hspace{-0.98in}&& \, \,  \quad \quad  \quad  \quad \quad 
a \, \, \,  +  b_1 \, x \, \,  +  b_2 \, y \,\,    +  b_3 \, {{p} \over {x\, y}} \,\,\, \, 
+  c_1 \, y\, {{p} \over {x\, y}} \,\,  +  c_2 \, x \, {{p} \over {x\, y}} \,\,  +  c_3 \, x\, y
\nonumber \\
\hspace{-0.98in}&& \quad \quad \quad \quad \quad \quad \quad  \,\,
\, \, \, +  d_1 \, x^2 \, y \, \,\,    + d_2 \, y^2 \, {{p} \over {x\, y}}
\,  \,  + d_3 \, {{p^2} \over {x^2\, y^2}}
  \, x  \, \,\,\, = \, \, \, 0, 
\end{eqnarray}
i.e. the  {\em ten-parameters bicubic}:
\begin{eqnarray}
\label{alg2bis}
\hspace{-0.98in}&&  \quad \quad  \quad  \quad \quad 
a \,x \, y^2 \, \,  + \, b_1 \, x^2 \, y^2 \,  + \, b_2 \, x \, y^3 \,   + \, b_3 \, p  \, y \,\,\, 
                    + \, c_1 \,  p \, y^2 \, + \, c_2 \, \, p x \, y \, + \, c_3 \, x^2 \, y^3
\nonumber \\
\hspace{-0.98in}&& \quad \quad \quad \quad \quad \quad \quad \quad 
 \, \, \, + \, \, d_1 \, x^3 \, y^3 \, \,   + \, \, d_2 \, \, y^3
 \,   + \, \, d_3 \, p^2 \,  \, \, = \, \, \, 0. 
\end{eqnarray}
Using formal calculations, one can easily calculate the genus of this selected planar algebraic curve (\ref{alg2bis}),
and find that this  planar algebraic curve  {\em is actually an elliptic curve}\footnote[3]{Generically, 
  the most general planar bicubic is {\em not} a genus-one algebraic curve. It is a genus-four curve.}
(genus-one). Again one can 
 find\footnote[5]{For the bicubic (\ref{alg2bis}) the calculation of the $\, j$-invariant using the command
   j$\_$invariant  using with(algcurves) in Maple, requires much more computing time.} the exact expression
 of the $\, j$-invariant of this elliptic curve
as a rational function of the ten parameters $\, a, \, b_1, \, b_2, \, \cdots, \, d_3$
in (\ref{Ratfoncplusplusplus}). One actually finds that
this $\, j$-invariant {\em is precisely the} $\, j$-invariant $\, j$
such that the Hauptmodul $\, {\cal H}  \,  \, = \,    \,  \,  {{1728} \over {j}}$ {\em is the exact expression}
(\ref{oftheform}). In other words, the {\em classical modular form} result (\ref{2F15HypformAplusplusplus})
could have been obtained, much more simply, calculating the $\, j$-invariant of an elliptic curve
canonically associated with the denominator of the rational function (\ref{Ratfoncplusplusplus}).

Thus, this ten-parameters result (\ref{2F15HypformAplusplusplus}) can again be seen as a simple
consequence of the fact that {\em there exists a family of ten-parameters bicubics} (see (\ref{alg2bis}))
{\em which are elliptic curves for generic values of the ten parameters}.

These preliminary calculations are a strong incentive to try to replace the
differential algebra calculations of the {\em creative telescoping}, by more intrinsic
algebraic geometry  calculations, or, at least, perform effective algebraic geometry  calculations to provide an
algebraic geometry interpretation of the exact results obtained from creative telescoping.

\subsection{Finding creative telescoping results from $\, j$-invariant calculations.}
\label{finding}

One might think that these results are a consequence of the simplicity of the denominators
of the rational functions (\ref{Ratfoncplusplus}) or (\ref{Ratfoncplusplusplus}),
being associated with biquadratics or selected bicubics.

 Let us consider a nine-parameters family of planar algebraic curves that are
 not biquadratics or (selected) bicubics:
\begin{eqnarray}
\label{notbicubics}
  \hspace{-0.98in}&& 
\,  \,  \,  \, \, \,
a_1 \, x^4 \, + a_2 \, x^3 \,  \, +  a_3 \, x^2  \, +   a_4 \, x  \, + a_5
+a_6 \, x^2 \, y \, +   a_7  \, y^2 \, + a_8  \, y \,  \,   + a_9 \, x \, y \,  \, \, = \, \, \, 0.                 
\end{eqnarray}
One can easily calculate the genus of this planar curve and see that this
genus is actually one for arbitrary values of the $\, a_n$'s.
Thus the planar curve (\ref{notbicubics}) {\em  is an elliptic curve
for generic values of the nine parameters} $\, a_1, \, \cdots \, , a_9$.
It is straightforward to see that the algebraic surface $\, S(x, \, y, \, z) \, = \, \, 0$, corresponding to  
\begin{eqnarray}
\label{notbicubicsz}
  \hspace{-0.98in}&&
\,  \,  \,  \,  \,  \,
a_1 \, x^4 \, + a_2 \, x^3 \, + a_3 \, x^2  \, +   a_4 \, x  \, + a_5
+a_6 \, x^2 \, y \, +   a_7  \, y^2 \, +  a_8  \, y \,   +  a_9 \, {{p} \over {z}}  \,\,  \, = \,\,  \, 0,                 
\end{eqnarray}
or
\begin{eqnarray}
\label{notbicubicsz}
  \hspace{-0.98in}&& 
 \,  z \cdot \, (a_1 \, x^4 \, +  a_2 \, x^3 \, +  a_3 \, x^2  \, +  a_4 \, x  \, +\, a_5
  +  a_6 \, x^2 \, y \, +    a_7  \, y^2 \, +  a_8  \, y) \,  \,  + \,  a_9 \, p  \, \, = \,\,  \, 0,                 
\end{eqnarray}
will automatically be such that its intersection with the hyperbola
$\, p \, = \, \, x\, y\, z \, $ gives back the elliptic curve
(\ref{notbicubics}).

Using this kind of ``reverse engineering'' yields to consider
the rational function in three variables $\, x$, $\, y$ and $\, z$
\begin{eqnarray}
\label{notbicubicszratio}
 \hspace{-0.98in}&&  \quad \,  \, \,
                    R(x, \, y, \, z) \, \, = \, \, \,
\nonumber \\
\hspace{-0.98in}&&  \quad        \,  \,  \,  \,  \,
{{ 1} \over { 1 \,  \, + \, z \cdot \,
(a_1 \, x^4 \, + \, a_2 \, x^3 \, + \,   \, a_3 \, x^2  \, + \,  a_4 \, x  \, +\, a_5
 + \, a_6 \, x^2 \, y \,  + \,  a_7  \, y^2 \, + \,  a_8  \, y) }},                   
\end{eqnarray}
which will be such that
{\em its denominator is canonically associated with an elliptic curve}. Again we can
immediately calculate the $\, j$-invariant of that  elliptic curve.
If one calculates the telescoper of this eight-parameters family of  rational functions
(\ref{notbicubicszratio}), one finds  that this telescoper is an order-two linear differential
operator with pullbacked hypergeometric solutions of the form
\begin{eqnarray}
\label{notbicubicszratioSOLU}
 \hspace{-0.98in}&&  \quad \quad  \quad \quad  \quad \quad  \quad \quad  \quad
{\cal A}(x)        \cdot \,
_2F_1\Bigl([{{1} \over {12}}, \, {{5} \over {12}}], \, [1], \, {\cal H}\Bigr), 
\end{eqnarray}
where $\, {\cal A}(x)$ is an algebraic function and, where again, the pullback-Hauptmodul
$\,  {\cal H} \, = \, \, 1728/j$, {\em precisely corresponds to
  the $\, j$-invariant of the  elliptic curve}.
In \ref{Clebsch}, we give another example of a (planar) elliptic curve corresponding
to the {\em intersection of two quadrics}\footnote[8]{Intersections of quadrics are
  well-known to give elliptic curves~\cite{Knaff,Quasi}.} where, again, one can get the
(creative telescoping) pullbacked $\, _2F_1$
result from a simple calculation of a $\, j$-invariant. 

More generally, seeking for planar elliptic curves,  one can look for planar algebraic curves
\begin{eqnarray}
\label{notbicubicszratiosum}
  \hspace{-0.98in}&&  \quad \quad  \quad \quad  \quad \quad \quad  \quad \quad
 \sum_{n=0}^{n=N}  \sum_{m=0}^{m=M} \,  \,  a_{m, \, n} \cdot \, x^n \, y^m \, \,  \, = \, \,  \, \, 0, 
\end{eqnarray}
defined by the set of $\,  a_{m, \, n}$'s  which are equal to zero, apart of $\, {\cal N}$ homogeneous
parameters $\,  a_{m, \, n}\, $ being, as in (\ref{alg1bis}) or  (\ref{alg2bis}) or
(\ref{notbicubicsz}), {\em independent parameters}. Finding such an $\, {\cal N}$-parameters family
of (planar) elliptic curves automatically provides an $\, {\cal N}$-parameters family
of rational functions such that their telescopers have a pullbacked $\, _2F_1$ hypergeometric solution
we can simply deduce from the $\, j$-invariant of that elliptic curve.

{\bf Question}: Recalling section \ref{ten}, is it possible to find
families of such (planar) elliptic curves {\em which depend on more than ten independent parameters}?  

Before addressing this question, let us recall the concept
of {\em birationally equivalent elliptic curves}. Let us consider for example 
the following monomial transformation:
\begin{eqnarray}
\label{205}
\hspace{-0.98in}&& \, \,  \quad  \quad \quad \quad \quad \quad\quad
(x, \, y) \quad  \quad \longrightarrow \quad  \quad  \quad
(x^{12}\, y^{11}, \,\, x^{205} \, y^{188}).
\end{eqnarray}
Its compositional inverse is the monomial transformation:
\begin{eqnarray}
\label{188}
\hspace{-0.98in}&& \, \, \quad   \quad \quad \quad \quad \quad\quad
(x, \, y) \quad  \quad \longrightarrow \quad  \quad  \quad
 \Bigr({{ x^{188} } \over {y^{11} }}, \,  \, {{ y^{12} } \over {x^{205} }}   \Bigr).  
\end{eqnarray}
This  monomial transformation (\ref{205}) is thus a  {\em birational}\footnote[1]{This transformation
is rational and its compositional inverse is also rational (here monomial).} transformation.
A birational transformation transforms an elliptic curve, like (\ref{notbicubics}),
into another elliptic curve {\em with the same} $\, j$-invariant: these two
elliptic curves are called {\em birationally equivalent}.
In the case of the  birational and monomial transformation (\ref{205}),
the elliptic curve (\ref{notbicubics}) is changed
into\footnote[9]{One can easily verify for particular values of the $\, a_k$'s, using
with(algcurves) in Maple, that the $\, j$-invariants of (\ref{notbicubics})
and  (\ref{notbicubics22}) are actually equal.}:
\begin{eqnarray}
\label{notbicubics22}
  \hspace{-0.98in}&& \quad \quad \quad \, 
a_1 \, x^{48} \, y^{44} \, \, +  a_2 \, x^{36}\, y^{33}
\,\,  +  a_3 \, x^{24} \, y^{22} \,  \, +   a_4 \, x^{12} \, y^{11} \,  \, + a_5
\nonumber \\
\hspace{-0.98in}&& \quad \quad \quad  \quad \,\,\,\, 
+ a_6 \, x^{229} \, y^{210} \,\, +  a_7  \,x ^{410} \,  y^{376}
\, + a_8  \,x^{205} \, y^{108} \,   +  a_9 \, x^{217}  \, y^{199}  \, \, \, = \,\,  \, 0.                 
\end{eqnarray}
With this kind of birational monomial transformation (\ref{205}), we see that
one can find {\em families of elliptic curves} (\ref{notbicubics22})
{\em of arbitrary large degrees}  in $\, x$ and $\, y$. Consequently one can find
nine or ten parameters families of
rational functions of {\em arbitrary large degrees} yielding
pullbacked $\, _2F_1$ hypergeometric functions. There is no constraint on the degree
of the planar algebraic curves (\ref{notbicubics22}): {\em the only relevant question is
the question of the maximum number of (linearly) independent parameters of families of planar 
elliptic curves}. 
In fact, {\em it is possible to show that the maximum number of independent parameters is actually ten}.
We sketch the demonstration\footnote[5]{ We thank Josef Schicho for providing this demonstration.}
in \ref{MaxParam}.

\subsection{Pullbacked $\, _2F_1$ functions for higher genus curves: monomial transformations.}
\label{PullbackMonomial}

We have already remarked in~\cite{DiagJPA,unabri} that once we have an exact result
for a diagonal of a rational function $\, R(x, \, y, \, z)$, 
we immediately get another exact result for the diagonal of the  rational function
$ \, R(x^n, \, y^n, \, z^n)$
for any positive integer $\, n$. As a result we obtain a new expression
for the diagonal changing $\, x$ into $\, x^n$. In fact, this is also a result
on the telescoper of the rational function $\, R(x, \, y, \, z)$:
the telescoper of the rational function $ \, R(x^n, \, y^n, \, z^n)$
is the $\, x \, \rightarrow \, x^n$ pullback of the telescoper of the rational
function $\, R(x, \, y, \, z)$. Having a pullbacked $\, _2F_1$ solution
for the telescoper of the rational function $\, R(x, \, y, \, z)$
(resp. the diagonal of the rational function $\, R(x, \, y, \, z)$), we will immediately
deduce a pullbacked $\,\, _2F_1\,$ solution
for the telescoper of the rational function $ \, R(x^n, \, y^n, \, z^n)$
(resp. the diagonal of the rational function $ \, R(x^n, \, y^n, \, z^n)$).

Along this line, let us change in the rational function (\ref{Ratfoncplusplus}),
$\, (x, \, y, \, z) \, $ into  $\, (x^2, \, y^2, \, z^2)$: 
\begin{eqnarray}
\label{Ratfoncplusplusx2}
  \hspace{-0.98in}&& \,  \quad 
 R_2(x, \, y, \, z) \, \, = \, \,
 \\
 \hspace{-0.98in}&& \, \,                  
 {{1} \over {
a \,\, \,+ \, b_1 \, x^2 \, + \, b_2 \, y^2 \,  + \, b_3 \, z^2 \,\,\,
 + \, c_1 \, y^2 \, z^2 \, + \, c_2 \, x^2 \, z^2 \,  + \, c_3 \, x^2 \, y^2 
\, \,\, + \, \, d \, y^4 \, z^2 \, \, \, + \, \, e \, z^2 \, x^4 }}.
\nonumber 
\end{eqnarray}
The diagonal of this new rational function (\ref{Ratfoncplusplusx2}) will
be the pullbacked $\, _2F_1$ exact expression (\ref{2F15HypformAplusplus})
where we change $ \, \, x \, \rightarrow \, \, x^2$. The intersection
of the algebraic surface corresponding to the vanishing condition of the
denominator of the new rational function (\ref{Ratfoncplusplusx2}),
with the hyperbola $\, p \, = \, \, x\, y\, z \, $ (i.e. $\, z \, = p/x/y$), 
is nothing but the equation (\ref{alg1bis}) where we have
changed $\, (x, \, y; \, p) \, $ into   $\, \, (x^2, \, y^2; \, p^2)$
\begin{eqnarray}
\label{alg1bisp2}
\hspace{-0.98in}&& \quad \quad  \quad   \, \, 
 a \, \, x^2 \, y^2 \, \,+ \, b_1 \, x^4 \, y^2 \, + \, b_2 \,\,  x^2 \, y^4 \,  + \, b_3 \, p^2 \,  \,
+ \, c_1 \, \, p^2 \, \, y^2 \,  \, + \, c_2\,  \, p^2 \, x^2 \, \, \,  + \, c_3\,  \, x^4 \, y^4
\nonumber \\
\hspace{-0.98in}&& \quad \quad \quad \quad \quad \quad \quad \quad 
\, \,\, + \, d \, p^2 \, \, y^4 \, \, \, + \, e \, p^2 \, x^4 \, \, \ \, = \, \, \, 0,                   
\end{eqnarray}
which is {\em no longer}\footnote[2]{If we perform the same calculations
  with the ten-parameters rational function  (\ref{Ratfoncplusplusplus})
  we get an algebraic curve of genus $\, 10$ instead of $\, 9$.}
{\em an elliptic curve but a curve of genus} $\, 9$.

With that example we see that classical modular form results, or pullbacked
$\, _2F_1$ exact expressions like (\ref{2F15HypformAplusplus}), can actually emerge
from {\em higher genus curves} like (\ref{alg1bisp2}). As far as these
diagonals, or telescopers, of rational function
calculations are concerned, higher genus curves
like (\ref{alg1bisp2}) must in fact be seen as ``almost'' elliptic curves up to
a  $\, x \, \rightarrow \, \, x ^n \, $ covering.

Such results for monomial transformations like
$ \,  (x, \, y, \, z)\, \rightarrow \,  (x^n, \, y^n, \, z^n) \, $
can, in fact, be generalised to more general
(non birational\footnote[1]{In contrast with transformations like (\ref{205}).})
monomial transformations.
This is sketched in \ref{monomialsec}.

\vskip .1cm 

\subsection{Changing the parameters into functions of the product $\,  \,  p \, = \, x\, y\, z$.}
\label{findinga}

All these results for many parameters families of rational functions can be {\em drastically generalised}
when one remarks that allowing any of these parameters
to be a {\em function} of the product $\, p \, = \, x\, y\, z \, $ also
yields to the previous pullbacked $\, _2F_1$ exact expression, like (\ref{2F15HypformAplusplus}),
{\em where the parameter is changed into that function of} $\, x$ (see \cite{DiagJPA}).  Let
us consider a simple (two-parameters) illustration
of this general result. Let us consider a subcase of the previous nine or ten parameters families, 
introducing the two parameters rational function:
\begin{eqnarray}
\label{ratexample}
\hspace{-0.98in}&& \quad \quad  \quad  \quad \quad  \quad   \quad    \quad 
  {{1} \over {
1 \,\, \,\,+2 \, x \, + \, b_2 \cdot \, y \, 
 + \, 5 \, y\, z \, \, + x \, z \,\,  +c_3 \cdot \, x\, y  }}.
\end{eqnarray}
The  diagonal of this rational function (\ref{ratexample})
is the pullbacked hypergeometric function:
\begin{eqnarray}
\label{2F1ratexample}
\hspace{-0.7in}&& \quad \quad \quad \quad   \quad  
{{1} \over { P_2(x)^{1/4}}} \cdot \, 
 _2F_1\Bigl([{{1} \over {12}}, \, {{5} \over {12}}], \, [1],
 \, \,  43200 \cdot \, x^4 \cdot \, {{ P_4(x) } \over {P_2(x)^3}}\Bigr),
\end{eqnarray}
where
\begin{eqnarray}
\label{ratexamplewhere}
 \hspace{-0.98in}&& \quad  \quad 
P_2(x) \, \, = \, \, \,\,
 1 \, \,\, \,  -8\cdot \, (b_2 \, +10) \cdot \, x  \, \, \,  \,
 +8 \cdot (2\, b_2^2 \, -20\, b_2 \, + 15\, c_3 \, +200) \cdot \, x^2, 
\end{eqnarray}
and 
\begin{eqnarray}
\label{ratexamplewhere2}
 \hspace{-0.98in}&& \quad  \quad  \,
P_4(x) \, \, = \, \, \,
-675\, c_3^4 \cdot \, x^4  \, \, \, 
+4 \, c_3^2 \cdot \, (b_2 +10) \cdot \, (4\, b_2^2-100\, b_2+45\, c_3+400) \cdot \, x^3
\nonumber \\
\hspace{-0.98in}&&   \quad  \,  \quad   \quad    \quad             
\, \, \, \, +(64\, b_2^4-32\, b_2^3\, c_3-8\, b_2^2\, c_3^2-1280\, b_2^3+1280\, b_2^2\, c_3
\nonumber \\
\hspace{-0.98in}&&   \quad  \, \,   \quad   \, \quad  \,  \,  \,  \,  \,  \,   \quad  \quad    \quad 
-460\, b_2\, c_3^2-5\, c_3^3+6400\, b_2^2-3200\, b_2\, c_3-800\, c_3^2) \cdot \, x^2
 \\
\hspace{-0.98in}&&  \, \,  \quad   \, \quad   \quad    \quad \quad    \quad \quad    \quad            
\, \, -(b_2 \, +10) \cdot \, (32\, b_2^2-16\, b_2\, c_3-c_3^2) \cdot \, x   \,\,\,\,
+2\, b_2\, \cdot \, (2\, b_2 \, -\, c_3), \nonumber 
\end{eqnarray}

Let us now consider the previous rational function (\ref{ratexample}) where
the two parameters $\, b_2$ and $\, c_3$ {\em become some rational functions of the product}
$\, p \, = \, x\, y\, z$, for instance:
\begin{eqnarray}
\label{ratexampleB2c3}
\hspace{-0.98in}&&    \quad \, \,  \, 
b_2(p) \, = \,   \,  {{1 \, +3\,p} \over {1 \, +7\, p^2 }},   \quad \quad \, 
c_3(p) \, = \, \,   {{1 \, +p^2} \over {1 \, +2\, p }} 
\quad \quad   \quad \quad\hbox{where:}  \quad  \, \, \quad
p \, = \, x\, y\, z.  
\end{eqnarray}
The new corresponding rational function becomes more involved
but one can easily calculate the telescoper of this  new  rational function of three variables
$\, x$, $\, y$ and $\, z$,
and find that it is, {\em again}, an order-two linear differential operator
having  the pullbacked hypergeometric solution (\ref{2F1ratexample})
{\em where}  $\, b_2$ {\em and} $\, c_3$ {\em are, now, replaced by} ($\, p$ is now $\, x$)
the functions:
\begin{eqnarray}
\label{ratexampleB2c3x}
\hspace{-0.98in}&& \quad \quad  \quad  \, \, \,  \,   \quad \quad \quad \,\,\,
b_2(x) \, = \,   \,  {{1 \, +3\,x} \over {1 \, +7\, x^2 }},
\quad \quad  \, \,\, \quad
c_3(x) \, = \, \,   {{1 \, +x^2} \over {1 \, +2\, x }}.     
\end{eqnarray}
In that case (\ref{ratexample}) with (\ref{ratexampleB2c3}), one
gets a diagonal which is the pullbacked hypergeometric solution
\begin{eqnarray}
\label{inthatcase}
\hspace{-0.98in}&& \quad \quad \quad \quad \quad \, \,
 (1\, +2\, x)^{1/4} \cdot \,  (1\, +7\, x^2)^{1/4} \cdot \, q_8^{-1/4}
 \nonumber \\
\hspace{-0.98in}&& \quad \quad \quad \quad \quad \,  \quad \quad \quad \, \, 
 \times     \, \,
 _2F_1\Bigl([{{1} \over {12}}, \, {{5} \over {12}}], \, [1], \, \, \,
 {{43200 \cdot \, x^4 \cdot \,  (1\, +7\, x^2)^2 \cdot \, q_{20} } \over {
 (1\, +2\, x) \cdot \, q_8^3 }}   \Bigr),                      
\end{eqnarray}
where: 
\begin{eqnarray}
\label{inthatcasewhere}
  \hspace{-0.98in}&& \quad \quad \quad
q_8  \, \, = \,   \,  \,
 5880\,{x}^{8} +156800\,{x}^{7} +71400\,{x}^{6} +35330\,{x}^{5} +19985\,{x}^{4}
 \nonumber \\
 \hspace{-0.98in}&& \quad  \quad \quad\quad \quad \quad
+1332\,{x}^{3}+1390\,{x}^{2}-86\,x \, +1,
\nonumber \\
\hspace{-0.98in}&&   \quad \quad \quad
q_{20}  \, \, = \,   \,  \,
-1620675\,{x}^{20} +1234800\,{x}^{18} +158332230\,{x}^{17} +153642195\,{x}^{16}
 \nonumber \\
\hspace{-0.98in}&& \quad \quad \ \quad \quad \quad  \quad
+427157990\,{x}^{15} +344201585\,{x}^{14} +367632300\,{x}^{13} +293263834\,{x}^{12}
\nonumber \\
  \hspace{-0.98in}&& \quad \quad \quad  \quad \quad  \quad \, \,
 +229496405\,{x}^{11} \,
+188180096\,{x}^{10}+107454499\,{x}^{9}+51936025\,{x}^{8}
\nonumber \\
  \hspace{-0.98in}&& \quad \quad\quad \quad  \quad  \quad \, \,
+21019296\,{x}^{7}+6259829\,{x}^{6}+1645018\,{x}^{5} +266619 \,{x}^{4}
\nonumber \\
\hspace{-0.98in}&& \quad  \quad \quad \quad \quad  \quad  \, \quad \quad \quad
   +40629\,{x}^{3}-1110\,{x}^{2}\,  \, -127\,x \,\, +2,
\end{eqnarray}
which is nothing but (\ref{2F1ratexample}) (with (\ref{ratexamplewhere}) and (\ref{ratexamplewhere2}))
where $\, b_2$ and $\, c_3$ have been replaced by the functions (\ref{ratexampleB2c3x}). 
Similar calculations can be performed
for more general rational functions (\ref{Ratfoncplusplus}) or  (\ref{Ratfoncplusplusplus}),
{\em when all the (nine or ten) parameters are more involved rational functions.}

When performing our creative telescoping symbolic calculations
using the HolonomicFunctions package~\cite{Koutschan}, such results
may look quite impressive.  From the algebraic geometry viewpoint,
it is almost tautological\footnote[2]{An algebraic geometer will probably see this
  as a trivial remark: diagonalization is an algebraic procedure and nothing really happens
  to the coefficients. Therefore if one replaces the coefficients by anything else, one will find those
replaced coefficients in the end result.}, if one takes for granted the result
of our previous subsections \ref{revisiting} and \ref{finding}, namely that the  
pullbacked hypergeometric solution of the telescoper corresponds to the Hauptmodul $\, 1728/j$, where
$\, j$ is the $\, j$-invariant of the elliptic curve corresponding to the intersection of
the algebraic surface corresponding to the vanishing condition of the
denominator, with the hyperbola $\, p \, = \, x\, y\, z$:   
this  calculation of the $\, j$-invariant is performed for $\, p$ fixed, and arbitrary (nine or ten) parameters
$ \, a, \, b_1, \, \cdots$ . It is clearly possible to force the parameters to be functions\footnote[1]{
The functions should be rational
functions if one wants to stick with diagonals and telescopers of {\em rational} functions,
but the result remains valid for {\em algebraic functions}, or {\em even transcendental  functions}
with reasonable series expansions at $\, x\, = \, 0$.}
of $\, p$, the $\, j$-invariant being changed accordingly. Of course,
in that case, the parameters in the rational function
are the same functions but of the product $\, p \, = \, x\, y \, z$.

{\em One thus gets pullbacked hypergeometric solutions (classical modular forms) for an (unreasonably ...) large set
of rational functions in three variables, namely the families of rational functions}
(\ref{Ratfoncplusplus}) or  (\ref{Ratfoncplusplusplus}), {\em but where, now, the  nine or ten
parameters are  nine, or ten, totally arbitrary rational functions of the product} $\, \, p \, = \, x\, y\, z$.

\vskip .2cm

{\bf Remark:} When the  rational function depends on parameters, one can straightforwardly deduce 
the solutions of the {\em telescoper} of the rational function where the  parameters are changed into functions.
In this example (see (\ref{ratexample}), (\ref{2F1ratexample}), (\ref{ratexampleB2c3}), (\ref{inthatcase})),
the solution  (\ref{2F1ratexample}) or (\ref{inthatcase})  of the {\em telescoper} of the rational function
is actually the diagonal of the rational function. 

\vskip .1cm
\vskip .1cm

We see experimentally that changing  the parameters  of the rational function
into functions, actually works for {\em diagonals} of rational functions.
Let us sketch the demonstration.

\subsubsection{Sketching the demonstration.}
\label{Demo}

Let us introduce the multi-Taylor expansion of
the rational function (\ref{ratexample}) where $\, b_2$ and $\, c_3$
are parameters  (not functions
of the product $\, p \, = \, \, x  \, y \, z$):
\begin{eqnarray}
\label{multiTaylorcalR1}
\hspace{-0.98in}&& \quad  \quad \,  \quad  \quad  \quad  \, 
R(x, \, y, \, z) \, \, \,= \, \, \,\,
\sum_m  \sum_n  \sum_l \,\, a_{m, \, n, \, l}\Bigl(b_2, \, \, c_3\Bigr)
\cdot x^n \, y^m \, z^l. 
\end{eqnarray}
The diagonal of this  rational function  (\ref{ratexample}) reads:
\begin{eqnarray} 
\label{diagmultiTaylorbis}
\hspace{-0.98in}&& \quad  \quad \,  \quad \quad  \quad \quad  \, 
Diag\Bigl(R(x, \, y, \, z)\Bigr) \, \,\, = \, \, \,
 \sum_m \,\,  a_{m, \, m, \, m}(b_2, \, \, c_3) \cdot x^n. 
\end{eqnarray}
Let us assume that we have an exact closed expression
$\, {\cal E}(b_2, \, c_3; \, x)$ for
this diagonal (\ref{diagmultiTaylorbis}), like the previous pullbacked
hypergeometric functions (\ref{2F15HypformAplusplus})
or (\ref{2F15HypformAplusplusplus}) (or possibly some Heun functions~\cite{malala},
or more involved exact expressions, like Appell functions, Lauricella functions, ...).

Let us assume that the coefficients $\, a_{m, \, m, \, m}$,
seen as functions of $\, b_2$ and $\, c_3$, have a multi-Taylor
expansion in  $\, b_2$ and $\, c_3$:
\begin{eqnarray} 
\label{amnlb7c3}
\hspace{-0.98in}&& \quad  \quad \,  \quad \quad  \quad \quad  \quad  \quad
a_{m, \, m, \, m}(b_2, \, \, c_3) \, \, = \, \, \,
\sum_{M, \, N} \, A_{M, \, n} \cdot \, b_2^{M}  \, c_3^{N}.
\end{eqnarray}

Let us now  assume that $\, b_2$ and $\, c_3$ are {\em functions of
  the product} $\, p \, = \, \, x \, y \, z$ (or more generally functions with
Taylor series expansions at $\, p\, = \, 0$).
The rational function (\ref{ratexample}), where $\, b_2$ and $\, c_3$
are now functions of the product $\, p \, = \, \, x \, y \, z$,
has the multi-series expansion:
\begin{eqnarray}
\label{multiTaylorcalR}
\hspace{-0.98in}&& \quad  \quad \,  \quad  \quad  \quad  \, 
{\cal R}(x, \, y, \, z) \, \, \,= \, \, \,\,
\sum_m  \sum_n  \sum_l \,\,
a_{m, \, n, \, l}\Bigl(b_2(p), \, \, c_3(p)\Bigr)
\cdot \,  x^n \, y^m \, z^l. 
\end{eqnarray}

Let us  assume that these two functions $\, b_2(p)$
and $\, c_3(p)$ both have a Taylor series expansion near
$\, p \, = \, 0$.

Consequently the coefficients $\, a_{m, \, m, \, m}$
in the multi-Taylor expansion (\ref{multiTaylorcalR}) have a
Taylor series expansion near $\, p \, = \, 0$:
\begin{eqnarray} 
\label{amnlb7c3}
\hspace{-0.98in}&& \quad  \quad \,  \quad \quad  \quad \quad  \quad
a_{m, \, m, \, m}\Bigl(b_2(p), \, \, c_3(p)\Bigr) \, \, = \, \, \, \,
\sum_{q} \, \alpha^{(q)}_{M, \, n} \cdot \, p^{q}.
\end{eqnarray}
The diagonal of the  rational function  (\ref{multiTaylorcalR}) is
actually (\ref{diagmultiTaylorbis}) where
the two parameters $\, b_2$ and $\, c_3$ are
changed into two functions $\, b_2(x)$ and $\, c_3(x)$
(like (\ref{ratexampleB2c3x})):
\begin{eqnarray} 
\label{diagmultiTaylorbister}
\hspace{-0.98in}&& \quad  \quad \,  \quad \quad  \quad \quad  \, 
Diag\Bigl({\cal R}(x, \, y, \, z)\Bigr) \, \,\, = \, \, \,
 \sum_m \,\,  a_{m, \, m, \, m}\Bigl(b_2(x), \, \, c_3(x)\Bigr) \cdot x^n. 
\end{eqnarray}
This multi-series (\ref{diagmultiTaylorbister}) has a Taylor
series expansion  which can be seen to be the Taylor series expansion of
the exact closed expression
$\, \, {\cal E}(b_2(x), \, c_3(x); \, x)$. 

\vskip .2cm

Of course this demonstration can be generalised to an arbitrary number of parameters
and for an arbitrary numbers of variables.

\vskip .2cm

\section{Creative telescoping on rational functions of more than three variables associated with products or foliations of elliptic curves}
\label{productelliptic}

Let us show that such an algebraic geometry approach of the creative telescoping can be generalised
to rational functions of {\em more than three} variables, when the vanishing condition of the denominator
can be associated with {\em products of elliptic curves}, or more generally, algebraic varieties with
{\em foliations in elliptic curves}.

\vskip .2cm 

$\, \bullet$  The telescoper of the rational function in the {\em four variables} $\, x$, $\, y$, $\, z$ and $\, w$
\begin{eqnarray}
\label{ratsurface}
 \hspace{-0.98in}&&  \quad     \,     \quad \quad    \, \,
{{ x \, y \, z} \over { (1+z)^2 \,\, \, \,
 -x \cdot \, (1-x) \cdot \, (x \, -x\,y\,z\,w) \cdot \, y \cdot \, (1-y) \cdot \, (y \, -x \, y \, z \,w)  }},  
\end{eqnarray}
gives an order-three {\em self-adjoint} linear differential operator which is, thus,
the {\em symmetric square} of an order-two linear differential operator. This
 order-two linear differential operator has
the pullbacked hypergeometric solution: 
\begin{eqnarray}
\label{surfacesol}
\hspace{-0.98in}&&  \quad   \quad   \, \,   \,   \,   \quad
{\cal S}_1 \, \, = \, \,\, \,  (1\, -x\, +x^2)^{-1/4} \cdot \,
 _2F_1\Bigl([{{1} \over {12}}, \, {{5} \over {12}}], \, [1], \, \,
 {{ 27} \over {4}} \cdot \, {\frac {{x}^{2} \cdot \, (1 \, -x)^{2}}{ ({x}^{2}-x+1)^{3}}}\Bigr). 
\end{eqnarray}
This  pullbacked hypergeometric solution (\ref{surfacesol}) can also be simply written:
\begin{eqnarray}
\label{surfacesolrewritten}
\hspace{-0.98in}&&  \quad   \quad   \, \,   \,   \quad \quad  \quad \quad   \quad \quad    \quad
  _2F_1\Bigl([{{1} \over {2}}, \, {{1} \over {2}}], \, [1], \, \, x\Bigr). 
\end{eqnarray}
In~\cite{malala} we underlined the difference between the {\em diagonal} of a rational function and
{\em solutions of the telescoper} of the same rational function. In this case, the  diagonal
of the rational function  (\ref{ratsurface}),
is zero\footnote[1]{The reason is that the integration takes place over a cycle homologically equivalent to the
  trivial cycle. The cycle becomes trivial after taking the limit $\, p \, \rightarrow \, 0$. Integrals over non vanishing cycles
  usually give logarithms of $\, p$, like the second solution to the hypergeometric function
  $\, _2F_1([1/2,1/2],[1], , x)$.}
and is thus different from the  pullbacked hypergeometric
solution (\ref{surfacesol}), which is
a {\em ``Period''}~\cite{KontZagier} of the algebraic variety corresponding to the denominator over some 
({\em non-vanishing}\footnote[2]{Diagonals of the rational functions correspond
  to periods over {\em vanishing cycles}~\cite{Igusa,ChristolPicard}.}) 
{\em cycle}. From now, we will have a similar situation in most of the 
following examples of this paper. 

This example is a simple illustration of what we expect for {\em products of elliptic curves}, or algebraic varieties with
{\em foliations in elliptic curves}. Introducing the product $\, p \, = \, xyzw$,
the vanishing condition of the denominator
of the rational function (\ref{ratsurface}) reads the surface $\, S(x, \, y, \, z) \, = \, \, 0$:
\begin{eqnarray}
\label{ratsurface_S}
 \hspace{-0.98in}&&  \quad     \quad    \quad \quad    \, \,
 (1+z)^2 \, \, \,\,
 -x \cdot \, (1-x) \cdot \, (x \, -p) \cdot \, y \cdot \, (1-y) \cdot \, (y \, -p)  \, \,\, = \,\, \, \, 0.  
\end{eqnarray}
For fixed $\, p$ and fixed $\, y$, equation (\ref{ratsurface_S}) can be seen as an
algebraic curve
\begin{eqnarray}
\label{ratsurfaceS}
 \hspace{-0.98in}&&  \quad    \quad   \quad    \quad \quad    \, \,
 (1+z)^2 \, \, \, - \, \lambda \cdot \, x \cdot \, (1-x) \cdot \, (x \, -p)  \, \, \, = \, \, \, \, 0
  \\
 \hspace{-0.98in}&&  \quad     \quad  \quad    \quad \quad   \quad \quad    \quad \quad    
 \hbox{for:} \quad  \quad  \quad \quad  \quad  \quad 
\lambda \, = \,  \, y \cdot \, (1-y) \cdot \, (y \, -p).
\nonumber 
\end{eqnarray}
For  fixed $\, p$ and fixed $\, y$, $\, \lambda$ can be seen as a constant, the algebraic
curve (\ref{ratsurfaceS}) being an {\em elliptic curve} with an obvious Weierstrass form:
\begin{eqnarray}
\label{ratsurfaceSZ}
 \hspace{-0.98in}&&  \quad    \quad \quad    \quad     \quad     \quad    \quad \quad    \, \,
 Z^2 \, \, \, \, - \,  \, x \cdot \, (1-x) \cdot \, (x \, -p) \,  \, \, = \, \, \, \, 0.
\end{eqnarray}
The $\, j$-invariant of (\ref{ratsurfaceS}), or\footnote[2]{A shift $\, z \rightarrow \, \, z+1 \, $
 or a rescaling $\, \, z^2  \rightarrow \, \,  z^2/\lambda \, $ does not change the $\, j$-invariant of the
 Weierstrass elliptic form.} (\ref{ratsurfaceSZ}), is well-known and yields the Hauptmodul $\, {\cal H}$:
\begin{eqnarray}
\label{ratsurfaceSZHaupt}
 \hspace{-0.98in}&&  \quad     \quad    \quad  \quad \quad    \quad \quad    \, \,
{\cal H} \, \, = \, \, \,{{1728} \over {j}} \, \, = \, \, \, 
{{ 27} \over {4}} \cdot \, {\frac {{p}^{2} \cdot \, (1 \, -p)^{2}}{ ({p}^{2}-p+1)^{3}}}.
\end{eqnarray}
For fixed $\, p$ and fixed $\, x$, equation (\ref{ratsurface_S}) can be seen as an
algebraic curve
\begin{eqnarray}
\label{ratsurfaceSx}
 \hspace{-0.98in}&&  \quad    \quad   \quad  \quad    \quad \quad    \, \,
 (1+z)^2 \, \,\, \, - \, \mu \cdot \, y \cdot \, (1-y) \cdot \, (y \, -p)
 \, \, \, = \, \,\,  \, 0
  \\
 \hspace{-0.98in}&&  \quad    \quad  \quad   \quad  \quad    \quad    \quad \quad    \, \,
 \hbox{for:} \quad  \quad  \quad  \quad \quad  \quad
 \mu \, = \,  \, x \cdot \, (1-x) \cdot \, (x \, -p),
 \nonumber 
\end{eqnarray}
which is also an elliptic curve with an obvious Weierstrass form
and the {\em same}  Hauptmodul  (\ref{ratsurfaceSZHaupt}).

More generally, the rational function of the {\em four variables} $\, x$, $\, y$, $\, z$ and $\, w$
\begin{eqnarray}
\label{ratsurfacemoregenera}
 \hspace{-0.98in}&&     \, \,\quad  \quad  \quad     \quad    \, \, 
{{ x \, y \, z} \over { (1+z)^2 \, \, \, \,
 -x \cdot \, (1-x) \cdot \, (x \, -R_1(p)) \cdot \, y \cdot \, (1-y) \cdot \, (y \, -R_2(p))  }},
\end{eqnarray}
where   $\, p\, = \, \, x\,y\,z\,w$,
and where  $\, R_1(p)$ and  $\, R_2(p)$ are two arbitrary rational functions
of the product $\, p\, = \, \, x\,y\,z\,w$,
yields a telescoper which has an {\em order-four} linear differential operator which is the
{\em symmetric product}\footnote[5]{This paper belonging to the symbolic computation literature and not
  pure mathematics for algebraic geometers, we use the standard Maple (DEtools) terminology of symmetric powers and symmetric
  products of linear differential operators~\cite{Weil}. Note that "symmetric product" is not a proper mathematical name
  for this construction on the solution space; it is a homomorphic image of the tensor product.
  The (Maple/DEtools) reason for choosing the name symmetric$\_$product is the resemblance with the function symmetric$\_$power. }  
of two order-two linear differential operators having respectively
the pullbacked hypergeometric solutions (\ref{surfacesol}) where $\, x$ is replaced by
$\, R_1(x)$ and  $\, R_2(x)$. These two  hypergeometric solutions thus have the two Hauptmodul pullbacks: 
\begin{eqnarray}
\label{ratsurfaceSZHaupt1}
 \hspace{-0.98in}&&  \quad     \quad    \quad  \quad \quad  \quad    \, \,
{\cal H}_1 \, \, = \, \, \,{{1728} \over {j_1}} \, \, = \, \, \, 
{{ 27} \over {4}} \cdot \, {\frac {{R_1(p)}^{2} \cdot \, (1 \, -R_1(p))^{2}}{ (R_1(p)^{2}-R_1(p)+1)^{3}}},
 \\
\hspace{-0.98in}&&  \quad     \quad   \quad   \quad \quad  \quad    \, \,
{\cal H}_2 \, \, = \, \, \,{{1728} \over {j_2}} \, \, = \, \, \, 
{{ 27} \over {4}} \cdot \, {\frac {{R_2(p)}^{2} \cdot \, (1 \, -R_2(p))^{2}}{ (R_2(p)^{2}-R_2(p)+1)^{3}}}.                 
\end{eqnarray}
A solution of the telescoper of (\ref{ratsurfacemoregenera}) is thus the {\em product} of these two pullbacked
hypergeometric functions. Let us give  two simple illustrations of this general result, with the two next examples.  

\vskip 2mm

$\, \bullet$ The telescoper of the rational function in the {\em four variables} $\, x$, $\, y$, $\, z$ and $\, w$
\begin{eqnarray}
\label{surface2}
 \hspace{-0.98in}&&  \quad     \quad   \, \,
{{x \, y \, z} \over { (1+z)^2 \, \, \,\,
  -x \cdot \, (1-x) \cdot \, (x \, -x\,y\,z\,w)
  \cdot \, y \cdot \, (1-y) \cdot \, (y \, -x^2 \, y^2 \, z^2 \,w^2)  }},  
\end{eqnarray}
gives an {\em order-four} linear differential operator which is the {\em symmetric product}
of {\em two order-two} linear differential operators having respectively
the pullbacked hypergeometric solution (\ref{surfacesol}) and the solution
(\ref{surfacesol}) where $\, x $ has been changed into $\, x^2$: 
\begin{eqnarray}
\label{surfacesol2}
\hspace{-0.98in}&&  \quad  \quad   \quad    \quad
 (1\, -x^2\, +x^4)^{-1/4} \cdot \,
_2F_1\Bigl([{{1} \over {12}}, \, {{5} \over {12}}], \, [1], \, \,
 {{ 27} \over {4}} \cdot \, {\frac { {x}^{4} \cdot \, (1 \, -x^2)^{2}}{({x}^{4}-x^2+1)^{3}}}\Bigr). 
\end{eqnarray}

$\, \bullet$ The telescoper of the rational function in the {\em four variables} $\, x$, $\, y$, $\, z$ and $\, w$
\begin{eqnarray}
\label{surface3}
 \hspace{-0.98in}&&    \,    \, \,   \quad    \quad   \, \,
{{ x \, y \, z} \over { (1+z)^2 \, \, \,  \,  \,
 -x \cdot \, (1-x) \cdot \, (x \, -x\,y\,z\,w) \cdot \, y \cdot \, (1-y) \cdot \, (y \, - 3 \, x \, y \, z \,w)  }},  
\end{eqnarray}
gives an {\em order-four} linear differential operator which is the {\em symmetric product}
of two order-two operators having respectively the pullbacked hypergeometric
solution (\ref{surfacesol})
and the solution  (\ref{surfacesol}) where the variable $\, x $ has been changed into $\, 3\, x$: 
\begin{eqnarray}
\label{surfacesol3}
\hspace{-0.98in}&&  \,\, \, \,  \,   \, \,  
{\cal S}_2 \, \, = \, \, \, \,  (1\, -3\, x\, + 9\, x^2)^{-1/4} \cdot \,
_2F_1\Bigl([{{1} \over {12}}, \, {{5} \over {12}}], \, [1], \, \,
{{ 243} \over {4}} \cdot \, {\frac {{x}^{2} \cdot \, (1 \, -3\, x)^{2}}{(1\, -3\, x\, + 9\, x^2)^{3}}}\Bigr). 
\end{eqnarray}

\vskip .2cm

\subsection{Creative telescoping on rational functions of five variables associated with products or foliations of three elliptic curves}
\label{productthreeelliptic}

Let us, now, introduce the rational function
in {\em five} variables $\, x$, $\, y$, $\, z$, $\, v$ and $\, w$
\begin{eqnarray}
\label{volume}
 \hspace{-0.98in}&&    \, \quad    \quad    \quad  \quad  \quad  \quad    \quad     \quad  \quad    \quad   \, \,
 {{ x \, y \, z \, v} \over { D(x, \, y, \, z, \, v, \, w)}},
\end{eqnarray}
where the denominator $\, \,  D(x, \, y, \, z, \, v, \, w) \, $ reads: 
\begin{eqnarray}
\label{volumewhere}
  \hspace{-0.98in}&&\,
D_p  \,\,= \, \,
\\
\hspace{-0.98in}&& \,\, \, \, \,  \,  \, \, \,
(1+v)^2 \,   \,   \,   \,  \,  -x \cdot \, (1-x) \cdot \, (x \, -p) \cdot \, y \cdot \, (1-y) \cdot \,
(y \, - 3 \, p) \cdot \, z \cdot \, (1-z) \cdot \, (z \, - 5 \, p), 
\nonumber \\
\hspace{-0.98in}&&
\quad \quad  \quad  \quad  \quad  \quad \hbox{where:}   \quad \quad \quad \quad \quad \quad
        p\, \, = \, \, \,  x\,y\,z\,v \,w.        \nonumber 
\end{eqnarray}
 The telescoper of the rational function (\ref{volume}) of {\em five} variables
 gives\footnote[2]{Such a creative telescoping calculation requires ``some'' computing time to achieve the result ...}
 an {\em order-eight} linear differential operator which is the {\em symmetric product}
of {\em three order-two} linear differential operators having respectively the pullbacked hypergeometric
solution (\ref{surfacesol}),
the solution  (\ref{surfacesol}) where $\, x $ has been changed into $\, 3\, x$,
namely (\ref{surfacesol3}), and the solution  (\ref{surfacesol}),
where $\, x $ has been changed into $\,\, 5 \, x$:
\begin{eqnarray}
\label{surfacesol5}
\hspace{-0.98in}&&  \, \,    
{\cal S}_3 \, \, = \, \, \,   \,  (1\, -5\, x\, + 25\, x^2)^{-1/4} \cdot \,
_2F_1\Bigl([{{1} \over {12}}, \, {{5} \over {12}}], \, [1], \, \,
{{ 675} \over {4}} \cdot \,  {\frac {{x}^{2} \cdot \, (1 \, -5\, x)^{2}}{(1\, -5\, x\, + 25 \, x^2)^{3}}}\Bigr). 
\end{eqnarray}
In other words, the order-eight linear differential telescoper of the rational function  (\ref{volume})
has the {\em product} $\, {\cal S} \, = \, \, {\cal S}_1 \cdot \,  {\cal S}_2 \cdot \,  {\cal S}_3$,
of (\ref{surfacesol}),  (\ref{surfacesol3})  and  (\ref{surfacesol5})
as a solution. From an algebraic geometry viewpoint this is a consequence of the fact that,
for fixed $\, p$, the algebraic variety $\, D_p \, = \, \, 0$, where $\, D_p$ is given by
(\ref{volumewhere}), can be seen, for fixed $\, y$ and $\, z$, 
{\em as an elliptic curve} $\, {\cal E}_1$ of equation $\, D_{y,z,p}(v, \, x) \, = \, \, 0$,
 for fixed $\, x$ and $\, z$
 as an elliptic curve $\, {\cal E}_2$ of equation $\, D_{x,z,p}(v, \, y) \, = \, \, 0$,
 and  for fixed $\, x$ and $\, y$ also
 as an elliptic curve $\, {\cal E}_3$ of equation $\, D_{x,y,p}(v, \, z) \, = \, \, 0$,
 the $\, j$-invariants $\, j_k$, $\, k \, =\, 1, \, 2, \, 3 \, $  of  these three elliptic
 curves $\, {\cal E}_k$ yielding (in terms of $\, p$), precisely,
 the three Hauptmoduls $\, {\cal H}_k \, = \, \, 1728/j_k$
\begin{eqnarray}
\label{Haupt123}
\hspace{-0.98in}&&  \,  \, 
{{ 27} \over {4}} \cdot \, {\frac {{x}^{2} \cdot \, (1 \, -x)^{2}}{ ({x}^{2}-x+1)^{3}}}, \, 
\quad      {{ 243} \over {4}} \cdot \, {\frac {{x}^{2} \cdot \, (1 \, -3\, x)^{2}}{(1\, -3\, x\, + 9\, x^2)^{3}}},  \, 
\quad  {{ 675} \over {4}} \cdot \,  {\frac {{x}^{2} \cdot \, (1 \, -5\, x)^{2}}{(1\, -5\, x\, + 25 \, x^2)^{3}}},     
\end{eqnarray}
occurring as pullbacks in the three $\, {\cal S}_k$'s of the solution
$\, {\cal S} \, = \, \, {\cal S}_1 \cdot \,  {\cal S}_2 \cdot \,  {\cal S}_3 \-, $
of the telescoper of (\ref{volume}).

\vskip .1cm

\subsection{Weierstrass and Legendre forms}
\label{Weierstrass}

The telescoper of the rational function in three variables
\begin{eqnarray}
\label{Weier3}
  \hspace{-0.98in}&& \, \,  \quad \quad  \quad \quad \quad \quad \quad \quad
 {{ x\, y} \over { (1+y)^2 \,  \,  \,   \, -x \cdot \, (1-x)\cdot \, (x \, -x\, y\, z)}} , 
\end{eqnarray}
associated\footnote[3]{The diagonal extracts the terms function of the product
$\, p \, = \, x\, y\, z$
in the multi-Taylor series.} with the {\em elliptic curve} in a {\em Weierstrass form}:
\begin{eqnarray}
\label{Weier3form}
  \hspace{-0.98in}&& \, \,  \quad \quad \quad \quad \quad \quad \quad \quad
  (1 \, +y)^2 \,  \, \, \, \,  -x \cdot \, (1-x)\cdot \, (x \, - p)   \,\, \,\,   \, = \,  \,  \,\,  \, 0, 
\end{eqnarray}
is the order-two linear differential operator
\begin{eqnarray}
\label{Weier3formL2}
  \hspace{-0.98in}&& \, \,  \quad   \quad  \quad \quad \quad
\, L_2 \,   \,= \, \, \,  \,
-1 \, \,  \, \,  + \, 4 \cdot \, (1 - 2\,x) \cdot \, D_x
\, \,   \, + \, 4 \cdot \, x \cdot \, (1\, -x) \cdot \, D_x^2,
\end{eqnarray}       
which has the hypergeometric solution:
\begin{eqnarray}
\label{Weier3formsol}
  \hspace{-0.98in}&& \, \,  \quad    \,   \quad  \quad  \quad    \, \,
 _2F_1\Bigl([{{1} \over {2}}, \, {{1} \over {2}}], \, [1], \, \, x\Bigr) 
  \\
  \hspace{-0.98in}&& \, \,   \,    \,  \quad   \quad  \quad    \,  \,  \, \quad  \quad \quad
 \, \, = \, \, \,\, 
 (1\, -x \, +x^2)^{-1/4}         \cdot \,
_2F_1\Bigl([{{1} \over {12}}, \, {{5} \over {12}}], \, [1], \, \, {{27} \over {4}} \cdot \,
{{ x^2 \cdot \, (1\, -x)^2} \over {(1\, -x \, +x^2)^3 }} \Bigr).
\nonumber 
\end{eqnarray}
The elliptic curve (\ref{Weier3form}) has the Hauptmodul
\begin{eqnarray}
\label{Weier3formHaupt}
\hspace{-0.98in}&& \, \,  \quad \quad   \quad \quad  \quad  \quad  \quad \quad  \quad  \quad
{\cal H} \, \, = \, \, \,   \, {{27} \over {4}} \cdot \,
 {{ p^2 \cdot \, (1\, -p)^2} \over {(1\, -p \, +p^2)^3}}.      
\end{eqnarray}
in agreement with the pullback in (\ref{Weier3formsol}).

\vskip .1cm

\subsubsection{K3 surfaces as products or foliations of two elliptic curves.}
\label{K3surfaces}

\vskip .1cm 

All the previous examples of this section correspond to denominators which are
algebraic varieties that can be seen as {\em Weierstrass elliptic curves}
for fixed values of all the variables except two. Let us show that one also gets
simple telescopers for rational functions with denominators which are
{\em algebraic varieties with some foliation in elliptic curves}\footnote[2]{Like K3 surfaces, or
three-fold Calabi-Yau manifolds.}.

\vskip .1cm 

$\, \bullet$ The telescoper of the rational function in {\em four} variables
\begin{eqnarray}
\label{K3}
  \hspace{-0.98in}&& \, \,   \,   \quad \quad \quad \quad  \quad   \quad
 {{ x\, y \, z} \over {
 (1+z)^2 \,  \,  \,   \,   \, -x \cdot \, (1-x) \cdot \, y \cdot \, (x \, -y) \cdot \, (y \, -x\, y\, z \, w) }},  
\end{eqnarray}
associated with the {\em $\, K_3$ surface} written in a
{\em Legendre form}\footnote[4]{Along this line see the first equation
  page 19 of~\cite{Algreen}.}
\begin{eqnarray}
\label{K3form}
  \hspace{-0.98in}&& \, \,  \quad  \quad \quad   \quad \quad\quad
 (1 \, +z)^2 \,   \,   \,  \, -x \cdot \, (1-x) \cdot \, y \cdot \, (x \, -y) \cdot \, (y \, - p)
 \, \, \, \,   \,  = \, \, \,  \,   \,\, 0, 
\end{eqnarray}
is an order-three {\em self-adjoint}\footnote[9]{The order-three
  linear differential operator is thus the symmetric square of an order-two linear differential operator.}
 linear differential operator $\, L_3$
\begin{eqnarray}
\label{L3sol}
  \hspace{-0.98in}&&  \quad  \quad \quad \quad \quad \quad \quad \quad
L_3 \, \,  \, = \,  \, \, \,  \,   x \cdot \, (2\, \theta \, +1)^3 \, \,  \,  \, - \, 8 \cdot \, \theta^3, 
\end{eqnarray}
which has the following $\, _3F_2$ solution (which is also,
because of a Clausen formula the square of a $\, _2F_1$ function):
\begin{eqnarray}
\label{L3sol}
  \hspace{-0.98in}&& \, \,  \, \quad \quad \quad \quad \quad 
  _3F_2\Bigl([{{1} \over {2}}, \, {{1} \over {2}}, \, {{1} \over {2}}], \, [1, \, 1], \, \,  \, x\Bigr)
\, \, = \, \, \,\,  _2F_1\Bigl([{{1} \over {4}}, \, {{1} \over {4}}], \, [1], \, \, x\Bigr)^2.
\end{eqnarray}
The $\, K_3$ surface (\ref{K3form}) can be seen as
{\em associated with the product of two Weierstrass elliptic curves}\footnote[1]{$K_3$
  surfaces {\em are not abelian varieties}, but they are ``close'' to  abelian varieties:
  they can be seen as essentially products of two elliptic curves.} of
Hauptmoduls respectively:
\begin{eqnarray}
\label{2Haupt}
  \hspace{-0.98in}&& \, \,   \quad \quad \quad \, \, 
{\cal H}_x   \, \, = \, \, \,  \, {{27} \over {4}} \cdot \,
{{ p^2 \cdot \, (1\, -p)^2} \over {(1\, -p \, +p^2)^3}},  \, \,   \, \,   \quad \quad             
{\cal H}_y   \, \,   \,= \, \, \,   \,  {{27} \over {4}} \cdot \,
 {{ y^2 \cdot \, (1\, -y)^2} \over {(1\, -y \, +y^2)^3}}.                   
\end{eqnarray}
This order-three linear differential operator $\, L_3$ is the
{\em symmetric square} of the order-two linear differential operator
\begin{eqnarray}
\label{M2L3sol}
  \hspace{-0.98in}&& \, \,  \quad \quad \quad \quad
 \, \, \, M_2 \, = \, \, \,\, -1 \,\, \,  \,  + \, 8 \cdot \, (2 - 3\,x) \cdot \, D_x
   \, \,\,  \,  +  16 \cdot \, x \cdot \, (1\, -x) \cdot \, D_x^2,
\end{eqnarray}
 which has the hypergeometric solutions:
\begin{eqnarray}
\label{K3formsol}
  \hspace{-0.98in}&& \, \,  \,  \,  \quad 
 _2F_1\Bigl([{{1} \over {4}}, \, {{1} \over {4}}], \, [1], \, \, x\Bigr) \, \, = \, \, \,
 \Bigl(1\, -{{x} \over {4}} \Bigr)^{-1/4}  \cdot \,
  _2F_1\Bigl([{{1} \over {12}}, \, {{5} \over {12}}], \, [1],
 \, \, -\,  {{ 27 \cdot \, x^2 } \over {(x \, -4)^3 }} \Bigr).
\end{eqnarray}
One thus finds that the telescoping procedure  associates to the  $\, K_3$ surface,
``encoded'' by $\, ({\cal H}_x, \, {\cal H}_y )$,  the Hauptmodul
given in  (\ref{K3formsol}):
\begin{eqnarray}
\label{whattelescopdoes}
\hspace{-0.98in}&& 
\Bigl({{27} \over {4}} \cdot \,
{{ p^2 \cdot \, (1\, -p)^2} \over {(1\, -p \, +p^2)^3}}, \, \, \,\, \,  {{27} \over {4}} \cdot \,
{{ y^2 \cdot \, (1\, -y)^2} \over {(1\, -y \, +y^2)^3}} \Bigr)
\, \, \,   \longrightarrow \, \,   \, \,\,  \, 
\Bigl(-\,  {{ 27 \cdot \, p^2 } \over {(p \, -4)^3 }},
 \, -\,  {{ 27 \cdot \, p^2 } \over {(p \, -4)^3 }}\Bigr). 
\end{eqnarray}

\vskip 3mm

\hskip -7mm{\bf Remark:}  The telescoper of
\begin{eqnarray}
\label{K3screw}
  \hspace{-0.98in}&& \, \,   \quad  \quad  \quad  \quad \quad \quad \quad
 {{ x\, y} \over {
 (1+z)^2 \,  \,  \,   \,  -x \cdot \, (1-x) \cdot \, y \cdot \, (x \, -y) \cdot \, (y \, -x\, y\, z \, w) }},  
\end{eqnarray}
{\em is a huge (48990 characters ...) order-eleven linear differential  operator}.  The telescoper of
\begin{eqnarray}
\label{K3screw2}
  \hspace{-0.98in}&& \, \,   \quad  \quad  \quad  \quad \quad \quad \quad
 {{ 1} \over {
 (1+z)^2 \,  \,  \,  \, -x \cdot \, (1-x) \cdot \, y \cdot \, (x \, -y) \cdot \, (y \, -x\, y\, z \, w) }},  
\end{eqnarray}
{\em is a huge (58702  characters ...) order-twelve linear differential  operator}.
  The telescoper of
\begin{eqnarray}
\label{K3screw3}
  \hspace{-0.98in}&& \, \,    \quad \quad  \quad \quad \quad \quad
 {{ x\, z \, w} \over {
 (1+z)^2 \,  \,  \,  \,  -x \cdot \, (1-x) \cdot \, y \cdot \, (x \, -y) \cdot \, (y \, -x\, y\, z \, w) }},  
\end{eqnarray}
{\em is a huge ( 59754  characters ...) order-eleven linear differential  operator}.
This raises the question of {\em how telescopers of rational functions are changed 
when  one modifies the numerator  of the rational function, keeping the same denominator}.
This is a quite involved question that we will address in forthcoming papers.

\vskip .2cm

$\, \bullet$ Let us now introduce the telescoper of the rational function in {\em four} variables
$\, x, \, y, \, z$ and $\, w$
\begin{eqnarray}
\label{K3BIS}
  \hspace{-0.98in}&& \quad  \quad  \quad  \quad  \quad 
 {{ x\, y \, z} \over {
 (1+z)^2    \, \, \, \,
  -x \cdot \, (1-x) \cdot \, y \cdot \, (x \, -y)  \cdot \, (y \, - 4 \cdot \, p \cdot \, (1 \, -p)) }},  
\end{eqnarray}
where $\, p$ denotes the product $\, p \, = \, x\, y\, z \, w$.
This rational function is nothing but (\ref{K3}) where $\, p$
has been changed into $\, 4 \cdot \, p \cdot \, (1 \, -p)$.
The telescoper  of the rational function of {\em four} variables (\ref{K3BIS})
is a {\em self-adjoint} order-three linear differential operator which is, thus, the {\em symmetric square}
of an order-two linear  differential operator. This order-two linear  differential operator
has the solution:
\begin{eqnarray}
\label{K3BISsol}
 \hspace{-0.98in}&& \quad  \quad  \quad  \quad
 _2F_1\Bigl([{{1} \over {2}}, \, {{1} \over {2}}], \, [1], \, \, x \Bigr)  
 \\
 \hspace{-0.98in}&& \, \, \, \,  \quad \quad \quad \quad \quad \, \, 
  \, \, = \, \, \,
 \Bigl(1\, -x \, +x^2 \Bigr)^{-1/4}  \cdot \,
  _2F_1\Bigl([{{1} \over {12}}, \, {{5} \over {12}}], \, [1],
 \, \, {{27} \over {4}} \cdot \,
 {{ x^2 \cdot \, (1\, -x)^2} \over {(1\, -x \, +x^2)^3}} \Bigr).
 \nonumber     
\end{eqnarray}
The relation of this result (\ref{K3BISsol}) with the previous result  (\ref{K3formsol})
corresponds to the following  identity for $\, \, X \, = \, 4 \cdot \, x \cdot \, (1\, -x)$:
\begin{eqnarray}
\label{K3formsolHauptymore}
\hspace{-0.98in}&& \, \,  \quad \, \,  \, \,  
_2F_1\Bigl([{{1} \over {4}}, \, {{1} \over {4}}], \, [1], \, \, X \Bigr) \, \, = \, \, \, \,
   _2F_1\Bigl([{{1} \over {4}}, \, {{1} \over {4}}], \, [1], \, \, 4 \, x \, (1\, -x) \Bigr) 
 \, \, = \, \, \,
 _2F_1\Bigl([{{1} \over {2}}, \, {{1} \over {2}}], \, [1], \, \, x \Bigr)
 \nonumber \\
 \hspace{-0.98in}&& \, \,  \quad \quad \quad \quad \quad  \, \, 
\, \, = \, \, \,
 \Bigl(1\, -{{X} \over {4}} \Bigr)^{-1/4}  \cdot \,
  _2F_1\Bigl([{{1} \over {12}}, \, {{5} \over {12}}], \, [1],
 \, \, -\,  {{ 27 \cdot \, X^2 } \over {(X \, -4)^3 }} \Bigr)
  \\
 \hspace{-0.98in}&& \, \,  \quad \quad \quad \quad \quad \, \, 
  \, \, = \, \, \,
 \Bigl(1\, -x \, +x^2 \Bigr)^{-1/4}  \cdot \,
  _2F_1\Bigl([{{1} \over {12}}, \, {{5} \over {12}}], \, [1],
 \, \, {{27} \over {4}} \cdot \,
 {{ x^2 \cdot \, (1\, -x)^2} \over {(1\, -x \, +x^2)^3}} \Bigr).
 \nonumber
\end{eqnarray}

\vskip .2cm 

We thus get exactly the same solution (\ref{surfacesol}) or (\ref{surfacesolrewritten})
than the one for the telescoper of the rational function (\ref{ratsurface}), where
the algebraic surface, corresponding
to the vanishing condition of the denominator, was clearly the product of two identical 
 elliptic curves with the {\em same} Hauptmodul (\ref{ratsurfaceSZHaupt}). 

\vskip .2cm 

{\bf Question:} Could it be possible that the two algebraic {\em surfaces}
\begin{eqnarray}
\label{D1}
\hspace{-0.98in}&& \, \, \quad  \quad \quad \, \, 
 (1+z)^2    \,  \,  \,  \, 
  -x \cdot \, (1-x) \cdot \, y \cdot \, (x \, -y)  \cdot \,
 \Bigl(y \, - 4 \cdot \, p \cdot \, (1 \, -p)\Bigr) \, \, = \, \,  \, 0, 
\end{eqnarray}
and
\begin{eqnarray}
\label{ratsurfaceSbis}
 \hspace{-0.98in}&&  \quad     \quad    \quad \quad    \, \,
 (1+z)^2 \, \, \,\, \, 
 -x \cdot \, (1-x) \cdot \, (x \, -p) \cdot \, y \cdot \, (1-y) \cdot \, (y \, -p)
 \, \,\, = \,\, \, \,  \, 0.  
\end{eqnarray}
be\footnote[1]{For algebraic curves, the situation is simpler since
two elliptic curves are {\em birationally equivalent} if and ony if they have the 
{\em same} j-invariant.} birationally equivalent?

\vskip .2cm 

\subsubsection{ Calabi-Yau three-fold manifolds  as foliation in three elliptic curves.}
\label{CalabiYau3fold}

 The telescoper of the rational function in {\em five} variables
$\,\, x, \, y, \, z, \, v$ and $\, w$
\begin{eqnarray}
\label{CalabiYauEEE}
  \hspace{-0.98in}&& \, \,  \quad \quad  \quad 
 {{ x\, y \, z \, v} \over {
 (1+w)^2 \, \,  \, \, \,  -x \cdot \, (1-x) \cdot \, y \cdot \, (x \, -y)
 \cdot \, z  \cdot \, (y \, -z) \cdot \, (z \, -x\, y\, z \, v \, w)     }}, 
\end{eqnarray}
associated\footnote[9]{The diagonal extracts the terms function of the product
$\, p \, = \, x\, y\, z \, v \, w \, $
in the multi-Taylor series.} with the {\em Calabi-Yau three-fold} written in a {\em Legendre form}
\begin{eqnarray}
\label{CalabiYauform}
  \hspace{-0.98in}&& \, \,  \,  \,  \, \quad \quad \, \, 
 (1 \, +w)^2 \, \, \,  \,  -x \cdot \, (1-x) \cdot \, y \cdot \,
(x \, -y) \cdot \, \, z  \cdot \,  (y \, - z)  \cdot \, (z \, -p)
\, \, \,\, = \,\,\,  \, 0, 
\end{eqnarray}
is an order-four (self-adjoint) linear differential  operator $\, L_4$
\begin{eqnarray}
\label{CalabiYauOper}
\hspace{-0.98in}&& \, \,  \quad \quad   \quad \quad  \quad \quad  \quad  \quad \quad
L_4 \,  \, \, = \, \,\, \,    16 \cdot \, \theta^4 \,\,  \, \, \, -x \cdot \, (2\, \theta \, +1)^4, 
\end{eqnarray}
which is a {\em Calabi-Yau operator}\footnote[5]{This linear differential operator is self-adjoint, its exterior
  square is of order five, it is MUM (maximum unipotent monodromy~\cite{Almkvist,Straten,IsingCalabi}), ...}
with the $\, _4F_3 \, $ solution:
\begin{eqnarray}
\label{CalabiYauformsol}
  \hspace{-0.98in}&& \, \,  \quad \quad \quad  \quad  \quad \quad \quad \quad 
 _4F_3\Bigl([{{1} \over {2}}, \, {{1} \over {2}}, \, {{1} \over {2}}, \, {{1} \over {2}}],
   \, [1, \, 1, \, 1],\, \, \, x \Bigr).
\end{eqnarray}
For $\, y$ and $\, z$ fixed, the $\,$ Calabi-Yau three-fold (\ref{CalabiYauform})
is foliated in {\em genus-one} curves
\begin{eqnarray}
\label{foliated}
  \hspace{-0.98in}&& \, \,  \quad \quad \quad  \quad \quad \quad \quad 
     (1 \, +w)^2 \, \, \, \,  - \lambda \cdot \, x \cdot \, (1-x) \cdot \,
(x \, -y) \, \, \,\, = \,\,\,  \, 0,             
\end{eqnarray}
where $\, \lambda$ is the constant expression ($p$ is fixed): 
\begin{eqnarray}
\label{foliatedwhere}
  \hspace{-0.98in}&& \, \, \quad \quad \quad \quad \quad  \quad \quad  \quad \quad 
 \lambda \, \,  \, = \, \, \, \,   y \cdot \,\, z  \cdot \,  (y \, - z)  \cdot \, (z \, -p).
\end{eqnarray}          
The Hauptmodul of these {\em genus-one} curves is {\em independent of} $\, p$ and $\, z$, reading:
\begin{eqnarray}
\label{CalabformsolHaupt}
\hspace{-0.98in}&& \, \,  \quad  \quad  \quad \quad  \quad  \quad  \quad  \quad  \quad 
 {\cal H}_{y,z} \, \, = \, \, \,  {{27} \over {4}} \cdot \,
 {{ y^2 \cdot \, (1\, -y)^2} \over {(1\, -y \, +y^2)^3}}. 
\end{eqnarray}
Similarly  for $\, x$ and $\, z$ fixed, the $\,$ Calabi-Yau three-fold (\ref{CalabiYauform})
is foliated in {\em genus-one} curves
\begin{eqnarray}
\label{foliatedxz}
  \hspace{-0.98in}&& \, \, \quad \quad \quad \quad \quad  \quad \quad 
(1 \, +w)^2 \, \, \, \,  - \mu \cdot \, y \cdot \, (x \, -y)  \cdot \,
(y \, -z) \, \, \,\, = \,\,\,  \, 0,             
\end{eqnarray}
where $\, \mu$ is the constant expression ($p$ is fixed): 
\begin{eqnarray}
\label{foliatedwhere}
  \hspace{-0.98in}&& \, \, \quad \quad \quad \quad  \quad \quad \quad  \quad \quad 
\mu \, \, = \, \,\, \,  x \cdot \,\, z  \cdot \,  (1 \, - x)  \cdot \, (z \, -p).
\end{eqnarray}  
The {\em genus-one} curves (\ref{foliatedxz}) can be written in a simpler Weierstrass form:
\begin{eqnarray}
\label{foliatedxzY}
  \hspace{-0.98in}&& \, \,\quad \quad  \quad \quad  \quad \quad 
(1 \, +w)^2 \, \, \,\,  - \rho \cdot \, Y \cdot \, \Bigl(1 \, - Y\Bigr)  \cdot \,
\Bigl( Y \, - \, {{z} \over {x}}\Bigr) \, \, \,\, = \,\,\,  \, 0,             
\end{eqnarray}
where the constant  $\, \rho$ reads $\, \rho \,  = \, \mu \cdot \, x^3$,
and the variable $\, y$ has been rescaled into $\, Y \, =  \, y/x$. The
Hauptmodul of these {\em genus-one} curves (\ref{foliatedxz}) is the same as
the Hauptmodul of the {\em genus-one} curves (\ref{foliated}), and corresponds to
expression  (\ref{CalabformsolHaupt}) where $\, y$ has been changed into $\, z/x$ (see
the canonical form (\ref{foliatedxzY})), namely:
\begin{eqnarray}
\label{CalabformsolHauptxz}
\hspace{-0.98in}&& \, \,  \quad  \quad  \quad \quad  \quad  \quad  \quad  \quad  \quad 
 {\cal H}_{x,z} \, \, = \, \, \,  {{27} \over {4}} \cdot \,
 {{ x^2 \cdot \, z^2 \cdot \, (x \, -z)^2} \over {(x^2\, -x\,z \, +z^2)^3}}. 
\end{eqnarray}
Similarly  for $\, x$ and $\, y$ fixed,  the $\,$ Calabi-Yau three-fold (\ref{CalabiYauform})
is foliated in {\em genus-one} curves,
\begin{eqnarray}
\label{foliatedxy}
  \hspace{-0.98in}&& \, \,  \,  \,  \, \quad  \quad \quad \quad \quad  \quad \, \, 
 (1 \, +w)^2 \, \, \,  \,  - \, \nu \cdot \, \, z  \cdot \,  (y \, - z)  \cdot \, (z \, -p)
\, \, \,\, = \,\,\,  \, 0, 
\end{eqnarray}
where $\, \nu$ reads:
\begin{eqnarray}
\label{foliatedxywhere}
  \hspace{-0.98in}&& \, \,  \,  \,  \, \quad \quad \quad  \quad \quad \quad \quad \quad   \, \, 
 \nu \, \,  \, = \, \,  \,  \, x \cdot \, (1-x) \cdot \, y \cdot \, (x \, -y).
\end{eqnarray}
A reduction to a canonical Weierstrass form similar to (\ref{foliatedxzY})
gives immediately the Hauptmodul of the {\em genus-one} curve (\ref{foliatedxy}) which reads: 
\begin{eqnarray}
\label{CalabformsolHauptxy}
\hspace{-0.98in}&& \, \,  \quad  \quad  \quad \quad  \quad  \quad  \quad  \quad  \quad 
 {\cal H}_{x,y} \, \, = \, \, \,  {{27} \over {4}} \cdot \,
 {{ y^2 \cdot \, p^2 \cdot \, (y \, -p)^2} \over {(y^2\, -y\,p \, +p^2)^3}}. 
\end{eqnarray}
The  {\em Calabi-Yau three-fold} (\ref{CalabiYauform}) thus has
a foliation in a triple of elliptic curves
$\, {\cal E}_1$, $\, {\cal E}_2$ and $\, {\cal E}_3$.

\vskip .2cm

\section{Creative telescoping of rational functions in three variables associated with genus-two curves
with split Jacobians}
\label{split}

\vskip .2cm 

In a paper~\cite{HeunJPA,malala}, dedicated to Heun functions that are solutions of telescopers of simple rational
functions of three and four variables,  we have 
obtained\footnote[1]{See equation (83) in section 2.2 of~\cite{malala}.}
an order-four telescoper
of a rational function of {\em three} variables, which is
the {\em direct sum of two order-two linear differential operators}, each having {\em classical modular forms}
solutions which can be written as  pullbacked
$\, _2F_1$ hypergeometric solutions. Unfortunately, the intersection of the
algebraic surface corresponding to the denominator 
of the rational function with the $ \, p \, = \, x\, y\, z \,\, $ hyperbola,
yields a {\em genus-two} algebraic curve. Note that this is a ``true'' genus-two curve:
it does not correspond to
the ``almost genus-one curves'' situation mentioned in subsection \ref{PullbackMonomial}. 

Let us try to understand, in this section, {\em how a genus-two curve can yield
two classical modular forms}. Let us first recall the results in section 2.2 of~\cite{malala}. 

\vskip .1cm

\subsection{Periods of extremal rational surfaces}
\label{subthree2}
Let us recall the rational function in just {\em three} variables~\cite{malala}:
\begin{eqnarray}
\label{Ratfonc4}
  \hspace{-0.98in}&& \quad  \quad \quad \quad \quad 
\, \, 
R(x, \, y, \, z)  \, \, \,  = \, \,  \quad 
 {{1} \over { 1 \, \,  \,\,  + x \, + y \, + z  \, \, \,+ x \,y \, + y \, z \,  \,  \, - x^3 \,y \,z }}.
\end{eqnarray}
Its telescoper is actually an {\em order-four} linear differential operator $\, L_4$
which, not only factorizes into {\em two order-two}  linear differential operators,
but is actually the {\em direct sum} (LCLM) of {\em two}\footnote[2]{
These two order-two  linear differential operators $\, L_2$ and $\, M_2$ are {\em not} homomorphic.}
{\em order-two}  linear differential operators
$\, L_4 \, = \, \, L_2 \oplus \, M_2$. These two (non homomorphic) order-two linear differential operators
have, respectively, the two pullbacked hypergeometric solutions:
\begin{eqnarray}
\label{Ratfonc4HeunSol}
\hspace{-0.98in}&& \quad  \quad
{\cal S}_1 \, \, = \, \, \, Heun \Bigl({{1} \over {2}}\, -{{i \sqrt{3}} \over {2}}, \,
{{1} \over {2}}\, -{{i \sqrt{3}} \over {6}}
, \, 1, \, 1, \, 1, \,1, \,\,\,
{{3} \over {2}} \cdot \, \Bigl(-3 \, + \, i \sqrt{3}  \Bigr) \cdot \,x \Bigr)
\end{eqnarray}
\begin{eqnarray}
\label{Ratfonc4HeunSolsuite}                 
\hspace{-0.98in}&& \quad   \quad \quad   \quad 
\, = \, \, \, \, \, 
  (1\, +9\, x)^{-1/4} \cdot \,      (1\, +3\, x)^{-1/4} \cdot \, (1\, +27 \, x^2)^{-1/4}
                  \nonumber     \\
\hspace{-0.98in}&& \quad \quad \quad \quad \quad \quad \quad  \quad \quad 
\times \,\,
  _2F_1\Bigl([{{1} \over {12}}, {{5} \over { 12 }}], \, [1], \, \,
 {{1728 \cdot \, x^3 \cdot \, (1 \, + 9 \, x \, + 27 \, x^2)^3 } \over {
  (1 \, +3 \, x)^3 \cdot \,   (1 \, +9\, x)^3  \cdot \, (1\, +27 \, x^2)^3 }} \Bigr),
\nonumber 
\end{eqnarray}
and: 
\begin{eqnarray}
\label{spurious}
\hspace{-0.98in}&& \quad  \quad  \quad \quad 
{\cal S}_2 \,\,  = \, \,  \,
{{1} \over { (1 \, + 4\, x \,-2\, x^2 \,-36\, x^3 \, + 81\, x^4 )^{1/4}}} \cdot \,
 \\
\hspace{-0.98in}&& \quad  \quad \quad \quad  \quad \quad   \quad        \times \, 
   \, _2F_1\Bigl([{{1} \over {12}}, {{5} \over { 12 }}], \, [1], \, \,
 {{1728 \cdot \, x^5 \, \cdot \, (1 +9 \, x \, +27 \, x^2)  \cdot \,(1 \, -2\,x)^2  } \over {
(1 \, + 4\, x \,-2\, x^2 \,-36\, x^3 \, + 81\, x^4)^3}}\Bigr).  \nonumber
\end{eqnarray} 
The diagonal of  (\ref{Ratfonc4}) {\em is actually the half-sum of the two series}
  (\ref{Ratfonc4HeunSol}) and (\ref{spurious}): 
\begin{eqnarray}
\label{spurioushalf}
  \hspace{-0.98in}&& \quad  \quad \quad  \quad \quad \quad  \quad \quad \quad 
Diag\Bigl( R(x, \, y, \, z)\Bigr)          \,   \, \, = \, \,  \, \,\,\,
 {{ {\cal S}_1 \,  + {\cal S}_2} \over { 2}}. 
\end{eqnarray}

As far as our algebraic geometry approach is concerned, the intersection of the
algebraic surface corresponding to the denominator of the rational function    (\ref{Ratfonc4})
with the hyperbola $\, p \, = \,\, \, x \, y\, z \, $ gives the planar algebraic curve
(corresponding to the elimination
 of the $\, z $ variable by the substitution $\, z = \, p/x/y$): 
\begin{eqnarray}
\label{planar}
\hspace{-0.98in}&& \quad  \quad  \quad  \quad \quad  \quad 
1 \, \,  \,\, \, + x \, + y \, \,  + {{ p} \over { x\, y}}
 \,   \, \, \,+ x \,y \,\, \,  + y \, {{ p} \over { x\, y}}
\,  \,  \,\, - x^3 \,y \,{{ p} \over { x\, y}}  \, \,\, \,  = \,  \, \,\, \, 0. 
\end{eqnarray}
One easily finds that this algebraic curve is (for $\, p$ fixed) a {\em genus-two} curve,
 and that this higher genus situation does not correspond to the "almost elliptic curves"
 described in subsection \ref{finding} namely an elliptic curve transformed by a
 monomial transformation. How a ``true'' {\em genus-two} curve can give two $\, j$-invariants,
namely  a telescoper
with two Hauptmodul pullbacked $\, _2F_1$ solutions?
We are going to see that the answer is that the Jacobian of this {\em genus-two}
curve\footnote[1]{An algebraic geometer will probably recall that it is very well-known that a genus two
curve {\em may} have Jacobian isogeneous to a product of elliptic curves. This is not the case in general.
The genus two curves that have a (nonconstant) map to an elliptic curve have this property. Our purpose
in section (\ref{split}) is to perform a creative telescoping calculation
in such a selected situation.
}  is in fact
isogenous to a product $\, {\cal E} \times \, {\cal E}' \,$ of
{\em two elliptic curves} (split Jacobian).

\vskip .2cm

\subsection{Split Jacobians}
\label{splitsub}

Let us first recall the concept of {\em split Jacobian} with a very simple example.
In~\cite{Kumar}, one has a crystal-clear example of a genus-two curve $\, C$
\begin{eqnarray}
\label{genustwo}
 \hspace{-0.98in}&&  \quad   \quad  \quad   \quad   \quad    \quad   \, \,
 y^{2} \, \,  \, - \, ({x}^{3}+420\,x-5600)  \cdot \, ({x}^{3}+42\,{x}^{2}+1120)
\, \,  \,\, =  \, \, \, \,  \,  0, 
\end{eqnarray}
such that its {\em Jacobian} $\, J(C)$ {\em is isogenous to a product of elliptic curves}
with $\, j$ invariants $\, j_1 \, = \, \, -2^7 \cdot \, 7^2 \, = \, \, -6272 \, $ and
$\, j_2 \, = \, \, -2^5 \cdot \, 7 \, \cdot \, 17^3 \, = \, \, -1100512$. These
two values correspond to the
following two values of the Hauptmodul $\, H\, = \, \, 1728/j$:   $\, H_1\, = \, \, -27/98$
and   $\, H_2\, = \, \, -54/34391$. 
Let us consider the {\em genus-one} elliptic curve
\begin{eqnarray}
\label{genusone}
  \hspace{-0.98in}&&   \quad   \, \quad   \quad  \quad  \quad   \quad    \, \,
   v^2 \, \, \, = \, \, \, u^{3} \, +4900\, u^{2} \, +7031500 \,u \, \, +2401000000, 
\end{eqnarray}
of $\, j$-invariant
$\, \, j \, = \, \, j_2 \, = \, \,-2^5 \cdot \, 7 \, \cdot \, 17^3 \, = \, \, -1100512$.
We consider the following
transformation\footnote[3]{This transformation
  is rational but {\em not birational}. If it were birational, then it 
would preserve the genus. Here, one goes from genus one to genus two.}:
\begin{eqnarray}
\label{morphism}
\hspace{-0.98in}&&    \, \quad \quad \,
u  \, \, = \, \,  \,
-\, {\frac {882000 \cdot \,(x-14)}{{x}^{3} +420\,x -5600}}, \quad \quad
v  \, \, = \, \,  \,
 \, {\frac { 49000 \cdot \, ({x}^{3} -21\,{x}^{2} -140) }{ ({x}^{3} +420\,x -5600)^{2}}} \cdot \, y. 
\end{eqnarray}
This change of variable (\ref{morphism}) actually transforms the  {\em elliptic curve} (\ref{genusone}) into
the {\em genus-two} curve (\ref{genustwo}).
{\em This provides a simple example of genus-two curve with split Jacobian through K3 surfaces}.

More generally, let us consider the Jacobian of a {\em genus-two} curve $\, C$. The Jacobian is simple
if it does not contain a proper abelian subvariety,
otherwise the Jacobian is reducible, or decomposable or ``split''. 
For this latter case, the only possibility for a genus-two curve is that
its Jacobian is {\em isogenous to a product} 
$\,\,  {\cal E} \times \,  {\cal E}'$ of {\em two elliptic curves}\footnote[2]{Along these lines, see
  also the concepts of Igusa-Clebsch invariants and Hilbert modular
  surfaces~\cite{Kumar,Kumar3,Mukamel,ElkiesKumar}. }.
Equivalently, there is a degree $\, n$ map
$\,\,  C \, \rightarrow \,  {\cal E} \, $ to some elliptic curves. Classically such
pairs\footnote[1 ]{One also has an anti-isometry Galois invariant  $\,\,  {\cal E}' \, \simeq \, {\cal E}$
  under Weil pairing. The decomposition corresponds to real multiplication by quadratic ring of discriminant $\, n^2$.}
$\, \, C, \,  {\cal E}\, $ arose in the {\em reduction of hyperelliptic integrals to elliptic ones}~\cite{Kumar}. The
$\, j$-invariants correspond, here,
to the {\em two elliptic subfields}: see~\cite{Kumar}.

\vskip .1cm

\subsection{Creative telescoping on rational functions in three variables associated with genus-two curves
with split Jacobians: a two-parameters example.}
\label{split}

Let us now consider the example with {\em two-parameters}, $a$ and $b$, given in section 4.5 page 12 of~\cite{Kumar}.
Let us substitute the rational parametrisation\footnote[5]{See also~\cite{Kuhn} section 6 page 48.}
\begin{eqnarray}
\label{morph}
\hspace{-0.98in}&& \quad  \quad \quad  \quad  \, 
u \, \, = \, \, \,  {\frac {{x}^{2}}{ {x}^{3} \, +a \, {x}^{2} \,  \, +b \, x \, +1 }},
\quad \quad \, \, \,  v \, \, = \, \, \,
{\frac { y  \cdot \, ( {x}^{3} \, -b \, x \, -2) }{ ({x}^{3} +a \, {x}^{2}\,\, +b \, x \, +1)^{2}}},
\end{eqnarray}
in the {\em elliptic curve} 
\begin{eqnarray}
\label{ellipticab} 
\hspace{-0.98in}&& \quad \quad \,
 R  \cdot \,{v}^{2} \, \, = \, \, \,\, \,
 R \cdot \, {u}^{3} \,\, \, +2 \cdot \, (a{b}^{2}\,  -6\,{a}^{2} \, +9\,b) \cdot \,  {u}^{2}
\, \, +  (12\,a -{b}^{2}) \cdot \,  u \,\,\,\, -4,              
\end{eqnarray}
where
\begin{eqnarray}
\label{whereR}
\hspace{-0.98in}&& \,  \quad  \quad  \quad \quad \quad  \quad \quad \quad 
R \, \, \, = \, \, \, \,
4 \cdot \, ({a}^{3}  \, +{b}^{3}) \,\, \,   \,  -{a}^{2}{b}^{2} \,\,  -18\,ab \,\,  \,  +27. 
\end{eqnarray}
This gives the {\em genus-two curve} $\,\, C_{a,\, b}(x, \, y) \, = \, \, 0 \,$ with:
\begin{eqnarray}
\label{genusTWO}
\hspace{-0.98in}&&  \, 
C_{a,\, b}(x, \, y) \, \,  \, = \, \,  \, \, \,
R \cdot \, y^2  \, \,  \, + (4\, x^{3} \, + b^{2} \, x^{2} \, +2\, b \, x \, +1)
\cdot \,  ({x}^{3} \, + a \, x^{2} \, +b \, x \, +1).          
\end{eqnarray}
The $ \, j$-invariant of the elliptic curve (\ref{ellipticab}) reads:
\begin{eqnarray}
\label{jinv}
  \hspace{-0.98in}&&   \quad \quad \quad
j \, \, = \, \, \,
 \,{\frac {16 \cdot \, ({a}^{2}{b}^{4}+12\,{b}^{5}-126\,a{b}^{3} +216 \,b{a}^{2} +405\,{b}^{2}-972\,a)^{3}}{
(4\,{a}^{3}+4\,{b}^{3} -{a}^{2}{b}^{2} \, -18\,ab\,
+27)^{2} \cdot \, (b -3)^{3} \cdot \, ({b}^{2}+3\,b+9)^{3}}}.
\end{eqnarray}
The Hauptmodul $ \, \, {\cal H} \, = \, 1728/j\, $ thus reads
\begin{eqnarray}
\label{Hauptab}
  \hspace{-0.98in}&&   \quad  \,\,
{\cal H} \, \,\, = \, \, \, 
\, {{ 108 \cdot \, (b-3)^{3}  \cdot \, (4\,{a}^{3}+4\,{b}^{3} -{a}^{2}{b}^{2}
-18\,ab+27)^{2} \cdot \,({b}^{2}+3\,b+9)^{3}}
\over {
(a^{2}{b}^{4} +12\,{b}^{5}-126\,a{b}^{3}+216\,b{a}^{2}+405\,{b}^{2} -972\,a)^3}}.
\end{eqnarray}
For $\,  \, \, b \, = \,  \, 3 \, +x \, \, $ this Hauptmodul (\ref{Hauptab}) reads
\begin{eqnarray}
\label{Hauptx}
  \hspace{-0.98in}&&   \quad  \quad  \quad    \quad    \quad   \quad    \quad    \quad   \quad  
{\cal H}_x \, \, = \, \, \,\, 
{\frac { 108  \cdot \,  {x}^{3} \cdot \, ({x}^{2}+9\,x+27)^{3} \cdot \, P_2^{2}}{ P_4^{3}}}, 
\end{eqnarray}
where:
\begin{eqnarray}
\label{HauptxP2P4}
  \hspace{-0.98in}&&  
P_2 \, = \, \,4\,{x}^{3}  \,  \, - (a-6)  \cdot \, (a+6)  \cdot \, {x}^{2}
\, \, -6 \cdot \, (a+6)  \cdot \, (a-3) \cdot \,  x \, \, + ( 4\,a \, +15)  \cdot (a \, -3)^{2},
\nonumber 
\\
\hspace{-0.98in}&&   
P_4 \, = \, \,\,
12\,{x}^{5} \,\, +({a}^{2}+180) \cdot \, {x}^{4} \,\, \, 
+6 \cdot \,   (2 \,  {a}^{2}   -21\,a  +180) \cdot \,  {x}^{3}
\nonumber \\
\hspace{-0.98in}&&   
 \quad  \quad \quad \quad  \quad \quad  \, \,
+27 \cdot  \, (2\, {a}^{2} \, -42\,a \, +135)  \cdot \, {x}^{2}
 \\
\hspace{-0.98in}&&   
 \quad  \quad \quad \quad  \quad \quad   \quad  \quad   \quad  \, \,
+162 \cdot \, (a \, -3)  \cdot \, (2\,a \, -15) \cdot \,  x
\,\,\,\,\, +729 \cdot \, (a \, -3)^{2}. \nonumber 
\end{eqnarray}

Let us consider the telescoper of the rational function of three variables
$\, x\, y/D_{a}(x,\, y, \, z)\, $ where the denominator
$\,\, D_{a}(x,\, y, \, z) \, \, $ is $\,\, C_{a,\, b}(x, \, y)\,$ given by  (\ref{genusTWO}),
but for $\, \,\,   b \, = \, \, 3 \, + \, \, x\, y\, z$: 
\begin{eqnarray}
  \label{ratgenus2}
 \hspace{-0.98in}&&  \quad 
D_{a}(x,\, y, \, z) \, \, = \, \, \, C_{a,\, 3\, +\, xyz}(x, \, y) 
\nonumber  \\
\hspace{-0.98in}&& \quad \quad \,  \, \,\, \, = \, \, \, \,
{x}^{6}{y}^{3}{z}^{3} \,+{x}^{7}{y}^{2}{z}^{2} \,+4\,{x}^{3}{y}^{5}{z}^{3}
+9\,{x}^{5}{y}^{2}{z}^{2}+6\,{x}^{6}yz+3\,{x}^{4}{y}^{2}{z}^{2}+36\,{y}^{4}{x}^{2}{z}^{2}
\nonumber  \\
  \hspace{-0.98in}&& \quad \quad  \quad\quad \quad \, \,\,
+6\,{x}^{5}yz+4\,{x}^{6}+27\,{x}^{4}yz
+9\,{x}^{5}+18\,{x}^{3}yz+108\,x{y}^{3}z+18\,{x}^{4}+3\,{x}^{2}yz
 \nonumber  \\
  \hspace{-0.98in}&& \quad \quad \quad  \quad \quad \quad \quad \quad \,
+32\,{x}^{3}+27\,{x}^{2}+135\,{y}^{2} \,\, +9\,x \,\, +1
\nonumber  \\
  \hspace{-0.98in}&& \quad \quad \, \,
+ \, ({x}^{6}{y}^{2}{z}^{2}+6\,{x}^{5}yz+2\,{x}^{4}yz+4\,{x}^{5}
-18\,x{y}^{3}z+9\,{x}^{4}+6\,{x}^{3}+{x}^{2}-54\,{y}^{2}) \cdot \,  a
\nonumber  \\
  \hspace{-0.98in}&& \quad \quad  \quad \quad \quad \quad \quad \quad \quad 
-{y}^{2} \cdot \, (xyz \, +3)^{2} \cdot \, {a}^{2}
\, \, \, \, + 4\, y^2 \cdot \, {a}^{3}.
\end{eqnarray}
This telescoper of the rational function
\begin{eqnarray}
  \label{ratgenus2Ra}
 \hspace{-0.98in}&&  \quad \quad \quad \quad \quad \quad \quad 
R_{a}(x,\, y, \, z) \, \, = \, \, \, 
\, {{ x\, y} \over { D_{a}(x,\, y, \, z)}},
\end{eqnarray}
is an {\em order-four} linear differential operator $\, L_4$ which
{\em is actually the direct-sum},
$\, L_4 \, = \, \, LCLM(L_2, \, M_2)\, = \, \, L_2 \oplus \, M_2$, 
of {\em two order-two} linear differential
operators,  having two pullbacked hypergeometric
solutions (see \ref{General}). One finds out that one of the two pullbacks
{\em precisely corresponds  to the Hauptmodul} $\, {\cal H}_x \, $
given by (\ref{Hauptx}).
This general case is detailed in \ref{General}.

Let us consider the $\, a \, = \, 3\, $ subcase\footnote[1]{The discriminant in
 $\, b\, $ of $\, \, 4\,{a}^{3}+4\,{b}^{3} -{a}^{2}{b}^{2} \, -18\,ab\,  +27\,\,  $
 reads: $\, (a-3)^3 \cdot \, (a^2 \, +3\,a \, +9)^3\, $, consequently
the exact expressions are simpler at $\, a=3$. }. For $\, a\, = \, 3$,
the Hauptmodul $ \, {\cal H} \, = \, 1728/j$,  corresponding to
the $\, j$-invariant (\ref{jinv}), reads:
\begin{eqnarray}
\label{Haupt}
  \hspace{-0.98in}&&   \quad \quad \quad \quad \quad \quad
{\cal H} \, \,\,\, = \, \, \,\,\,
 {\frac { 4 \cdot  \, (b \, -3)  \cdot \, (4\,b+15)^{2} \cdot \,
({b}^{2} \, +3\,b \, +9)^{3} }{ (b \, +6)^{3} \cdot \, (4\,{b}^{2}+3\,b-18)^{3}}}.
\end{eqnarray}
This Hauptmodul becomes for $\,\,\, b \, = \,\, 3 \,\, +x$
\begin{eqnarray}
\label{Hauptxx}
  \hspace{-0.98in}&&   \quad \quad \quad \quad \quad \quad \quad
{\cal H} \, \,\,\, = \, \, \,\,\,
\, {\frac { 4 \cdot \, x \cdot \,  (27+4\,x)^{2} \cdot  \,
({x}^{2}+9\,x+27)^{3}}{ (9+x)^{3} \cdot \, (4\,{x}^{2}+27\,x+27)^{3}}}. 
\end{eqnarray}

The telescoper of the rational function
(\ref{ratgenus2Ra}) with  $\,\, D_{a}(x,\, y, \, z)$  given by (\ref{ratgenus2})
for $\, a \, = \, 3$,  is an {\em order-four}
linear differential operator which is the direct-sum of two
order-two linear differential
operators $\, L_4 \, = \, \, LCLM(L_2, \, M_2)\, = \, \, L_2 \oplus \, M_2$, these
{\em two order-two} linear differential  operators having  the pullbacked
hypergeometric solutions
\begin{eqnarray}
\label{Hauptx1}
\hspace{-0.98in}&&   \quad \quad  \quad \quad \quad \quad 
 (27\, + \, 4 \, x)^{-1/2}     \cdot \,  x^{-5/4}  \cdot \,
  _2F_1\Bigl([{{1} \over {12}}, \,{{5} \over {12}}], \, [1], \, \,
 1 \, + \, {{ 27} \over { 4 \, x}}  \Bigr), \quad
\end{eqnarray}
for $\, L_2$, and 
\begin{eqnarray}
\label{Hauptx2}
\hspace{-0.98in}&&   \quad \quad \quad 
{{3+x} \over { (9+x)^{1/4} \cdot \, (4\,{x}^{2}+27\,x+27)^{1/4} \cdot \, {x}^{3/2} \cdot (27+4\,x)^{1/2} }} 
\nonumber \\
\hspace{-0.98in}&&   \quad \quad \quad \quad \quad  \, \,
 \times \, \,
 _2F_1\Bigl([{{1} \over {12}}, \,{{5} \over {12}}], \, [1], \, \,
 \, {\frac {  4 \cdot \, x \cdot \, (27 \,+4\,x)^{2} \cdot \, ({x}^{2} \, +9\,x \, +27)^{3}}{
 (9+x)^{3} \cdot \, (4\,{x}^{2} \, +27\,x \, +27)^{3}}}\Bigr),  
\end{eqnarray}
for $\, M_2$, {\em where we see clearly that the Hauptmodul in} (\ref{Hauptx2})
{\em is precisely the  Hauptmodul} (\ref{Hauptxx}). 
The {\em Jacobian of the genus-two curve is a split  Jacobian corresponding to the product}
$\, {\cal E}_1 \, \times {\cal E}_2$
 {\em of two elliptic curves}, the $\, j$-invariant of the second elliptic curve being  (\ref{jinv}),
when the $\, j$-invariant of the first elliptic curve  reads
\begin{eqnarray}
\label{j2}
\hspace{-0.98in}&&   \quad \quad \quad \quad \quad \quad \quad  \quad \quad \quad \quad
j_1 \, \, = \, \, \, {{ 6912 \, x } \over { 27 \, \,+ 4 \, x}},              
\end{eqnarray}
corresponding to the Hauptmodul
$\,\, \, 1728/j_1  \, = \, \, \,  1 \, + \, {{ 27} \over { 4 \, x}}\,\,  $ in (\ref{Hauptx1}). 
This second invariant is, as it should, {\em exactly the 
$\, j$-invariant of the second elliptic curve} $\, {\cal E}'$, given page 48 in~\cite{Kuhn}:
\begin{eqnarray}
\label{j2Kuhn}
\hspace{-0.98in}&&   \quad \quad \quad \quad \quad \quad \quad
j({\cal E}') \, \, = \, \, \,
   \,{\frac { 256 \cdot \, (3\,b  \,  -{a}^{2})^{3}}{
 4\,{a}^{3}c \, \, -{a}^{2}{b}^{2} \, \, -18\,abc \, \, +4\,{b}^{3} \, +27\,{c}^{2}}},            
\end{eqnarray}
for $\,\,  c\, = \, \, 1$, $\, a \, = \, \, 3  \, $ and $\,\,  b \, = \, \, 3 \, + \, x$.

\vskip.2cm

This subcase $\, a\, = \, 3$ is very special. The general case of arbitrary value $\, a$ is
sketched in  \ref{General}.

\vskip.2cm

For some general facts\footnote[2]{Explicit calculations require to use
various tools in Magma: AnalyticJacobian, EndomorphismRing, ToAnalyticJacobian, FromAnalyticJacobian, ...} 
about algebraic curves,
their Jacobians and algebraic correspondences, see~\cite{correspondences,ShimuraCorrespondences,ModularCor} and
page 301 in~\cite{Cambridge}: the multiplier equation
is seen to contain the modular equation as a particular case. In~\cite{Klein} Felix Klein
{\em also defined the associated idea of modular correspondence}\footnote[5]{See also F. Klein and R. Fricke in~\cite{Klein2} }.  

\vskip .2cm

\subsection{Creative telescoping on rational functions of three variables associated with genus-two curves
with split Jacobians: a simple example}
\label{splitgood}

Let us consider another simpler example of  {\em genus-two} curve with pullbacked
$\, _2F_1$ solution (not product of pullbacked $\, _2F_1$) of the telescoper.

\vskip .1cm

Let us consider the {\em genus-two algebraic curve} $\,\, C_p(x, \, y) \, = \, \, 0 \,\, $
given in lemma 7 of~\cite{Shaska2} (see also ~\cite{Shaska,Shaska3})
\begin{eqnarray}
\label{splitx5}
\hspace{-0.98in}&&   \quad  \quad  \quad \quad \quad \quad \quad \quad
 C_p(x, \, y) \, \, = \,\,  \,  \,\,  x^5 \, +x^3 \,\, \, + p \cdot \, x \, \, \,   -y^2,       
\end{eqnarray}
where $\, C_p(x, \, y)\, $ is  given in lemma 7 of~\cite{Shaska2}.
Let us introduce the rational function  $\, x\, y/D(x, \, y, \, z)\, $  where the denominator
$\,  D(x, \, y, \, z) \, \, $ is given by:
\begin{eqnarray}
\label{ratsplit}
\hspace{-0.98in}&&   \quad  \quad \,  \quad \quad \quad
D(x, \, y, \, z) \, \, = \, \, \, \,  C_{p\, = \, xyz}(x, \, y) \, \, \, = \,\, \,  \, \,
x^5 \, +x^3 \, \, + x^2\, y\, z \, \,  \,   -y^2.       
\end{eqnarray}
The telescoper of this rational function is an order-two linear differential operator
which has the two hypergeometric solutions
\begin{eqnarray}
\label{solratsplit1}
  \hspace{-0.98in}&&   \quad \quad \quad  \quad \quad  \quad \quad \quad \quad \quad
   x^{-1/4} \cdot \,  _2F_1\Bigl([{{1} \over {8}}, \, {{5} \over {8}}], \, [{{3} \over {4}}], \, \, 4\, x\Bigr)  
\end{eqnarray}
which is a Puiseux series at $\, \, x \, = \, \, 0 \, \, $ and: 
\begin{eqnarray}
\label{solratsplit2}
  \hspace{-0.98in}&&   \quad \quad \quad \quad  \quad \quad  \quad \quad \, \quad
 x^{-1/4} \cdot \,  _2F_1\Bigl([{{1} \over {8}}, \, {{5} \over {8}}], \, [1], \,\, 1 \, -4\, x\Bigr).
\end{eqnarray}
These two  hypergeometric solutions
can be rewritten as\footnote[1]{The fact that $\, _2F_1\Bigl([{{1} \over {8}}, \, {{5} \over {8}}], \, [1], \, z\Bigr)$
  can be rewritten as $\, _2F_1\Bigl([{{1} \over {12}}, \, {{5} \over {12}}], \, [1], \, H(z)\Bigr)$ where the Hauptmodul
$\, H(z)$ is solution of a quadratic equation is given in equation (H.14) of Appendix H of~\cite{malala}.} 
\begin{eqnarray}
\label{solratsplit3}
  \hspace{-0.98in}&&   \quad \quad \quad \quad \quad \quad \quad \quad \quad 
 {\cal A}(x) \cdot \,  _2F_1\Bigl([{{1} \over {12}}, \, {{5} \over {12}}], \, [1], \, {{1728}  \over {J}} \Bigr),
\end{eqnarray}
where the $\,j$-invariant $\, J$, in the Hauptmodul $\,\, 1728/J \,$ in (\ref{solratsplit3}), 
{\em corresponds exactly to the degree-two elliptic subfields}
\begin{eqnarray}
\label{solratsplit2subfield}
  \hspace{-0.98in}&&   \quad \quad  
{J}^{2} \,  \, \, \, -128 \cdot \,{\frac { (2000\,{x}^{2} \, +1440\,x \,+27) }{ (1 \, -4\,x)^{2}}}  \cdot \,  J
\, \, \,\,  -4096 \cdot \,{\frac { (100\,x-9)^{3}}{(1 \, -4\,x)^{3}}}
 \, \, \,  = \,  \, \,  \, \, 0, 
\end{eqnarray}
given in the first equation of page 6 of~\cite{Shaska2}.

Of course, if we change $\, p$ into $\, p \, \rightarrow  \, (1\, -p)/4 \, $
in subsection \ref{split}, the telescoper
of  the rational function  $\, \, x\, y/D(x, \, y, \, z)\, $ where the denominator
$\,  D(x, \, y, \, z)\, \,$ is given by:
\begin{eqnarray}
\label{ratsplitter}
\hspace{-0.98in}&&   \quad  \, \, \,   \, \, \,   \, 
D(x, \, y, \, z) \, \, = \, \, \, \,  C_{p\, = \, xyz}(x, \, y) \, \, \, = \,\, \, \, \, \,
x^5 \, +x^3 \, \,\,  +\Bigl({{ 1\, -x\, y\, z} \over {4}}\Bigr) \cdot x \, \,  \,  \,   -y^2,       
\end{eqnarray}
is the  order-two linear differential operator corresponding
to the $\, x \, \rightarrow  \, (1\, -x)/4$ pullback of the previous one. It has the two
hypergeometric solutions
\begin{eqnarray}
\label{solratsplitter}
  \hspace{-0.98in}&&   \quad \quad  \quad   \quad \quad  \quad \quad \quad \quad
 (1\, -x)^{-1/4} \cdot \,
_2F_1\Bigl([{{1} \over {8}}, \, {{5} \over {8}}], \, [1], \, \,  \,\, \, x\Bigr),
\end{eqnarray}
and:
\begin{eqnarray}
\label{solratsplit1ter}
 \hspace{-0.98in}&&   \quad \quad \quad \quad   \quad \quad \quad \quad \quad
 (1\, -x)^{-1/4} \cdot \,
 _2F_1\Bigl([{{1} \over {8}}, \, {{5} \over {8}}], \, [{{3} \over {4}}], \, \,  \, \, 1\, -\, x\Bigr).  
\end{eqnarray}

\vskip .1cm

{\bf Remark:} In contrast with the previous example of subsection \ref{split}
where we had two $\, j$-invariants corresponding to the {\em two order-two} linear differential
operators $\, L_2$ and $\, M_2$ of the direct-sum decomposition
of the order-four telescoper, we have, here,  {\em just one  order-two}  telescoper,
which is enough to ``encapsulate'' two  $\, j$-invariants (\ref{solratsplit2subfield}). 
One  order-two  linear differential operator is enough because the two 
$\, j$-invariants are Galois-conjugate (see (\ref{solratsplit2subfield})). 

\vskip .2cm

\subsection{Creative telescoping on rational functions of three variables associated with genus-two curves
with split Jacobians: another simple example}
\label{splitgoodbis}

Another similar example of {\em genus-two algebraic curve} $\, C_p(x, \, y) \, = \, \, 0 \, $
given in equation (5) of lemma 4 of~\cite{Shaska2} 
\begin{eqnarray}
\label{splitbis}
\hspace{-0.98in}&&   \quad  \quad  \quad \quad \quad \quad \quad \quad
 C_p(x, \, y) \, \,\, = \,\,  \,\, \, \,  x^6 \, +x^3 \, \, + p \, \, \, \,  -y^2,       
\end{eqnarray}
with a split  Jacobian, yields
an order-two telescoper for the corresponding rational function,
with pullbacked hypergeometric solutions, where, again, the $\,j$-invariant $\, J$,
in the Hauptmodul $\,\, 1728/J \,$  {\em corresponds exactly to the degree-two elliptic subfields}
of the split Jacobian of the {\em genus-two} curve. 
More details are sketched  in \ref{splitgoodbisapp}. 

\vskip .3cm 

\subsection{Creative telescoping on rational functions of three variables associated with genus-two curves
with split Jacobians: another example}
\label{splitgoodbister}

Let us now recall the paper~\cite{Diarra} by
K. Diarra. Similarly to subsection \ref{split},
we consider the {\em one-parameter} $a$-example given in section 6 page 52 of~\cite{Diarra}.
Let us substitute the rational parametrisation
\begin{eqnarray}
\label{morph2}
\hspace{-0.98in}&& \quad  \quad   \quad \quad \quad \quad  \quad  \quad  \, 
u \, \, = \, \, \,  -{\frac {{x}^{2}+a}{{x}^{2}-1}}, 
\quad \quad \,  \, \quad \, \,
v \, \, = \, \, \, {\frac {y}{ \left( {x}^{2}-1 \right) ^{2}}}, 
\end{eqnarray}
in the {\em elliptic curve} $\, \, P(u,  \, v) \, \, = \, \, \, 0$
\begin{eqnarray}
\label{ellipticabbis} 
  \hspace{-0.98in}&& \quad \quad  \quad \quad  \quad \quad  \quad  \quad \quad \,
{u}^{3} \, +1 \,  \,  \,  - \, \,  {v}^{2}
 \,  \, \, = \, \, \,\, \, 0,
\end{eqnarray}
which $\, j$-invariant is $\, 0$ (and thus the Hauptmodul is $\, \infty$).
This gives the {\em genus-two curve} $\,\, C_{a}(x, \, y) \, = \, \, 0 \,$ with:
\begin{eqnarray}
\label{genusTWObis2}
\hspace{-0.98in}&&  \quad \quad \quad 
C_{a}(x, \, y) \, \,  \, = \, \,  \, 
 \\
\hspace{-0.98in}&&  \, \quad \quad \quad \quad \quad \quad 
{y}^{2} \,  \, \, +  (a \, +1)  \cdot \, ({x}^{2} \, -1)  \cdot \,
\Bigl(3\,{x}^{4} \, +3 \cdot \, (a \, -1) \cdot \, {x}^{2} \, +{a}^{2}  -a \, +1\Bigr).
 \nonumber 
\end{eqnarray}
Let us denote $\, {\cal C}(x, \, y, \, z)$ the previous polynomial (\ref{genusTWObis2})
where the parameter $\, a$ becomes the product $\, a \, = \, \,  \, x\, y\, z$.
The telescoper of the rational function 
\begin{eqnarray}
\label{telescgenusTWObis}
  \hspace{-0.98in}&&  \, \quad \quad    \quad \quad  \quad \quad  \quad \quad  \quad
 R(x, \, y, \, z) \, \, = \, \, \, \, {{ x\, y} \over {  {\cal C}(x, \, y, \, z)}}, 
\end{eqnarray}          
is an {\em order-two} linear differential operator
\begin{eqnarray}
\label{telescgenusTWObisordertwo1}
\hspace{-0.98in}&&  \, \quad  \quad    \,  
L_2 \, \, = \, \,  \,  \, \,  3 \cdot \, (3\,x \, -1) \,\,\,
+4 \cdot \, (4\, x^2 \, -x \, +1) \cdot \, D_x\,  \, +4  \cdot \,(1 \, +x^3)  \cdot \,D_x^2, 
\end{eqnarray}         
which has the following $\, _2F_1$ hypergeometric solutions
\begin{eqnarray}
\label{solgenusTWObisordertwo1}
  \hspace{-0.98in}&&  \,   \quad  \quad \quad  \quad
(1 \, +x)^{-1}  \cdot \,   (1 \,-2\,x)^{-1/2}  \cdot \,
 _2F_1\Bigl([{{1} \over {4}}, \, {{3} \over {4}}], \, [1], \, -{{3} \over {(1\, -2\, x)^2}} \Bigr), 
\end{eqnarray}
or:
\begin{eqnarray}
\label{solgenusTWObisordertwo12}
  \hspace{-0.98in}&&  \, \quad    \quad \quad  \quad
(1 \, +x)^{-1}  \cdot \,   (1 \,-2\,x)^{-1/2}  \cdot \,
 _2F_1\Bigl([{{1} \over {4}}, \, {{3} \over {4}}], \, [1],
      \, {{4  \cdot \, ({x}^{2}-x+1) } \over {  (1 \,-2\,x)^2  }} \Bigr), 
\end{eqnarray}
and:
\begin{eqnarray}
\label{solgenusTWObisordertwo21}
  \hspace{-0.98in}&&  \, \quad \quad \quad  \quad
 (1 \, +x)^{-5/4}  \cdot \,   (x\, -2)^{-1/4} 
\nonumber \\
\hspace{-0.98in}&& \quad \quad \quad \quad \quad  \quad \quad  \quad 
\times \, \, \,
 _2F_1\Bigl([{{1} \over {12}}, \, {{5} \over {12}}], \, [1], \,
  -\, {\frac { 81 \cdot \, ({x}^{2} \, -x \, +1)^{2}}{ 4 \cdot \, (x\, -2)^{3}  \cdot \, (x\, +1 )^{3}}} 
 \, \Bigr). 
\end{eqnarray}
Again the telescoper is {\em not an  order-four} operator 
but an {\em order-two operator}. This is a consequence of the fact that, among the
two $\, j$-invariants of the split Jacobian, one is trivial ($j \, = \, 0$).
Note that the Hauptmodul in (\ref{solgenusTWObisordertwo21}) is simply related to
the  Hauptmodul $\, \, 1728\, z/(z\, +16)^3 \, $ in~\cite{Maier1} (see $\, N= \, 2$
in Table 4 page 11 in~\cite{Maier1}):
\begin{eqnarray}
\label{relMaier}
 \hspace{-0.98in}&&  \,  \,   \, 
  -\, {\frac { 81 \cdot \, ({x}^{2} \, -x \, +1)^{2}}{ 4 \cdot \, (x\, -2)^{3}  \cdot \, (x\, +1 )^{3}}}
  \, \, = \, \, \, {{ 1728 \, z} \over {(z\, +16)^3 }} \quad   \, \quad \hbox{when:}  \quad  \quad
  z \, \, = \,\, \, {{ -\, 48 } \over {{x}^{2} \, -x \, +1 }}. 
\end{eqnarray}      
Consequently we can also write (see $\, N= \, 2$ in Table 5 page 12 in~\cite{Maier1})
the solution of the order-two telescoper (\ref{telescgenusTWObisordertwo1})
in terms of the alternative Hauptmodul:
\begin{eqnarray}
\label{relMaier}
 \hspace{-0.98in}&&  \,   \quad
  \,  \, {\frac { 972 \cdot \, ({x}^{2} \, -x \, +1)}{ (16 \,x^2 \, -16\,x \, +13  )^{3}}}
  \, \, = \, \, \, {{ 1728 \, z^2} \over {(z\, +256)^3 }} \quad  \,  \quad \hbox{with:}  \quad  \quad
  z \, \, = \,\, \, {{ -\, 48 } \over {{x}^{2} \, -x \, +1 }}. 
\end{eqnarray}   
This alternative writing of the solution reads:
\begin{eqnarray}
\label{solgenusTWObisordertwo21altern}
  \hspace{-0.98in}&&  \, \quad \quad \, 
 {{(16\,x^2 \, -16\,x \, +13 )^{-1/4}} \over {1\, + \, x}}  \, \cdot \,
 _2F_1\Bigl([{{1} \over {12}}, \, {{5} \over {12}}], \, [1], \,
 {\frac { 972 \cdot \, ({x}^{2} \, -x \, +1)}{ (16\,x^2 \, -16\,x \, +13  )^{3}}}
 \, \Bigr). 
\end{eqnarray}   

\vskip .3cm
\vskip .2cm 

\section{Rational functions with tri-quadratic denominator and $\, N$-quadratic denominator.}
\label{triquadratic}

We try to find telescopers of rational functions  corresponding to (factors of) linear differential 
operators of small orders, for instance order-two linear differential operators 
with pullbacked $\, _2F_1$ hypergeometric functions, classical modular forms,
or their modular generalisations (order-four Calabi-Yau linear
differential operators~\cite{Straten}, etc ...). As we saw in the previous sections,
this corresponds to the fact that the denominator  of these rational functions
is associated with an elliptic curve, or products of elliptic curves, with
K3 surfaces or  with threefold Calabi-Yau manifolds corresponding to algebraic varieties
with foliations in elliptic curves\footnote[4]{Even if K3 surfaces, or threefold Calabi-Yau manifolds,
are {\em not} abelian varieties, the Weierstrass-Legendre forms introduced in the previous section,
amounts to saying that K3 surfaces can be ``essentially viewed'' (as far as creative telescoping is concerned)
as  foliation in two elliptic curves, and threefold Calabi-Yau manifolds as foliation in three elliptic curves.}.
Since this paper tries to reduce the {\em differential algebra} creative telescoping calculations to
{\em effective algebraic geometry} calculations\footnote[3]{One has birational automorphisms in projective spaces~\cite{Bedford,Angles},
  but since this paper is dedicated to (efficient) formal calculations we work exclusively in affine coordinates (see
  for instance (\ref{involutionbirat}), (\ref{triquadrarewritt1}), (\ref{involutionbiratbis}) below). For algebraic geometers
  an ellitic curve is a smooth complete genus 1 curve with a choice of a base point. 
  Here our elliptic curves are, in fact, an affine piece of a genus 1 curve with no base point, but this does not really
  matter because the $\, j$-invariant which is all we care about in this kind of creative telescoping calculations, is determined by
that much information.} 
and structures, we want to
focus on rational functions  with denominators that correspond to {\em selected}
algebraic varieties~\cite{Quasi,JMP}, beyond algebraic varieties corresponding to products
of elliptic curves or foliations in elliptic  curves\footnote[8]{K3 surfaces,
threefold Calabi-Yau manifolds, higher curves with  split Jacobian corresponding to products
of elliptic curves, ...},  namely algebraic varieties
{\em with an infinite number of birational automorphisms}\footnote[1]{The best explicit illustration
of this situation emerges in integrable models~\cite{Quasi,JMP,BMV1,BMV2}}. This
{\em infinite number of birational symmetries}, excludes algebraic varieties
of the ``general type'' (with  {\em finite} numbers\footnote[2]{There are even precise bounds for the number
  of automorphisms. The upper bound is $\, 84 \, (g-1)$ for curves of genus $\, g$ and these
  bounds have been extensively studied in higher dimensions~\cite{Corti,Szabo,Hacon}.} of birational symmetries).
For algebraic surfaces, this amounts to discarding the  surfaces of the ``general type''
which have Kodaira dimension 2, focusing on Kodaira dimension one (elliptic surfaces),
or  Kodaira dimension zero (abelian surfaces, hyperelliptic surfaces, K3 surfaces, Enriques surfaces),
or even Kodaira dimension $ \, -\infty$ (ruled surfaces, rational surfaces).  

In contrast with algebraic {\em curves} where one can easily, and very efficiently, calculate
the genus of the curves to discard the algebraic curves of higher genus and, in the case of genus-one,
obtain the $\, j$-invariant using formal calculations\footnote[5]{Use with(algcurves) in Maple
  and the command ``genus'' and ``j$\_$invariant''. }, it is, in practice, quite difficult to see
for higher dimensional algebraic varieties, that the  algebraic variety is not
of the ``general type'',  because it has an {\em infinite number of birational symmetries}.
For  these ``selected cases'' we are interested in, calculating the generalisation
of the $\, j$-invariant (Igusa-Shiode invariants, etc ...) is quite  hard.

Along this line we want to underline that there exists a remarkable set of algebraic surfaces,
namely the algebraic surfaces corresponding to tri-quadratic equations: 
\begin{eqnarray}
\label{triquadra}
\hspace{-0.98in}&& \, \,     \quad  \quad   \quad  \quad  \quad  \quad 
\sum_{m=0,1,2} \, \sum_{n=0,1,2} \, \sum_{l=0,1,2} \, a_{m,n,l} \cdot \, x^m \, y^n \, z^l
\, \, \, \,  = \, \,  \, \, \,  0, 
\end{eqnarray}
depending on $\,\, 27 = \, 3^3 \,$ parameters $ \,  a_{m,n,l}$.
More generally, one can introduce algebraic varieties corresponding to $\, N$-quadratic equations: 
\begin{eqnarray}
\label{Nquadra}
\hspace{-0.98in}&& \, \, \,  \, \,    \, \,  \, \, 
\sum_{m_1=0,1,2} \, \sum_{m_2=0,1,2} \, \cdots \,  \, \sum_{m_N=0,1,2} \,
a_{m_1, \, m_2, \cdots, \, m_N} \cdot \, x_1^{m_1} \, x_2^{m_2} \, \cdots \,  x_N^{m_N}
\, \, \,  = \, \,  \, \, \,  0. 
\end{eqnarray}
With these tri-quadratic (\ref{triquadra}), or $\, N$-quadratic  (\ref{Nquadra}) equations,
we will see, in subsection \ref{subtriquadratic2}, that we have {\em automatically}
(selected) algebraic varieties that are not
of the ``general type''  having an {\em infinite number of birational symmetries}, which is precisely
our requirement for the denominator of rational functions with remarkable telescopers\footnote[1]{
Telescopers with factors of small enough order, possibly yielding classical modular forms, Calabi-Yau operators, ...
Rational functions with denominators of the ``general type'' will yield
telescopers of very large orders.}. 

Let us first, as a warm-up, consider, in the next subsection, a remarkable example
of tri-quadratic (\ref{triquadra}), where the underlying foliation in elliptic curves
is crystal clear.

\vskip .1cm 

\subsection{Rational functions with tri-Quadratic denominator simply corresponding to elliptic curves.}
\label{subtriquadratic}

Let us first recall the tri-quadratic equation in three variables
$\, x$, $\, y$ and $\, z$:
\begin{eqnarray}
\label{surfacetriquadr}
  \hspace{-0.98in}&&     \quad   \, \,\,    \quad  \quad    \, \,   \, \,
{x}^{2}{y}^{2}{z}^{2} \, \,\, \,\, -2 \cdot \,M \cdot \, x y z \cdot \, (x \, +y \, +z)
\, \, \, \, \, +4 \cdot \,M \cdot \, (M \, +1) \cdot \, xyz
\nonumber \\
\hspace{-0.98in}&&  \quad     \,     \quad  \quad  \quad  \quad   \, \, \, \,
+{M}^{2} \cdot \,  ({x}^{2}+{y}^{2}+{z}^{2})  \,  \,  \, \,
-2\,{M}^{2} \cdot \, (xy \, +xz \, +yz)
\, \, \, = \, \,\, \, 0,                
\end{eqnarray}
already introduced in Appendix C of~\cite{Youssef}.  
This algebraic surface, symmetric in $\, x$, $\, y$ and $\, z$, can be seen
for $\, z$ (resp.  $\, x$ or $\, y$)
fixed, as an {\em elliptic curve} which $\, j$-invariant is {\em independent of $\, z$} and reads 
\begin{eqnarray}
\label{jsurfacetriquadr}
\hspace{-0.98in}&&  \quad     \quad     \quad    \quad     \quad   \quad    \quad  \quad    \quad  \, \,
j \, \, = \, \, \,
256 \cdot \,{\frac { ({M}^{2}-M+1)^{3}}{{M}^{2} \cdot \, (M-1)^{2}}},
\end{eqnarray}
the corresponding Hauptmodul reading:
\begin{eqnarray}
\label{Hauptsurfacetriquadr0}
\hspace{-0.98in}&&  \quad     \quad     \quad    \quad   \quad      \quad  \quad    \quad  \quad   \, \,
{\cal H} \, \, = \, \, \,   {\frac {27 \cdot \,{M}^{2} \cdot \,  (M-1)^{2}}{4\, \cdot \, ({M}^{2}-M+1)^{3}}}.
\end{eqnarray}
This corresponds to the fact that this algebraic surface (\ref{surfacetriquadr})
can be seen as a product of two times the same elliptic curve with the $\, j$-invariant (\ref{jsurfacetriquadr})
or the Hauptmodul (\ref{Hauptsurfacetriquadr0}). This is a consequence of the fact that, introducing 
$\, x \, = \, sn(u)^2$,  $\, y \, = \, sn(v)^2$ 
and  $\, z \, = \, sn(u \, +v)^2$, and $\, M \, = \, \, 1/k^2$, 
this algebraic surface (\ref{surfacetriquadr}) corresponds to
the well-known formula for the
{\em addition on elliptic sine}\footnote[5]{See equation (C.3) in Appendix C of~\cite{Youssef}.}:
\begin{eqnarray}
\label{additionellipticsinus}
\hspace{-0.95in}&& \quad  \quad  \quad \quad  \quad  
sn(u \,+ \,  v) \,\,\, = \, \, \, \, {{ 
sn(u) \, cn(v) \, dn(v) \,\,  + \, \, sn(v) \, cn(u) \, dn(u) 
} \over { 1 \, \, - k^2 \, sn(u)^2 \,  sn(v)^2
}}.
\end{eqnarray}
For $\, M \, = \, \, x \, y\, z\, w$, the LHS of the tri-quadratic equation (\ref{surfacetriquadr})
yields a polynomial of {\em four variables }$\, x$,  $\, y$, $\, z$ and $\, w$, that we denote
$\,\,  T (x, \, y, \, z, \, w)$:
\begin{eqnarray}
\label{surfacetriquadr2bis}
  \hspace{-0.98in}&&
 \,  \quad     \quad         \quad    T (x, \, y, \, z, \, w)
\, \,\,  = \, \,\,  \,
 \\
\hspace{-0.98in}&&  \quad       \quad    \quad    \quad  \quad  \,\,  \, \, \,        
{x}^{2}{y}^{2}{z}^{2} \, \, \,\, -2  \,  \cdot \, x^2 y^2 z^2 \, w \cdot \, (x \, +y \, +z)
 \,  \, \,\, +4  \cdot \, (x y z w \, +1) \cdot \, x^2 y^2 z^2 \, w 
\nonumber \\
\hspace{-0.98in}&&  \quad     \quad  \quad     \quad     \quad  \,\,  \,  \quad  \quad   \, \, \, \,
 + x^2 y^2 z^2  w^2 \cdot \,  ({x}^{2}+{y}^{2}+{z}^{2})
\,  \, \,  \, -2\, x^2 y^2 z^2  w^2 \cdot \, (xy \, +xz \, +yz).\nonumber 
\end{eqnarray}
The telescoper of the rational function in {\em four variables} $\, x, \, y, \, z$ and $\, w$,
\begin{eqnarray}
\label{Ratsurfacetriquadr2bis}
  \hspace{-0.98in}&&  \quad   \quad     \, \, \quad
 \quad \quad    \quad  \quad    \quad  \quad    \quad   \, \,
     {{x \,y \,z } \over {  T (x, \, y, \, z, \, w)  }}, 
\end{eqnarray}
is an order-three (self-adjoint) linear differential operator which is
the {\em symmetric square}
of the order-two linear differential operator 
having the following pullbacked $\, _2F_1$ hypergeometric solution: 
\begin{eqnarray}
\label{solRatsurfacetriquadr2bis1}
 \hspace{-0.98in}&&  \quad   \quad   \quad  \quad     \quad  \quad  
 x^{-1/2}  \cdot \,({x}^{2}-x+1)^{-1/4} \cdot \,
 \nonumber \\
 \hspace{-0.98in}&&  \quad  \quad   \quad   \quad    \quad  \quad     \quad  \quad   \quad  \quad  
\times \, \,     _2F_1\Bigl([{{1} \over {12}}, \, \, {{5} \over {12}}], \, [1], \, \,
 {\frac {27 \cdot \,{x}^{2} \cdot \,  (x-1)^{2}}{4\, \cdot \, ({x}^{2}-x+1)^{3}}} \Bigr).                      
\end{eqnarray}                  
As it should the Hauptmodul in (\ref{solRatsurfacetriquadr2bis1}) is the same as
the Hauptmodul (\ref{Hauptsurfacetriquadr0}).  
The algebraic surface (\ref{surfacetriquadr}) can be seen
as the product of {\em two times the same elliptic curve} with the
 Hauptmodul (\ref{Hauptsurfacetriquadr0}). As expected the
 solution of the order-three telescoper is the {\em square} of the pullbacked $\, _2F_1$ hypergeometric
 function (\ref{solRatsurfacetriquadr2bis1}) with that  Hauptmodul.

More generally, we can also introduce 
 the tri-quadratic equation of three variables
$\, x$, $\, y$ and $\, z$ and two parameters $\, M$ and $\, N$:
\begin{eqnarray}
\label{surfacetriquadr2}
  \hspace{-0.98in}&&  \quad   \quad     \quad    \quad   \, \,
{x}^{2}{y}^{2}{z}^{2} \, \,\,\, \, -2 \,M \cdot \, x y z \cdot \, (x \, +y \, +z)
\, \, \, \, + N \cdot \, xyz
  \\
\hspace{-0.98in}&&  \quad \quad      \quad  \quad    \quad  \quad    \quad   \, \,  \, \,
+{M}^{2} \cdot \,  ({x}^{2}+{y}^{2}+{z}^{2})  \,  \,  \,   \, -2\,{M}^{2} \cdot \, (xy \, +xz \, +yz)
\, \,\,  = \, \, \, \, 0.    \nonumber              
\end{eqnarray}
This surface, symmetric in $\, x$, $\, y$ and $\, z$,
can be seen for $\, z$ (resp.  $\, x$ or $\, y$)
fixed as an elliptic curve which $\, j$-invariant is, again, {\em independent of $\, z$} and reads 
\begin{eqnarray}
\label{jsurfacetriquadr2}
\hspace{-0.98in}&&  \quad     \quad    \quad  \quad    \quad    \quad  \quad    \quad  \quad    \quad   \, \,
j \, \, = \, \, \,
{\frac { (48 \,{M}^{3}-{N}^{2})^{3}}{ {M}^{6} \cdot \, (64\,{M}^{3}-{N}^{2}) }}.
\end{eqnarray}
the corresponding Hauptmodul reading:
\begin{eqnarray}
\label{Hauptsurfacetriquadr}
\hspace{-0.98in}&&  \quad     \quad    \quad  \quad   \quad     \quad  \quad    \quad  \quad    \quad   \, \,
{\cal H} \, \, = \, \, \, \,
{\frac { 1728  \cdot  \, {M}^{6} \cdot \, (64\,{M}^{3}-{N}^{2}) }{ (48\,{M}^{3}-{N}^{2})^{3}}}.
\end{eqnarray}
Let us consider the following change of variables $\, M\, = \, \, m^2 \, $ and
$\, \,  N \, = \, \, 8 \cdot \, m^3 \, + p$ in  (\ref{surfacetriquadr2}). 
For $\, p \, = \, \, x \, y\, z\, w$, the LHS of the tri-quadratic equation (\ref{surfacetriquadr2})
yields a polynomial in {\em four variables} $\, x$,  $\, y$, $\, z$ and $\, w$, that we denote
$\,  {\cal T}_{m} (x, \, y, \, z, \, w)$:
\begin{eqnarray}
\label{surfacetriquadr2bis1}
  \hspace{-0.98in}&&  \,  \quad \quad   \quad   \quad    
{\cal T}_{m} (x, \, y, \, z, \, w) \, \,\,  = \, \,\,  \,
 \nonumber \\
\hspace{-0.98in}&&  \quad     \quad  \quad  \quad    \quad  \quad      \quad           
{x}^{2}{y}^{2}{z}^{2} \, \,\,  \,  -2\, m^2 \cdot \, x y z \cdot \, (x \, +y \, +z)
\, \,  \,\, + ( 8\cdot \, m^3 \, + \, \, x \, y\, z\, w) \cdot \, xyz
 \nonumber \\
\hspace{-0.98in}&&  \quad     \quad    \quad  \quad    \quad  \quad  \quad    \quad   \, \,  \, \, \, \,
+{m}^{4} \cdot \,  ({x}^{2}+{y}^{2}+{z}^{2})  \, \,  \,  \,  -2\,{m}^{4} \cdot \, (xy \, +xz \, +yz).
\end{eqnarray}
For  $\, z$ (resp.  $\, x$ or $\, y$) fixed the corresponding
Hauptmodul (\ref{Hauptsurfacetriquadr}) reads:
\begin{eqnarray}
\label{Hauptsurfacetriquadr2bis1}
\hspace{-0.98in}&&  \quad     \quad    \quad  \quad    \quad  \quad    \quad   \, \,  \, \,
{\cal H} \, \, = \, \, \, \,
\, {\frac {1728 \cdot \, {m}^{12} \cdot \,  p \cdot \,
(16\,{m}^{3} \, +p) }{ (16\,{m}^{6} \, +16\,{m}^{3} \cdot \, p \, +{p}^{2})^{3}}}.
\end{eqnarray}
The telescoper of the rational function in {\em four variables} $\, x, \, y, \, z$ and $\, w$,
\begin{eqnarray}
\label{Ratsurfacetriquadr2bis1}
  \hspace{-0.98in}&&  \quad     \quad   \quad      \quad
 \quad \quad    \quad  \quad    \quad  \quad    \quad   \, \,
     {{x \,y \,z } \over {   {\cal T}_{m} (x, \, y, \, z, \, w)  }}, 
\end{eqnarray}
is an order-three  (self-adjoint) linear differential operator which is the {\em symmetric square}
of an order-two linear differential operator 
having the following pullbacked $\, _2F_1$ hypergeometric solution: 
\begin{eqnarray}
\label{solRatsurfacetriquadr2bis11}
  \hspace{-0.98in}&&  \quad  \, \, \,   \quad  \quad     \quad  \quad  
 (16\,{m}^{6} \, +16\,{m}^{3} \cdot \, x \, +{x}^{2})^{-1/4} \cdot \,
 \nonumber \\
 \hspace{-0.98in}&&  \, \,  \quad  \, \,   \, \,  \quad  \quad     \quad  \quad   \quad  \quad  
\times \, \,     _2F_1\Bigl([{{1} \over {12}}, \, \, {{5} \over {12}}], \, [1], \, \,
{\frac {1728 \cdot \, {m}^{12} \cdot \,  x \cdot \, (16\,{m}^{3} \, +x) }{
(16\,{m}^{6} \, +16\,{m}^{3} \cdot \, x \, +{x}^{2})^{3}}}  \Bigr). 
\end{eqnarray}
As it should the Hauptmodul in (\ref{solRatsurfacetriquadr2bis11}) is the same as
the Hauptmodul (\ref{Hauptsurfacetriquadr2bis1}). 
The algebraic surface (\ref{surfacetriquadr2}) can be seen
as the product of {\em two times the same elliptic curve} with the
 Hauptmodul (\ref{Hauptsurfacetriquadr}) (or (\ref{Hauptsurfacetriquadr2bis1})). As expected the
 solution of the order-three telescoper is the  {\em square}  of the pullbacked $\, _2F_1$ hypergeometric
 function (\ref{solRatsurfacetriquadr2bis11}) with the  Hauptmodul (\ref{Hauptsurfacetriquadr2bis1}).

\vskip .1cm 

{\bf Remark:}  Let us perform some deformation of the rational function (\ref{Ratsurfacetriquadr2bis}),
changing the first $\, -2 \, $ coefficient in (\ref{surfacetriquadr2bis}) into a $\, -3$ coefficient.
The polynomial $\,\,  T (x, \, y, \, z, \, w)$:
\begin{eqnarray}
\label{surfacetriquadr2bis2}
 \hspace{-0.98in}&&\,  \quad  \quad      \, \, 
T (x, \, y, \, z, \, w)
\, \,\, \,  = \, \,\,  \,
 \\
\hspace{-0.98in}&&  \quad       \quad    \quad    \,   \quad  \,\,  \, \, \,   \, \,          
{x}^{2}{y}^{2}{z}^{2} \, \, \,\, -3  \,  \cdot \, x^2 y^2 z^2 \, w \cdot \, (x \, +y \, +z)
 \,  \, \,\, +4  \cdot \, (x y z w \, +1) \cdot \, x^2 y^2 z^2 \, w 
\nonumber \\
\hspace{-0.98in}&&  \quad     \quad     \quad     \quad   \quad  \,\,  \,  \quad  \quad   \, \, \, \,
 + x^2 y^2 z^2  w^2 \cdot \,  ({x}^{2}+{y}^{2}+{z}^{2})
\,  \, \,  \, -2 \cdot \, x^2 y^2 z^2  w^2 \cdot \, (xy \, +xz \, +yz).\nonumber 
\end{eqnarray}
The telescoper of the rational function in {\em four variables},
\begin{eqnarray}
\label{Ratsurfacetriquadr2bister}
\hspace{-0.98in}&&  \quad   \quad     \, \, \quad    \quad \quad    \quad  \quad    \quad  \quad    \quad   \, \,
     {{x \,y \,z } \over {  T (x, \, y, \, z, \, w)  }}, 
\end{eqnarray}
is an (irreducible) of (only) {\em order-four} linear differential operator $\, L_4$ which is non-trivially
homomorphic to its adjoint\footnote[2]{Its exterior square has a rational solution. However
  this order-four linear differential operator is not MUM
  (maximum unipotent monodromy~\cite{Almkvist,Straten,IsingCalabi})}.
A priori, we cannot exclude the fact that $L_4$ could be
homomorphic to the symmetric cube of a second-order linear differential
operator, or to a symmetric product of two second-order operators.
Furthermore, it could also be, in principle, that these second-order operators
admit classical modular forms as solutions (pullbacks of special $_2F_1$
hypergeometric functions). However, these options can both be excluded by
using some results from differential Galois theory~\cite{Singer}, specifically
from \cite[Prop.~7, p.~50]{Person} for the symmetric cube case, and from
\cite[Prop.~10, p.~69]{Person} for the symmetric product case, see
also~\cite[\S3]{Hoeij}. Indeed, if $ \, L_4$ were either a symmetric cube or a
symmetric product of order-two operators, then its symmetric square would
contain a (direct) factor of order 3 or 1. This is ruled out by a
factorization procedure which shows that the symmetric square of~$L_4$ is
(LCLM-)irreducible.

This example does not correspond to an addition
formula like (\ref{additionellipticsinus}), but the polynomial
$\,\,  T (x, \, y, \, z, \, w)$
still corresponds to a tri-quadratic (and thus an algebraic variety
with an infinite number of birational automorphisms). 

\vskip .2cm 

\subsection{Rational functions with tri-quadratic denominator: another example.}
\label{subtriquadratic2another}

The telescoper of the rational function in four variables $\, x, \, y, \, z$ and $\, w$,
\begin{eqnarray}
\label{Ratsurfacetriquadr2bisterquat}
\hspace{-0.98in}&& \quad \quad \, \, \quad \quad \quad \quad  \quad \quad  \quad  \quad   \, \,
     {{x \,y \,z } \over {  T (x, \, y, \, z, \, w)  }}, 
\end{eqnarray}
where the polynomial $\, T (x, \, y, \, z, \, w)$ almost corresponds to
the tri-quadratic (\ref{surfacetriquadr2}):
\begin{eqnarray}
\label{surfacetriquadr1}
  \hspace{-0.98in}&&  \quad  \quad  \quad  \quad \, \, \, 
T (x, \, y, \, z, \, w) \, \, = \, \, \,\, \,
{x}^{2}{y}^{2}{z}^{2} \, \,\, \,\, -2 \cdot \, x y z \cdot \, (x \, +y \, +z)
\, \, \, \, \, + 8 \cdot \, xyz
\nonumber \\
\hspace{-0.98in}&&  \quad    \quad  \quad  \quad \,  \quad   \quad   \, \, \, \, \, \,
+ ({x}^{2}+{y}^{2}+{z}^{2})  \,  \,  \, \, \, \,
-2 \cdot \, (xy \, +xz \, +yz) \, \, \, \, + x\, y\, z\, w.
\end{eqnarray}
is an {\em order-four} linear differential operator non-trivially\footnote[5]{The intertwiners
 are of order one and order three. This order-four linear differential operator is not MUM
  (maximum unipotent monodromy~\cite{Almkvist,Straten,IsingCalabi}).}
homomorphic to its adjoint. 

\vskip .2cm

{\bf Remark:} If one (slightly) changes the first coefficient $\, -2$
into $\, -3$ in (\ref{surfacetriquadr1}) 
\begin{eqnarray}
\label{surfacetriquadrbis}
  \hspace{-0.98in}&&  \quad  \quad  \quad  \quad  \quad  
T (x, \, y, \, z, \, w) \, \, = \, \, \,\,
{x}^{2}{y}^{2}{z}^{2} \, \,\, \,\, -3 \cdot \, x y z \cdot \, (x \, +y \, +z)
\, \, \, \, \, + 8 \cdot \, xyz
\nonumber \\
\hspace{-0.98in}&&  \quad    \, \,   \quad  \quad  \quad  \quad  \quad   \, \, \, \, \, \,
+ ({x}^{2}+{y}^{2}+{z}^{2})  \,  \,  \, \, \, \,
-2 \cdot \, (xy \, +xz \, +yz) \, \, \,\, + x\, y\, z\, w.
\end{eqnarray}
one obtains an {\em order-six} telescoper for  rational function
of four variables (\ref{Ratsurfacetriquadr2bisterquat}).
This order-six linear differential operator is non trivially\footnote[2]{The intertwiners
 are of order three and order five. This order-six linear differential
  operator is not MUM~\cite{Almkvist,Straten,IsingCalabi}.}
homomorphic to its adjoint. 

\vskip .2cm 

\subsection{Rational functions with tri-quadratic denominator.}
\label{subtriquadratic2}

Let us consider the most general {\em tri-quadratic surface} 
\begin{eqnarray}
\label{triquadraBIS}
\hspace{-0.98in}&& \, \,  \quad  \, \quad  \quad   \quad  \quad  \quad  \quad 
\sum_{m=0,1,2} \, \sum_{n=0,1,2} \, \sum_{l=0,1,2} \, a_{m,n,l} \cdot \, x^m \, y^n \, z^l
\, \, \,\,  = \, \,  \, \, \,  0, 
\end{eqnarray}
depending on $\,\, 27 = \, 3^3 \,$ parameters $ \,  a_{m,n,l}$. It 
can be rewritten as:
\begin{eqnarray}
\label{triquadrarewritt}
\hspace{-0.98in}&& \, \,  \quad  \quad  \, \quad   \quad   \quad   \quad \quad 
A(x, \, y) \cdot \, z^2 \, \,\,  +  B(x, \, y) \cdot \, z \, \,\, \, +  C(x, \, y)
\, \, \, = \, \,\,   \, \,  0. 
\end{eqnarray}
It is straightforward to see that condition (\ref{triquadrarewritt}) is preserved
by the {\em birational involution} $\, I_z$ 
\begin{eqnarray}
\label{involutionbirat}
\hspace{-0.98in}&& \, \,  \quad \quad   \quad   \quad    \quad  
I_z: \quad \quad \Bigl(x, \, \, y, \, \, z\Bigr)
\quad \, \,  \,   \longrightarrow  \,  \,  \,  \,  \,  \,  \, \, \, \,  
\Bigl(x, \, \,\,  y, \, \, \, \, {{ C(x, \, y)} \over { A(x, \, y) }} \cdot {{1 } \over { z}}\Bigr), 
\end{eqnarray}
and we have of course two other similar {\em birational involutions} $\, I_x$ and $\, I_y$ that single out
$\, x$ and $\, y$ respectively. The (generically) {\em infinite-order} birational transformations
$\, K_x\, = \,  I_y \cdot \, I_z$, $\, K_y\, = \,  I_z \cdot \, I_x \, $
and $\, K_z\, = \, I_x \cdot \, I_y \, $ are birational symmetries of the
surface (\ref{triquadraBIS}) or (\ref{triquadrarewritt}). 
They are related by $\,\,  K_x \cdot \, K_y \cdot \, K_z \, = \, \, identity$. Note
that the {\em birational transformation} $\, K_x$
{\em preserves} $\, x$. The iteration of the (generically) {\em infinite-order} birational transformation $\, K_x$
gives {\em elliptic curves}. Since equation (\ref{triquadraBIS}) or (\ref{triquadrarewritt})
is preserved by $\, K_x$, {\em which also preserves} $\, x$, the equation of the {\em elliptic curves} corresponding
to the  iteration\footnote[1]{The birational transformation $\, K_x$ maps the elliptic curve
  onto itself (self-map). One can  use the iteration of the birational transformation $\, K_x$  to  actually visualise
  the elliptic curve~\cite{Quasi,Noether}.}  of $\, K_x$ is (\ref{triquadraBIS})
{\em for fixed values of} $\, x$. Equation
(\ref{triquadraBIS}), for fixed values of $\, x$, is a (general) biquadratic curve in $\, y$
and $\, z$ and is thus {\em an elliptic curve depending on} $\, x$. Therefore one has a canonical foliation of the
algebraic surface (\ref{triquadraBIS}) in elliptic curves. Of course the  iteration of $\, K_y$ (resp. $\, K_z$)
also yields elliptic curves, and similarly yields two other foliations in elliptic curves.

\vskip .2cm

We have a foliation in two families of elliptic curves $\,{\cal E}$ and $\, {\cal E'}$
of the surface. Consequently, this
tri-quadratic surface (\ref{triquadraBIS}), having an {\em infinite set}
of {\em birational automorphisms},  an {\em infinite set} of {\em birational symmetries},
cannot be of the ``general type'' (it has Kodaira dimension less than $\, 2$).

\vskip .1cm
\vskip .1cm

\subsection{Rational functions with tri-Quadratic denominator: Fricke cubics examples associated with Painlev\'e VI equations}
\label{Painleve}

Let us consider more very simple examples of tri-quadratic surfaces that occur in different domains of mathematics and physics.

Among the {\em Fricke families} of cubic surfaces, the family~\cite{Fricke,Boalch,CantatLoray} 
\begin{eqnarray}
\label{Fricke}
\hspace{-0.98in}&& \quad  \quad \quad \quad \quad 
x\, y \, z \, \, \, \, + x^2 \, +y^2 \, +z^2 \, \, \, +b_1\, x \, +b_2\, y \,  +b_3\, z \,\, \,  \, +c
 \, \, = \,\,  \, \, 0,
\end{eqnarray}
of affine cubic surfaces parametrised by the four constants $(b_1,\,b_2,\,b_3,\,c)$
is known~\cite{Boalch} to be a deformation of a $\, D_4$ singularity which occurs at the
symmetric (Manin's) case  $ \, b_1=\,b_2=\,b_3= \, -8, \,c\, = \, 28$.

Among the symmetric $\, b_1=\,b_2=\,b_3 \, $ cases
some selected sets of  the four constants  $\, (b_1,\,b_2,\,b_3,\,c)$ emerge: the Markov cubic
$\, b_1=\,b_2=\,b_3\, = c \, = 0$, Cayley's nodal cubic  $\, b_1=\,b_2=\,b_3\,\, = 0, \, \, c \, = -4$,
Clebsch diagonal cubic  $\, b_1=\,b_2=\,b_3\,\, = 0, \, \, c \, = -20$, and Klein's cubic 
$\, b_1=\,b_2=\,b_3\,\, = -1, \, \, c \, = 0$.

Some of these  symmetric cubics play can be seen
as the {\em monodromy manifold} of the {\em Painlev\'e VI equation} (see equation (1.7) in ~\cite{MazoVidu},
see also equations (1.2) and (1.4) in~\cite{CantatLoray}): 
the Picard-Hitchin cases $\, (0,0,0, \, 4)$,  $\, (0,0,0,\, -4)$,  $\, (0,0,0, \, -32)$, the Kitaev's cases
$\, (0,0,0, \, 0)$, $\, (-8,-8,-8, \, -64)$, and especially the Manin's case  $\, (-8,-8,-8, \, 28)$.

\vskip .2cm

Let us consider the  Picard-Hitchin example $\, (0,0,0,\, -4)$ as a denominator of a rational
function~\cite{Boalch}. Let us consider the rational function in three variables
$\, x, \, y$ and $\, z$~\cite{Boalch}:
\begin{eqnarray}
\label{RatPainleve}
  \hspace{-0.98in}&& \quad  \quad \quad \quad \quad  \quad  \quad 
\, \, 
R(x, \, y, \, z)  \, \, \,  = \, \,  \quad 
 {{1} \over {  x^2  + y^2 \, + z^2  \,  \,\, + x\,y\,z \,\, \,  \, -4}}.
\end{eqnarray}
The telescoper of the rational function (\ref{RatPainleve})
is actually an {\em order-two} linear differential operator $\, L_2$
\begin{eqnarray}
\label{RatPainleveL2}
  \hspace{-0.98in}&& \quad  \quad \,  
 L_2   \, \,  = \, \, \, \, \,
 2 +x  \,\, \,   + (3\,x^2 \, +14\, x \, -8) \cdot \, D_x
 \, \,  \,  + \,   x \cdot \, (x+8) \cdot \, (x-1) \cdot \, D_x^2,           
\end{eqnarray}                 
which has the pullbacked hypergeometric solution\footnote[1]{Note
  the emergence of the pullback
  $\, -27\, x^2/(x-4)^3$ that we already saw in
  (\ref{K3formsol}) and in (\ref{K3formsolHauptymore}).}: 
\begin{eqnarray}
  \hspace{-0.98in}&&  \quad \,\,\,
 {{1} \over {x\, +2}} \cdot \,
_2F_1\Bigl([{{1} \over {3}}, \, {{2} \over {3}}], \, [1],
\, \, {{ 27 \, x} \over {(x\, +2)^3}} \Bigr)  \, \, = \, \, \,
 -\, {{2} \over {x\, -4}} \cdot \,
_2F_1\Bigl([{{1} \over {3}}, \, {{2} \over {3}}], \, [1],
 \, \,  -\, {{ 27 \, x^2} \over {(x\, -4)^3}} \Bigr)  \nonumber 
  \\
\label{RatPainleveSol}
 \hspace{-0.98in}&&   \quad \quad \quad \quad  
\, \, = \, \, \, \Bigl((x\, +2) \cdot \, ({x}^{3}+6\,{x}^{2}-12\,x+8) \Bigr)^{-1/4} 
 \\
 \hspace{-0.98in}&& \quad  \, \quad \quad \quad \quad \quad \quad \quad \quad  
                    \times \, 
 _2F_1\Bigl([{{1} \over {12}}, \, {{5} \over {12}}], \, [1], \, \,
\,{\frac {1728 \cdot \, {x}^{3}\cdot \,  (x+8)  \cdot\, (x-1)^{2}}
   { (x \, +2)^{3} \cdot ({x}^{3}+6\,{x}^{2}-12\,x+8)^{3}}} \Bigr)
 \nonumber
\end{eqnarray}
\begin{eqnarray}
  \label{RatPainleveSolsol}
 \hspace{-0.98in}&&   \quad \quad \quad \quad  
\, \, = \, \, \, 2 \cdot \, \Bigl((x\, -4) \cdot \, (({x}^{3}+12\,{x}^{2}+48\,x-64) \Bigr)^{-1/4} 
  \\
  \hspace{-0.98in}&& \quad  \, \quad \quad \quad \quad \quad \quad \quad \quad  
   \times \, 
 _2F_1\Bigl([{{1} \over {12}}, \, {{5} \over {12}}], \, [1], \, \,
 \, - \,{\frac {1728 \cdot \,  {x}^{6} \cdot \, (x\, -1)  \cdot \, (x \, +8)^{2}}{
 (x \, -4)^{3} \cdot \, ({x}^{3}+12\,{x}^{2}+48\,x-64)^{3}}}\Bigr).
\nonumber 
\end{eqnarray}

\vskip .2cm

{\bf Remark:} Note that the two Hauptmoduls in  (\ref{RatPainleveSol}) and (\ref{RatPainleveSolsol})
are related by the involution $\, x \, \, \longleftrightarrow \, \, -8/x$. This symmetry of the problem
corresponds to the fact that the order-two telescoper $\, L_2$ is simply conjugated to its pullback
by $\, x \, \, \longrightarrow \, \, -8/x$.
\vskip .2cm
Eliminating $\, z \, = \, p/x/y \, $ in the denominator of (\ref{RatPainleve})
gives the {\em genus-four} algebraic curve:
\begin{eqnarray}
\label{RatPainleveSoleq}
 \hspace{-0.98in}&& \quad  \quad \quad \quad \quad \quad \quad
 {x}^{4}{y}^{2} \,\,\, +{x}^{2}{y}^{4} \,\, \, +(p\, -4) \cdot \, {x}^{2}{y}^{2}\, \,\, +{p}^{2}
\,  \,\,\, =  \, \, \,   \,0.  
\end{eqnarray}
The question is to see whether the Jacobian of this {\em genus-four} algebraic curve
(\ref{RatPainleveSoleq}) could also correspond to a split Jacobian, with a $\, j$-invariant
corresponding to the Hauptmoduls in (\ref{RatPainleveSol}) or (\ref{RatPainleveSolsol}).

\vskip .2cm
\vskip .1cm

More generally  the symmetric rational function in three variables
$\, x, \, y$ and $\, z$~\cite{Boalch}:
\begin{eqnarray}
\label{RatPainleveSymc}
  \hspace{-0.98in}&& \quad  \quad \quad \quad \quad  \quad  \quad 
\, \,  \, 
R(x, \, y, \, z)  \, \, \,  = \, \,  \quad 
 {{1} \over {  x^2  + y^2 \, + z^2  \,\, \,  \, + x\,y\,z \,\, \,  \, +c}},
\end{eqnarray}
which takes into account the other Picard-Hitchin
cases\footnote[2]{As well as the Markov cubic
  $\, b_1=\,b_2=\,b_3\, = c \, = 0$, Cayley's nodal cubic
  $\, b_1=\,b_2=\,b_3\,\, = 0, \, \, c \, = -4$,
and Clebsch diagonal cubic  $\, b_1=\,b_2=\,b_3\,\, = 0, \, \, c \, = -20$ cases.
} $\, (0,0,0, \, 4)$,  $\, (0,0,0,\, -4)$,  $\, (0,0,0, \, 32)$, 
also has an {\em order-two} telescoper which has a simple pullbacked
hypergeometric solution:
\begin{eqnarray}
\label{RatPainleveSolc}
  \hspace{-0.98in}&& \quad  \quad \quad
 {{1} \over {x\, +c}} \cdot \,
 _2F_1\Bigl([{{1} \over {3}}, \, {{2} \over {3}}], \, [1],
       \, \, -\, {{ 27 \, x^2} \over {(x\, +c)^3}} \Bigr) 
   \\
 \hspace{-0.98in}&&  \quad \quad \quad \quad \quad
   \, = \, \, \, \,  \,\,   p_6(x)^{-1/6} \cdot \, 
 _2F_1\Bigl([{{1} \over {12}}, \, {{7} \over {12}}], \, [1], \, \,
 \,  {{ 1728 \cdot \, x^6 \cdot \,  p_3(x)  } \over { p_6(x)^2}} \Bigr)
  \nonumber   \\
 \hspace{-0.98in}&&  \quad \quad \quad \quad \quad
\, = \, \, \, \, \, \,  (x\, +c)^{-1/4} \cdot \,   q_3(x)^{-1/4} \cdot \,  
 _2F_1\Bigl([{{1} \over {12}}, \, {{5} \over {12}}], \, [1], \, \,
\,  -\, {{ 1728 \cdot \, x^6 \cdot \,  p_3(x)  } \over { (x\, +c)^3 \cdot \, q_3(x)^3 }} \Bigr),
\nonumber  
\end{eqnarray}
where\footnote[1]{The values $\, c=0$ and $\, c= \, -4$ are the only values such that
  the discriminant in $\, x$ of $\, p_3(x)$ can be zero.}:
\begin{eqnarray}
\label{RatPainleveSolwhere}
  \hspace{-0.98in}&& \quad  \quad  \quad  
p_3(x)  \, \,  = \, \,      \,  \,  x^3 \,  \,  +3 \cdot \, (c+8) \cdot \, x^2
\,  \, +3 \cdot \, c^2 \cdot \, x \,  \,  +c^3,
 \nonumber \\
  \hspace{-0.98in}&& \quad  \quad  \quad
q_3(x)  \, \,  = \, \,   \,     \,
 x^3 \,  \,  +3 \cdot \, (c+9) \cdot \, x^2  \,  \, +3 \cdot \, c^2 \cdot \, x \,  \,  +c^3,             
                     \nonumber \\
 \hspace{-0.98in}&& \quad  \quad  \quad    p_6(x)  \, \,  = \, \,   \,     \,
x^6 \, \,  \,  +6 \cdot \, (c+6) \cdot \, x^5 \, \,  \,  +(216 \, +108\, c \,\,  +15\, c^2) \cdot \, x^4 \,
 \nonumber \\
 \hspace{-0.98in}&& \quad  \quad  \quad \quad  \quad \quad  \quad 
+(20\, c \, +108) \cdot \, c^2 \cdot \, x^3     \,   \, \,
+(15\, c \, +36) \cdot \, c^3  \cdot \, x^2 \, \,  \,  +6 \cdot \, c^5\, x \, \,   \, +c^6.      
\nonumber 
\end{eqnarray}
Eliminating $\, z \, = \, p/x/y \, $ in the denominator of (\ref{RatPainleve})
gives the {\em genus-four} algebraic curve:
\begin{eqnarray}
\label{RatPainleveSymc}
 \hspace{-0.98in}&& \quad  \quad \quad \quad \quad \quad \quad \quad
 {x}^{2}{y}^{2} \cdot (x^2 \, + \, y^2) \, \,\, \, +(p\, +c) \cdot \, {x}^{2}{y}^{2} \, \, \,\, +{p}^{2}
\,  \,\,\, =  \, \, \,   \,0.  
\end{eqnarray}
Again, the question is to see whether the Jacobian of this {\em genus-four}
algebraic curve (\ref{RatPainleveSymc})
could also correspond to a split Jacobian, with a $\, j$-invariant
corresponding to the Hauptmodul in (\ref{RatPainleveSolc}).

\vskip .2cm

{\bf Remark:} Note after~\cite{Boalch} that the value $\, c \, = \, -4 \, \, $ is particular. It is such that
the denominator
\begin{eqnarray}
\label{deffxyz}
 \hspace{-0.98in}&& \quad  \quad \quad \quad \quad \quad  \quad
\, f(x,\, y, \, z; \,  c) \, \,  = \, \,  \, \,
 x^2  + y^2 \,  + z^2  \,\, \,  \,  \, + x\,y\,z  \,  \,\, \,  \, +c,
\end{eqnarray}
when transformed by the simple quadratic transformation 
\begin{eqnarray}
\label{RatPainleveSymctransf}
 \hspace{-0.98in}&& \quad  \quad \quad \quad \quad  \quad
 (x, \, y, \, z) \, \quad  \quad  \longrightarrow
  \quad \quad  \quad \Bigl(2 \, -x^2, \,\,\, 2 \, -y^2, \, \,\,2\, -z^2\Bigr), 
\end{eqnarray}
factorises nicely:
\begin{eqnarray}
\label{Rafactonice}
\hspace{-0.98in}&& \,   \quad \, \,  \, 
f\Bigl(2 \, -x^2, \, 2 \, -y^2, \, 2\, -z^2; \, -4\Bigr)
\, \, = \,\,  \, f(x,\, y, \, z; \,  -4) \cdot \, f(-x,\, y, \, z; \,  -4).
\end{eqnarray}
In other words we have an {\em endomorphism of the Cayley cubic surface}.

\vskip .2cm 

\subsubsection{Singular symmetric Fricke surface.}
\label{singular}

Let us consider the rational function
\begin{eqnarray}
\label{singularFricke}
\hspace{-0.98in}&& \, \,   \,   \quad \quad\quad \quad
 R(x, \, y, \, z) \, \, = \, \, \,
  {{1} \over {xyz \, + x^2 \, +y^2 + \, z^2 \,\,  + b\cdot \, (x\, +y\, +z) \,\,  +c }}. 
\end{eqnarray}
The vanishing condition of the denominator of (\ref{singularFricke}) is a symmetric Fricke surface
which, according to Lemma 9 in~\cite{Boalch}, is singular for
\begin{eqnarray}
\label{singularFrickecond}
  \hspace{-0.98in}&& \, \,   \,   \quad \quad \quad \quad \quad \quad\quad \quad
   b^2\, -8\, b \, -16 \, -4 \, c \, \, = \, \, \, 0,  
\end{eqnarray}
and:
\begin{eqnarray}
\label{singularFrickecond2}
  \hspace{-0.98in}&& \, \,   \,   \quad \quad \quad \quad  \quad \quad \quad
   4\, b^3 \, -3\, b^2\, -6\, b \, c  \,  + \, c^2 \, + 4 \, c \, \, = \, \, \, 0.
\end{eqnarray}
For instance, for $\, b \, = \, -8$, the first condition (\ref{singularFrickecond})
gives $\, c \, = \, 28$ (i.e. Manin's case) and the second condition
(\ref{singularFrickecond2}) gives  $\, c \, = \, 28$  and  $\, c \, = \, -80$. 

The calculation of the telescoper of (\ref{singularFricke}) in the singular
$\, (b, \, c) \, = \, \, (-8, \, 28)$
case gives an (irreducible) {\em order-four} linear differential operator which is (non-trivially)
{\em homomorphic to its adjoint}\footnote[2]{The intertwiners are of order-two. The exterior  square
 of that operator has a simple rational
 solution $\, (x^2 \, +39\,x-168)/(x+343)/x/(x-8)^2/(x-9)$. We have a similar result for
$\, (b, \, c) \, = \, \, (-8, \, -80)$, the exterior  square
 of that operator having the rational
 solution $ \, p_3(x)/x/(x+64)/(x-125)/q_3(x)$,
 where $\, p_3(x)= \, x^3-149\,x^2+34080\,x-3010560 \, $ and $\, q_3(x)=x^3-349\,x^2+38656\,x-1032192 $.}.

\vskip .1cm 

\vskip .2cm 

\subsection{Rational functions with $\, N$-Quadratic denominator.}
\label{subtriquadraticN}

The calculations of subsection \ref{subtriquadratic2} $\,\, $ can straightforwardly 
be generalised to $\, N$-quadratic equations,
writing the $\, N$-quadratic (\ref{Nquadra}) as
\begin{eqnarray}
\label{triquadrarewritt1}
\hspace{-0.98in}&& \, \,   \,   \quad \quad    \quad     \quad    \quad    \quad   
A(x_1, \, x_2, \, \cdots, \, x_{N-1}) \cdot \, x_N^2 \, \,\,\,\,
 +B(x_1, \, x_2, \, \cdots ,\, x_{N-1}) \cdot \, x_N \, \, \,
\nonumber \\
\hspace{-0.98in}&& \, \,  \quad \quad \quad \quad \quad \quad    \quad   \quad \quad \quad \quad   
 +  C(x_1, \, x_2, \, \cdots, \, x_{N-1})  \,  \, \, \,\, = \, \,\,   \, \,  0,
\end{eqnarray}
and introducing the 
 {\em birational involution} $\, I_N$ 
\begin{eqnarray}
\label{involutionbiratbis}
\hspace{-0.98in}&& \, \,  \quad   \quad  
I_N: \, \quad \quad   \Bigl(x_1, \, x_2, \, \cdots, \, x_{N}\Bigr)
 \\
 \hspace{-0.98in}&& \, \,  
 \, \quad \quad \quad \quad  \quad  \quad  \,   \, \,   \longrightarrow  \,  \, \,   \,   \quad   \,   
 \Bigl(x_1, \, x_2, \, \cdots, \, x_{N-1}, \, \, \, \,
  {{ C(x_1, \, x_2, \, \cdots, \, x_{N-1})} \over { A(x_1, \, x_2, \, \cdots \, x_{N-1}) }}
    \cdot {{1 } \over { x_N}}\Bigr).
 \nonumber
\end{eqnarray}
Similarly to subsection \ref{subtriquadratic2}, we can introduce $\, N$ involutive birational transformations
$\, I_m$ and consider the products of two such involutive birational transformations
$\, K_{m,n} \, = \, \, I_m \cdot I_n$.
These  $\, K_{m,n}$'s are (generically) infinite order birational
transformations preserving the $\, N-2$ variables that are not $\, x_m$ and $\, x_n$.

Using such remarkable $\, N$ variables algebraic varieties, with an {\em infinite set of birational
automorphisms}, one can build rational functions of $\, N+1 \, $ variables, any of the parameter
of the  algebraic variety, becoming an arbitrary rational\footnote[1]{Or even an arbitrary algebraic function
  of the product $\, p \, = \, x_1 \, x_2  \, \cdots \, x_N$, or a transcendent series
analytic at $\, p \, = \, 0$.} function of the product $\, p \, = \, x_1 \, x_2  \, \cdots \, x_N \, $
in order to build the denominator of the rational function. The telescopers of such 
rational functions is seen (experimentally using creative telescoping) to be  of substantially  smaller
order than the one for rational functions where their  denominators are,
after reduction by $\, p \, = \, x_1 \, x_2  \, \cdots \, x_N$, associated with algebraic varieties
of the ``general type''.

\vskip .1cm

\section{Telescopers of rational functions of several variables}
\label{several}

In our previous paper~\cite{HeunJPA,malala}, dedicated to Heun functions that are solutions of telescopers
of simple rational functions of (most of the time) four variables, we have obtained many 
order-three telescopers having square of pullbacked $\, _2F_1$ hypergeometric solutions. Recalling
sections \ref{productelliptic}, \ref{Weierstrass}, or even \ref{split} in~\cite{malala}, it is natural to imagine, 
for these examples in~\cite{HeunJPA,malala} yielding square of pullbacked $\, _2F_1$ hypergeometric functions,
a scenario where, after elimination of the fourth variable ($w \, = \, p/x/y/z$)
in the denominator of the rational function
of four variables, the corresponding algebraic surface $\, S(x, \, y, \, z) \, = \, \, 0$,
in the remaining three variables,
could be seen as $\, K_3$ surface (\ref{K3form}) which can be seen as
{\em associated with the product of two times the same elliptic curve}, or
other ``Periods~\cite{KontZagier} of extremal rational surfaces''
scenario.  Some other cases of similar rather simple rational
functions of four variables, yield order-two telescopers with  pullbacked $\, _2F_1$ hypergeometric functions
(but not square or products of pullbacked $\, _2F_1$ hypergeometric functions).

\vskip .1cm

$\bullet$ Let us consider the rational function in {\em four} variables $\, x, \, y, \, z, \, u$: 
\begin{eqnarray}
\label{Ratfonc4first}
  \hspace{-0.98in}&&  \,  \quad  \quad \quad \quad \, \, 
\, \,  
R(x, \, y, \, z, \, u)  \, \, \,  = \, \,  \,   \,  
 {{1} \over { 1 \,\, \,\,+3\, y \,\,+z \,\, \,\,+9 \,y \,z \,
\, \, \,  + 11 \,z^2\,y  \, \, +3 \,u \,x}}.
\end{eqnarray}
The telescoper of this rational function of four variables is an order-two linear differential
operator $\, L_2$  which has the pullbacked hypergeometric solution:
\begin{eqnarray}
\label{telescRatfonc4u}
\hspace{-0.98in}&& \quad \, \, \, \, 
 (1\, - \,  2592\, x^2)^{-1/4} 
 \\
\hspace{-0.98in}&& \quad \quad \,  \, \, \,\,\times \,
 _2F_1\Bigl([{{1} \over {12}}, \, {{5} \over {12}}], \, [1], \,\, 
 - \, {{  419904 \cdot \, x^3 \cdot \, (5 \, -12\, x \, -19440 \, x^2 \, +2665872\, x^3) } \over {
(1 \, - \, 2592\, x^2)^3}} \Bigr).
 \nonumber 
\end{eqnarray}
The diagonal of (\ref{Ratfonc4first}) is the expansion of this pullbacked hypergeometric
function (\ref{telescRatfonc4u}):
\begin{eqnarray}
\label{telescRatfonc4uexp}
  \hspace{-0.98in}&& \quad \, \, \, \,
 1 \, \, +648\,{x}^{2} \, \, -72900\,{x}^{3} \, \,
 +1224720\,{x}^{4} -330674400\,{x}^{5} \, +23370413220\,{x}^{6} \,
  \\
  \hspace{-0.98in}&& \quad \quad \, \,           \,
\, -1276733858400\,{x}^{7} \, +180019474034400\,{x}^{8}
 \, -12013427240614800\,{x}^{9} \, \, \, + \, \, \, \cdots \nonumber 
\end{eqnarray}
If one considers the intersection of the vanishing condition of the denominator of (\ref{Ratfonc4first})
with the hyperbola $\, p \, = \, \, x \, y \, z \, u$, eliminating for instance $\, u \, = \, p/x/y/z$
in the vanishing condition of the denominator of (\ref{Ratfonc4first}), one gets a condition,
{\em independent of} $\, x$, which corresponds to  a {\em genus-one} curve
\begin{eqnarray}
\label{g1}
\hspace{-0.98in}&& \quad \, \, \quad \quad \quad \quad  \quad \quad 
 11\,{y}^{2}{z}^{3} \, \, +9\,{y}^{2}{z}^{2}\,  \, +3\,{y}^{2}z \,\,  +y{z}^{2} \, +yz \,\,   \, +3\,p
 \, \,  = \,  \, \, 0.
\end{eqnarray}
The Hauptmodul of this elliptic curve (\ref{g1}) reads:
\begin{eqnarray}
\label{Hauptg1}
  \hspace{-0.98in}&& \quad \quad \quad \quad  \quad \, \,
{\cal H} \,  \, = \, \, \,
- \, {{  419904 \cdot \, p ^3 \cdot \, (5 \, -12\, p \, -19440 \, p^2 \, +2665872\, p^3)} \over {
(1 \, - \, 2592\, p^2)^3}},         
\end{eqnarray}
which corresponds precisely to the Hauptmodul pullback in (\ref{telescRatfonc4u}).

\vskip 3mm

$\bullet$ Let us, now, generalize the rational function (\ref{Ratfonc4first})
of {\em four} variables $\, x, \, y, \, z, \, u$,
introducing the rational function of $\, N+3\,$ variables
$\, x, \, y, \, z, \, u_1, \, u_2, \, \cdots,  \, u_N$:
\begin{eqnarray}
\label{Ratfonc4firstN}
  \hspace{-0.98in}&&  \,  \quad \quad  \quad \quad 
\, \,  
R(x, \, y, \, z, \,  u_1, \, u_2, \, \cdots, \, u_N) 
 \\
\hspace{-0.98in}&&  \, \quad  \quad  \quad \, \quad \quad \quad \, \, \, \,  = \, \, \,  \,  \,  \,   \,
{{1} \over { 1 \,\, \,\,+3\, y \,\,+z \,\, \,\,+9 \,y \,z \,
\, \, \,  + 11 \,z^2\,y \, \,   \, \, +3 \, x \cdot \, u_1 \, u_2 \, \, \cdots \, u_N }}.
\nonumber
\end{eqnarray}
The telescoper of this rational function of  $\, N+3$ variables is the same order-two
telescoper as for (\ref{Ratfonc4first}), 
which has the pullbacked hypergeometric solution (\ref{telescRatfonc4u}).
Again one can verify that the diagonal  of (\ref{Ratfonc4firstN}) is the expansion (\ref{telescRatfonc4uexp})
of the pullbacked hypergeometric
function\footnote[2]{A pure algebraic geometer will probably consider this result
  as trivial from the computational point of view, saying
  that the variety is a fiber bundle over a family of elliptic curves with constant fiber (see also below).}
(\ref{telescRatfonc4u}). 
If one considers the intersection of the vanishing condition
of the denominator of (\ref{Ratfonc4firstN})
with the hyperbola $\, p \, = \, \, x \, y \, z \, u_1 \, u_2 \, \, \cdots \, u_N $, eliminating
for instance $\, u_N \, = \, p/x/y/z/u_1/\cdots/u_{N-1} \, $
in the vanishing condition of the denominator of (\ref{Ratfonc4firstN}), one gets again a condition,
{\em independent of} $\, x$, which corresponds to  a {\em genus-one} curve (\ref{g1}):
\begin{eqnarray}
\label{g1N}
\hspace{-0.98in}&& \quad \quad \quad \quad \quad \quad \quad \, \,
 11 \,{y}^{2}{z}^{3} \, \, +9\,{y}^{2}{z}^{2} \, \, +3\,{y}^{2}z \,\,
 +y{z}^{2} \,\,  +yz \,\, \,   \,  +3\,p \, \,  \, \, = \,\,  \, 0.
\end{eqnarray}
The Hauptmodul of this elliptic curve (\ref{g1N}), or (\ref{g1}) reads again
the Hauptmlodul (\ref{Hauptg1}) 
which corresponds precisely to the Hauptmodul pullback in (\ref{telescRatfonc4u}).

\vskip .1cm
\vskip .1cm

{\bf Remark:}  Recalling subsections \ref{nine} and \ref{ten}, one can consider
a nine-parameters biquadratic in two variables, or a selected ten-parameters bicubic
like (\ref{Ratfoncplusplusplus}), where the parameters are now {\em functions} of the product
of $\, N$-variables $\, p \, = \, \, x_1 \, x_2  \, \cdots \, x_N$. This will yield
an algebraic variety of $\, N$ variables (that are not on the same footing) that
will automatically be foliated in elliptic curves. 

Simple other examples are displayed in \ref{limitapp}, and one sees (experimentally)
that the Hauptmodul of the
pullbacked $\, _2F_1$ hypergeometric functions can be seen as corresponding to some
$\, x \, \rightarrow \, 0$ limit of Hauptmoduls of the elliptic curves foliating the previous
algebraic surface. In contrast with the other examples and results of this paper, we have no
algebraic geometry interpretation of this experimental result yet.

\vskip .2cm
\vskip .2cm

\section{Conclusion}
\label{Conclusion}

We have shown that the results we had obtained on diagonals of
nine and ten parameters families of rational functions, 
using creative telescoping yielding classical modular forms expressed as
pullbacked $\, _2F_1$ hypergeometric functions~\cite{DiagJPA,unabri}, can be obtained 
much more efficiently calculating the $\, j$-invariant of
an {\em elliptic curve canonically associated with the denominator of
the rational functions}.  In the case where creative telescoping yields
pullbacked $\, _2F_1$ hypergeometric functions,
we generalize this result to other families of rational functions
of three, and even more than three, variables,
when the denominator can be associated with products of elliptic curves or foliation
in terms of elliptic curves, or when the
denominator is associated with a {\em genus-two curve} with {\em a split Jacobian}
corresponding to {\em products of elliptic curves}.

\vskip .1cm

We have seen different scenarii. In the first cases, we have
considered denominators corresponding
to {\em products} of elliptic curves: in these cases the solutions
of the telescoper were {\em  products} of pullbacked $\, _2F_1$ hypergeometric
functions. We have also considered denominators corresponding to {\em genus-two} curves
with {\em split Jacobians isogenous to products of two elliptic curves},
and in these  cases the solutions of the telescoper were {\em sums} of
two pullbacked $\, _2F_1$ hypergeometric functions, sometimes
one  pullbacked $\, _2F_1$ hypergeometric function being enough to
describe the two Galois-conjugate $\, j$-invariants (see  \ref{splitgood}).
We also considered denominators corresponding to algebraic varieties
with elliptic foliations, the Hauptmodul pullback in the pullbacked $\, _2F_1$
hypergeometric functions emerging from a selected elliptic
curve of the foliation ($x\,=\, 0$, see \ref{limitapp2}, \ref{limitapp1}).
We also encountered denominators corresponding to algebraic manifolds with
an infinite set of birational automorphisms and 
elliptic curves foliation yielding, no longer
classical modular forms represented as pullbacked $\, _2F_1$
hypergeometric functions, but more general modular structures
associated with selected linear differential operators
like Calabi-Yau linear differential operators~\cite{Almkvist,Straten}
and their generalisations. 

\vskip .1cm

The creative telescoping method on a rational function is a way
to find the periods of an algebraic variety over {\em all possible cycles}\footnote[2]{Not only
the  {\em vanishing cycles} corresponding to  {\em diagonals} of rational functions.}. The fact that
the solution of the telescoper corresponds to ``Periods"~\cite{KontZagier} {\em over all possible cycles}
is a simple consequence of the fact that creative telescoping  corresponds to {\em purely differential
algebraic manipulations} on the integrand {\em independently of the cycles},
thus {\em being blind to analytical details}. In this paper, we show that the final result emerging from
differential algebra procedure (which can be cumbersome when the result depends
on nine or ten parameters), can be obtained almost instantaneously from a more fundamental
intrinsic pure algebraic geometry approach, calculating the $j$-invariant of some
canonical elliptic curve. This corresponds to a shift Analysis $\, \rightarrow \, $  Differential Algebra 
$\, \rightarrow \, $  Algebraic Geometry. Ironically, algebraic geometry studies of more involved
algebraic varieties than product of elliptic curves, foliation in elliptic curves
(Calabi-Yau manifolds, ...) is often a tedious 
and/or difficult task (finding Igusa-Shiode invariants, ...), and formal calculations tools are not always
available or user-friendly. For such involved algebraic varieties the creative telescoping may then
become a simple and efficient tool to perform effective algebraic geometry studies.

\vskip .3cm 
\vskip .2cm

{\bf Acknowledgments.}  
J-M. M. would like to thank G. Christol for many enlightening 
discussions on diagonals of rational functions. 
J-M. M. would like to thank the School of Mathematics and 
Statistics of Melbourne University where
part of this work has been performed. 
S. B. would like to thank the LPTMC and the CNRS for kind support.
We thank Josef Schicho for providing the demonstration of the results 
of \ref{MaxParam}.
We would like to thank  A. Bostan
for useful discussions on creative telescoping.
Y. A. and C.K. were supported by the Austrian Science Fund (FWF): F5011-N15.
We thank the Research Institute for Symbolic Computation,
 for access to the RISC software packages.  We thank M. Quaggetto 
for technical support. 
This work has been performed 
without any support of the ANR, the ERC or the MAE, or any PES of the CNRS.

\appendix

\section{Diagonals of rational functions and Picard-Fuchs equations}
\label{Griffiths}

For simplicity let us consider a rational function of three variables.  
The diagonal of a rational function of three variables is obtained through
its multi-Taylor expansion~\cite{Short,Big}
\begin{eqnarray}
\label{multiTaylorapp}
\hspace{-0.98in}&& \quad  \quad \, \quad  \quad \quad  \quad  \quad  \, 
R(x, \, y, \, z) \, \, \,= \, \, \,\,
  \sum_m  \sum_n  \sum_l \,\, a_{m, \, n, \, l} \cdot x^m \, y^n \, z^l, 
\end{eqnarray}
by extracting the "diagonal" terms, i.e. the  powers of the product $\, p \, = \, \, x y z$:
\begin{eqnarray}
\label{diagmultiTaylorapp}
\hspace{-0.98in}&& \quad  \quad \, \quad  \quad \quad  \quad \quad  \, 
Diag\Bigl(R(x, \, y, \, z)\Bigr) \, \,\, = \, \, \,
 \sum_m \,\,  a_{m, \, m, \, m} \cdot p^m. 
\end{eqnarray}
Such diagonals are closely related to the integrals of rational
functions. For example $\, Diag\Bigl(R(x, \, y, \, z)\Bigr)$
is the constant term (in $\, y, \, z$) in the infinite expansion
\begin{eqnarray}
\label{multiTaylorscrewed}
\hspace{-0.98in}&& \quad   \, \quad  \quad \quad  \quad  \quad  \, 
R\Bigl({{p} \over {y \, z}}, \, y, \, z\Bigr) \, \, \,= \, \, \,\,
  \sum_{m, \, n, \, l \, \ge \,  0}  \,\, a_{m, \, n, \, l} \cdot p^m \, \, y^{n-m}\,  \, z^{l-m}, 
\end{eqnarray}
which can be represented by the integral~\cite{Lairez}
\begin{eqnarray}
\label{representedintegral}
\hspace{-0.98in}&& \quad  \quad \, \quad  \quad \quad  \quad  \quad  \, 
  {{1} \over {(2 \, \pi\, i)^2}} \, \oint \oint  \,
  R\Bigl({{p} \over {y \, z}}, \, y, \, z\Bigr) \, \, {{dy } \over {y}} \, \wedge {{dz } \over {z}}.
\end{eqnarray}
The  diagonal (\ref{diagmultiTaylorapp}) is also the constant term  (in $\, y, \, z$)
of
\begin{eqnarray}
\label{multiTaylorconstant}
\hspace{-0.98in}&&  \quad \, \quad  \quad \quad  \quad  \quad  \, 
R\Bigl({{p} \over {y }}, \, {{y} \over {z}}, \, z\Bigr) \, \, \,= \, \, \,\,
  \sum_{m, \, n, \, l \, \ge \,  0}   \,\, a_{m, \, n, \, l } \cdot p^m \, \, y^{n-m}\,  \, z^{l-n}, 
\end{eqnarray}
wich is of the form
\begin{eqnarray}
\label{multiTaylorconstant}
\hspace{-0.98in}&&  \quad \, \quad  \quad \quad  \quad  \quad  \, 
                   {{1} \over {(2 \, \pi\, i)^2}} \, \oint \oint  \, {{ N_p(y, \, z)} \over {D_p(y, \, z)}}
                   \, \, {{dy } \over {y}} \, \wedge {{dz } \over {z}},
\end{eqnarray}
where the numerator $\, N_p(y, \, z)$ and the denominator $\, D_p(y, \, z)$ are polynomials.
it is well-known that such integrals satisfy a linear differential equation with respect to $\, p$
having rational functions in $\, p$ as coefficients, called the Picard-Fuchs
equation\footnote[2]{The order of this linear differential equation is generally equal to the rank of the algebraic
 deRham cohomology of  $\, D_p(y, \, z) \, = \, 0$. For curves of genus $\, g$ this rank is $\, 2\, g$.}. 
the problem of determining such linear differential equations has been started by Griffiths~\cite{Griffiths}
with the assumption that the variety  $\, D_p(y, \, z) \, = \, 0$ is smooth, but later techniques were developed
to include examples with singular points~\cite{Lairez,Lairez2}. The linear differential equations
(Gauss-Manin systems, telescopers)
occuring in integrable models~\cite{LGF2,High, Khi6} are of order much larger than
order two\footnote[5]{Since Felix Klein it is well-known that
  the Picard-Fuchs equation corresponding to the (Weierstrass) elliptic curve corresponds to the hypergeometric function
$\, _2F_1([1/12,5/12],[1],1/J)$.} 
and almost never correspond to smooth varieties.  Creative telescoping\footnote[1]{For a detailed
  introduction to creative telescoping~\cite{Zeilberger} see for instance~\cite{Chyzak}.}
and more specifically the programs~\cite{Koutschan} corresponding to
a fast approach to creative telescoping~\cite{Koutschan3}, are a powerfull way to find these
linear differential operators annihilating these diagonal of rational functions in the cases
emerging naturally in theoretical physics, integrable models, enumerative combinatorics, for which the order
of the linear differential operators
is quite large~\cite{LGF2,High,Khi6}
and the  variety  $\, D_p(y, \, z) \, = \, 0$ is (most of the time) not a smooth one. All the pedagogical
(but non-trivial) examples of telescopers displayed in this paper can be viewed by an algebraic geometer
as a presentation of examples of families of varieties and their Picard-Fuchs equations. 

\section{A simple example corresponding to planar elliptic curves obtained as intersection of quadrics}
\label{Clebsch}

Let us consider the rational function  in three variables $\, x, \, y$ and $\, z$
\begin{eqnarray}
\label{ratioshift}
  \hspace{-0.98in}&& \quad   \quad   \quad   \quad   \quad   \quad   \quad   \quad  
\, \,  R(x, \, y, \, z)  \, \, \,  = \, \,  \, \, {{x \, y^2 } \over { D(x, \, y, \, z) }},
\end{eqnarray}
where:
\begin{eqnarray}
\label{ratioshiftD}
  \hspace{-0.98in}&& \quad \quad \quad \, \,\, 
D(x, \, y, \, z) \, \, \, = \, \, \, \, \,  4\,{x}^{4} \cdot \, xyz
 \, \, \, +16\,{y}^{2}{x}^{2} \, +16\,x{y}^{3} \, +16\,{y}^{4}
\, +32\,y{x}^{2} \, +40\,x{y}^{2}
\nonumber \\
 \hspace{-0.98in}&& \quad \quad  \quad   \quad   \quad \quad \quad \quad \quad
\, +40\,{y}^{3} \, +15\,{x}^{2}\, +25\,yx\, +41\,{y}^{2}\,\,  +40\,y\, +25, 
\end{eqnarray}
which corresponds (with $\, p \, = \, x\, y\, z$) to the {\em  elliptic curve}
\begin{eqnarray}
\label{ellipticshift}
  \hspace{-0.98in}&& \quad \, \quad  \quad \,
 C_p(x, \, y) \,\, = \, \,\, \,\, 
4\,{x}^{4} \cdot \, p \, \, \, \,
+16\,{y}^{2}{x}^{2} \, +16\,x{y}^{3} \, +16\,{y}^{4} \, +32\,y{x}^{2} \,
+40\,x{y}^{2}
\nonumber \\
\hspace{-0.98in}&& \quad \quad \quad \quad  \,  \quad \quad  \quad
\, \, \, +40\,{y}^{3} \, +15\,{x}^{2}\, +25\,yx\, +41\,{y}^{2}\, +40\,y \,\, +25
\, \, \, \, = \, \, \, \,  \, 0, 
\end{eqnarray}
corresponding to the intersection (elimination of $\, u$ at  $\, z \, = \, 1$)
of the two {\em quadrics} 
\begin{eqnarray}
\label{intersectionshift}
  \hspace{-0.98in}&&  \quad  \quad
p \cdot \, u^{2} \, +uz \, +x^{2} \, +y \, x \, +y^{2} \, +z^{2}      \, \, = \, \, \, 0,
\quad  \quad \, \,  4\, uy \, +5\,uz \, +2\, x^{2}   \, \, = \, \, \, 0.
\end{eqnarray}
The $\, j$-invariant of elliptic curve (\ref{ellipticshift}) reads:
\begin{eqnarray}
\label{J_nvshift}
  \hspace{-0.98in}&& \,\quad \quad \quad
J \,  \, \, = \, \,  \, \, {\frac {27  \cdot  \, (3523 \, +10496\,p)^{3}}{
6724 \cdot \, (2686976 \,  p^3 \,   - 1614336 \, p^2 \,  + 4051257 \, p \, - 470096 ) }}.        
\end{eqnarray}
The telescoper of the rational function (\ref{ratioshift}) is an order-three
linear differential operator which can be factorized as
\begin{eqnarray}
\label{L3shift}
\hspace{-0.98in}&& \quad  \quad \quad \quad  \quad \quad  \quad \quad 
 L_3 \, \,\, = \, \, \,\, L_2 \cdot \,
 \Bigl(D_x \, \, + {\frac {41\,x+8}{2 \, x \cdot \, \left( 41\,x-4 \right) }}  \Bigr), 
\end{eqnarray}   
where the order-two linear differential operator $\, L_2$ 
 is homomorphic to an order-two linear differential  operator $\, Z_2$ 
such that
\begin{eqnarray}
\label{Z2shift}
  \hspace{-0.98in}&&   \quad \quad \quad
L_2 \cdot \, \rho(x) \cdot \, X_1      \, = \, \, Y_1 \cdot \, Z_2
 \quad  \quad \quad \quad \quad \quad
\hbox{where:}
\\
  \hspace{-0.98in}&& \quad \quad \quad
\rho(x) \, \, = \, \, \,  {\frac {5}{16\, \, x \cdot \, (41\,x-4)  \, (6400\,x-11281) }}, 
\nonumber   \\
\hspace{-0.98in}&& \quad \quad \quad
 X_1 \, =  \, \, (839680\,{x}^{3} -16606384\,{x}^{2} -6835099\,x +2350480) \cdot \, D_x
\nonumber \\
\hspace{-0.98in}&&  \quad \quad \quad \quad \quad \quad
\, \, + \, 656 \cdot \, (2720\,{x}^{2}-9447\,x-2096),
\nonumber 
\\
\hspace{-0.98in}&& \quad \quad \quad
 Z_2 \, =  \, \, (2686976\,{x}^{3} -1614336\,{x}^{2} +4051257\,x -470096) \, (6400 \, x-11281) \cdot \, D_x^2
\nonumber \\
\hspace{-0.98in}&&  \quad \quad \quad \quad \quad \quad
\, \, + \, (34393292800\, x^3-101267079168\, x^2+36422648832\, x
 \\
\hspace{-0.98in}&&  \quad \quad \quad \quad \quad \quad
\, \, -42693615817) \cdot \, D_x + \, 1968 \cdot \, (1638400\, x^2-6531584\, x+79633),
\nonumber 
\end{eqnarray}
where the order-two operator $Z_2$ has the pullbacked $\, _2F_1$ hypergeometric solution
\begin{eqnarray}
\label{soluZ2shift}
  \hspace{-0.98in}&& \, \, \, \quad\quad\quad\quad\quad\quad
\Bigl(3523 \, +10496 \, x  \Bigr)^{-1/4} \cdot \, \,
 _2F_1\Bigl([{ {1} \over {12} }, \, {{5} \over {12}}], \, [1], \, \,
 \,   {\cal H}  \Bigr), 
\end{eqnarray}                
where the Hauptmodul $\, {\cal H} $ reads:
\begin{eqnarray}
\label{HauptZ2shift}
  \hspace{-0.98in}&& \quad \quad \quad \quad \quad \quad \quad \quad \quad 
{\cal H} \,  \, = \, \, \,
1 \, \, -  \,{\frac {27  \cdot \, (95457 -262400 \,x)^{2}}{ (3523 \, +10496\, \, x)^{3}}}, 
\end{eqnarray}              
which is nothing but the Hauptmodul associated to (\ref{J_nvshift})
(with, of course, $\, p$ changed into $\, x$).

\vskip .1cm 
\vskip .1cm 

\section{Maximum number of parameters for families of planar elliptic curves.}
\label{MaxParam}

We have seen, in section \ref{deducing}, that the previous results on diagonals of
nine or ten parameters families of
rational functions of three variables  being pullbacked $\, _2F_1$ hypergeometric
functions (and in fact classical  modular forms) can actually be seen as
corresponding to the (well-known in integrable models and integrable mappings)
fact that the most general biquadratic corresponding to  {\em elliptic curves}
is a {\em nine-parameters} family and that the most general ternary cubic
corresponding to elliptic curves is a {\em ten-parameters} family.
One can, for instance recall page 238 of~\cite{Silverman},
which amounts to considering
the collection of all cubic curves in $\, \mathbb{C} P_2$ with the homogeneous equation  
\begin{eqnarray}
\label{homo}
  \hspace{-0.98in}&& \quad   \quad \quad \quad \, \, \, \quad \quad
 a \, x^3 \,\, + b \, x^2\, y\, \, + \, c \, x \, y^2\, \, + \, d \,  y^3
  \, \,+ \,  e \, x^2 \, z \,\, +  f \, x \, z^2 \,\, +  g \, y^2\, z
 \nonumber \\
\hspace{-0.98in}&& \quad \quad \quad \quad \quad \quad \quad \quad  \quad \quad 
\, + h \, y \, z^2 \,\,\,  +i \, z^3 \,\,\,  + \, j \, x \, y \, z
 \,\, \,\, = \, \,\, \, 0,
\end{eqnarray}
and the associated problems of passing through nine given points. 
One can also recall the  ternary cubics in~\cite{Sadek,Poonen}
and other problems of elliptic curves of high rank~\cite{HighRank}
(see the concept of Neron-Severy rank). 

Since the rational functions of three variables we consider are essentially encoded
by the denominator of these rational functions, and in the cases we have considered,
the emergence of  pullbacked $\, _2F_1$ hypergeometric functions
(and in fact classical modular forms) corresponds to the fact that the intersection of
these denominators with the hyperbola $\, p \, = \, x \, y \, z$ corresponds to 
elliptic curves, one sees that these rational functions are essentially
classified by the possible $\, n$-parameters families $\, P(x, \, y) \, = \, \, 0 \, $
of elliptic curves.

If one considers a polynomial
\begin{eqnarray}
\label{homo}
\hspace{-0.98in}&& \quad   \quad \,   \quad \quad \quad  \quad \quad  \quad \, \, 
P(x,y)\, \,  =  \, \,\,  \,  \sum_m \, \sum_n \,  a_{m,n} \cdot  \, x^m \,  y^n, 
\end{eqnarray}
with generic coefficients $\, a_{m,n}\in \, \mathbb{C}$, then the genus
of the algebraic curve defined by~$P$ is determined by the support
  $\, supp(P) = \{(m,n) \in \mathbb{N}^2 : a_{m,n} \neq 0\}$.
More precisely, the genus equals the number of interior integer lattice points
inside the convex hull of $supp(P)$~\cite{Khovanskii} (see also the
discussion in~\cite{Khovanskii2}). For example, the support of the
ten-parameters family~\eqref{alg2bis} consists of the following $10$ points
in~$\mathbb{N}^2$:
\[
  (0,0), (0,1), (0,2), (0,3), (1,1), (1,2), (1,3), (2,2), (2,3), (3,3)
\]
which form a right triangle of side length~$3$. Only one of these points is an
interior point, namely $(1,2)$, hence the genus is~$1$.

Therefore we may ask: which integer lattice polytopes exist which have exactly
one interior point and what is the largest such polytope?  Not surprisingly,
the answer is known: there are (up to transformations like translation,
rotation, shearing) exactly $16$ different polytopes with a single interior
point~\cite{rabinowitz} (see also Figure 5, page 548 in~\cite{schicho}), 
the above-mentioned right triangle being the one with the
highest total number of lattice points.

This shows that {\em there cannot be a family of elliptic curves with more
than ten parameters}.

\vskip .1cm

\section{Monomial transformations preserving pullbacked hypergeometric results}
\label{monomialsec}

More generally, recalling subsection 4.2 in~\cite{unabri}
and  subsection 4.2 page 17 in~\cite{DiagJPA},
let us consider the monomial transformation
\begin{eqnarray}
\label{monomial}
 \hspace{-0.98in}&&  \,  \,  \, \quad \quad  \quad  \quad 
(x, \, y, \, z) \, \,\, \quad \,   \longrightarrow \,  \, \, \,\quad  \, 
M(x, \, y, \, z) \,\,\, = \, \,\,  (x_M, \, y_M, \, z_M) 
\nonumber \\
\hspace{-0.98in}&& \quad \quad \quad \quad \quad  \,  \,  \,\,  \,
\, \, \, \, = \, \, \, 
 \Bigl(x^{A_1} \cdot \, y^{A_2} \cdot \, z^{A_3}, 
\, \, \,   x^{B_1} \cdot \, y^{B_2} \cdot \, z^{B_3},
 \, \, \, x^{C_1} \cdot \, y^{C_2} \cdot \, z^{C_3} \Bigr), 
\end{eqnarray}
where the $\, A_i$'s,  $\, B_i$'s and   $\, C_i$'s are positive integers such that 
$\, A_1 \, = \, A_2 \, = \, A_3 \,$ is excluded (as well as 
 $\, B_1 \, = \, B_2 \, = \, B_3$ 
as well as  $\, C_1 \, = \, C_2 \, = \, C_3$),
and that the determinant of the $\, 3 \, \times \, 3$ matrix~\cite{DiagJPA,unabri}
\begin{eqnarray}
\label{3x3}
\hspace{-0.95in}&& \,\quad \,\quad \, \quad \quad \quad \,\quad \,\quad   \, 
\quad \quad \quad 
\left[ \begin {array}{ccc} 
                  A_1&B_1&C_1 \\ 
\noalign{\medskip}A_2&B_2&C_2 \\ 
\noalign{\medskip}A_3&B_3&C_3 
\end {array} \right],  
\end{eqnarray}
is not equal to zero\footnote[9]{We want the rational function 
$\, \tilde{{\cal R}} \, = \,\,  \,{\cal R}(M(x, \, y, \, z))$ 
deduced from  the monomial transformation (\ref{monomial}) 
to remain a rational function of {\em three} variables
and not of two, or one, variables. }, and that:
\begin{eqnarray}
\label{defdiagBIScond}
\hspace{-0.90in}&&\quad \quad    \quad   \quad  \, \, \,\,
\, A_1\, +B_1\, +C_1 \,\, = \, \, \,\, A_2\, +B_2\, +C_2 \, \, = \, \,\,  \, A_3\, +B_3\, +C_3.
\end{eqnarray}
We will denote by $\, n\, = \, A_i\, +B_i\, +C_i \, \, $ the integer\footnote[1]{Note a typo in
the footnote 28 page 17 of ~\cite{DiagJPA} as well as in the second footnote page 18 in~\cite{unabri}.
The sentence has been truncated. One should read:
For $\, n \, = \, 1$, 
the $\, 3 \times 3$ matrix (\ref{3x3}) is stochastic and  transformation (\ref{monomial}) 
is  a {\em birational  transformation} if the determinant of the
 matrix (\ref{3x3}) is $\, \pm \, 1$.}
in these three equal sums (\ref{defdiagBIScond}).
The condition (\ref{defdiagBIScond}) is introduced in order to impose that 
the product\footnote[5]{Recall that taking the diagonal of a rational function 
of three variables extracts, in the multi-Taylor expansion, only the terms 
that are $\, n$-th power of the {\em product}  $\, x\, y\, z$.} 
of $\, x_M\, y_M\, z_M$ is {\em an integer power of the product} of $\, x\, y\, z$:    
$\, \,x_M\, y_M\, z_M\, \, = \,\,\, (x\, y\, z)^n$.

If we take a  rational function  ${\cal R}(x, \, y, \, z)$ in  three variables 
and  perform such a monomial transformation (\ref{monomial}) 
$\, (x, \, y, \, z) \, \rightarrow \, M(x, \, y, \, z)$, 
on this rational function $ \,{\cal R}(x, \, y, \, z)$, we 
get  another rational function  that we  denote by
$\, \tilde{{\cal R}} \, = \,\,  \,{\cal R}(M(x, \, y, \, z))$. 
Now the diagonal of $\, \tilde{{\cal R}}$ 
is the diagonal of  $ \,{\cal R}(x, \, y, \, z)$
where we have changed $\, x$ into $\, x^n$:
\begin{eqnarray}
\label{defdiag2BIS}
\hspace{-0.7in}&&\, \, \, \, 
\Phi(x) \, \, = \, \, \, \, Diag\Bigl({\cal R}\Bigl(x, \, y, \, z \Bigr)\Bigr),
 \quad \quad \, \, \, 
 Diag\Bigl(\tilde{{\cal R}}\Bigl(x, \, y, \, z \Bigr)\Bigr)
\,  \, = \, \,  \,  \Phi(x^n).
\end{eqnarray}

\vskip .1cm 
\vskip .1cm

\section{Telescopers of  rational functions associated with a split Jacobian: the general case.}
\label{General}

Following calculations in subsection \ref{split}, we have seen that
for arbitrary values of the parameter $\, a$ and $\,\, b \, = \, \, 3 \, + \, x$, 
the Hauptmodul of one of the two elliptic curves of the split Jacobian reads
(see eq. (\ref{Hauptx})):
\begin{eqnarray}
\label{Hauptxapp}
  \hspace{-0.98in}&&   \, \quad    \quad   \quad   \quad    \quad   \quad    \quad    \quad   \quad  
{\cal H}_x \, \,\,\, = \, \, \,\, 
{\frac { 108 \cdot \, {x}^{3} \cdot \, ({x}^{2}+9\,x+27)^{3} \cdot \, P_2^{2}}{ P_4^{3}}}, 
\end{eqnarray}
where:
\begin{eqnarray}
\label{HauptxP2P4app}
  \hspace{-0.98in}&&  
P_2 \, = \, \,4\,{x}^{3}  \,  \, - (a-6)  \cdot \, (a+6)  \cdot \, {x}^{2}
\, \, \, -6 \cdot \, (a+6)  \cdot \, (a-3) \cdot \,  x \, \,\, + ( 4\,a \, +15)  \cdot (a \, -3)^{2},
\nonumber 
\\
\hspace{-0.98in}&&   
P_4 \, = \, \, 12 \,{x}^{5} \, \, \, +({a}^{2}+180) \cdot \,{x}^{4} \, \,\, 
 +6 \cdot \,  (2 \, {a}^{2}-21\,a+180) \cdot \,  {x}^{3}
\nonumber \\
\hspace{-0.98in}&&   \quad \quad \quad  \quad \quad  \,
+27 \cdot  \, (2\,{a}^{2} \, -42\,a \, +135)  \cdot \, {x}^{2}
\\
\hspace{-0.98in}&&   
 \quad  \quad \quad \quad  \quad \quad \quad \quad   \,
+162 \cdot \, (a \, -3)  \cdot \, (2\,a \, -15) \cdot \,  x \,\,\,\,
+729 \cdot \, (a \, -3)^{2}.\nonumber 
\end{eqnarray}

Let us consider the telescoper of the rational function of three variables
$\, x\, y/D_{a}(x,\, y, \, z)$ where the denominator
$\,\, D_{a}(x,\, y, \, z) \, $ is $\,\, C_{a,\, b}(x, \, y)\,$ given by  (\ref{genusTWO})
for $\, b \, = \, 3 \, + \, x\, y\, z$, namely (\ref{ratgenus2}).
We have calculated this telescoper for an arbitrary value of the parameter $\, a$.

This telescoper is an order-four linear differential operator $\, L_4$ which is actually the
{\em direct-sum} of {\em two order-two}
linear differential
operators $\, L_4 \, = \, \, LCLM(L_2, \, M_2)\, = \, \, L_2 \oplus \, M_2$, these two
order-two linear differential  operators having respectively the head polynomials $\, H(L_2)$ and  $\, H(M_2)$: 
\begin{eqnarray}
\label{head}
\hspace{-0.98in}&&\quad \,\,
H(L_2)  \, \, \,  = \, \, \, \,
4 \cdot \, x \cdot \, (x+3)^2 \cdot \,
\Bigl( a \cdot \, x^2 \, +\, 3 \cdot \, (2\, a \, -15) \cdot \, x \, +\, 3 \cdot \, (4\, a+15) \cdot \, (a-3) \Bigr)
\nonumber \\
\hspace{-0.98in}&&       \quad   \quad  \quad     \quad  \quad      \quad  \quad \quad \quad                 
\times  \, (x^2 \, +9\, x \, +27) \cdot \, P_2^2,
 \\
  \hspace{-0.98in}&& \quad \,\,
 H(M_2)  \, \,  \, = \, \, \,  \,
 4 \cdot \, \Bigl( a \cdot \, x \,\, + 3 \cdot \, (a \, -3)   \Bigr) \cdot \, P_2^2.                  
\end{eqnarray}
These two order-two linear differential operators {\em cannot be homomorphic} since they do not
have exactly\footnote[2]{They share, however, the singularities $\, P_2 \, = \, 0$.} the same singularities.
In view of the $\,\,  x^2 \, +9\, x \, +27\, \, $ term, we expect the order-two
linear differential operator $\, L_2 \, \, $
to have a  pullbacked $\, _2F_1$ hypergeometric solution with a Hauptmodul (\ref{Hauptxapp}),
or, at least, related to  (\ref{Hauptxapp}) by some isogeny (modular correspondence)~\cite{Klein,Klein2}.

Unfortunately the Maple command ''hypergeometricsols" of Mark van Hoeij~\cite{MarkvanHoeij} has
not been able to find the
pullbacked $\, _2F_1$ hypergeometric solution
of $\, L_2$ for values different from $\, a \, = \, 3$ and $\, a \, = \, \infty$.
For  $\, a \, = \, \infty \, $ we find the following  pullbacked $\, _2F_1$ hypergeometric
solution for $\, L_2$:
\begin{eqnarray}
\label{solainfty}
  \hspace{-0.98in}&&\quad \quad \quad \quad \quad 
{\cal S} \, \, = \, \, \,\,
\Bigr( {{ (x \, +\, 3)^3} \over { x^3 \, +9\, x^2 \, +27 \, x \, +3}}  \Bigl)^{1/4}
 \\
 \hspace{-0.98in}&& \quad \quad \quad \quad \quad \quad \quad \quad
\times \,
_2F_1\Bigl([{{1} \over {12}}, \, {{5} \over {12}}], \, [1], \,\,
 \,{\frac { 1728 \cdot \,  x \cdot \, (x^{2} \, +9\,x \, +27) }{
(x\, +3)^{3} \cdot \, (x^{3} \, +9\,{x}^{2} \, +27\,x \, +3)^{3} }}
 \Bigr).
\nonumber            
\end{eqnarray}
In the  $\, a \, = \, \infty \, $ limit the Hauptmodul (\ref{Hauptxapp})
reads:
\begin{eqnarray}
\label{Hauptainfty}
  \hspace{-0.98in}&&\quad \quad \quad \quad \quad \quad \quad \quad \quad
{\cal H} \, \,\,\, = \, \, \,\,\, \,
 {\frac {1728  \cdot \, {x}^{3} \cdot \, ({x}^{2}+9\,x+27)^{3}}{
 (x\, +9)^{3} \cdot \, (x \, +3)^{3} \cdot \, ({x}^{2} \, +27)^{3}}},             
\end{eqnarray}
to be compared with the pullback in the  pullbacked $\, _2F_1$ hypergeometric
solution  (\ref{solainfty}). This does not seem to match at first sight.
In fact, we have a remarkable identity: the   pullbacked $\, _2F_1$ hypergeometric
solution (\ref{solainfty}) {\em can also be written}:
\begin{eqnarray}
\label{solainfty2}
  \hspace{-0.98in}&&\quad \quad \quad \quad \quad
{\cal S} \, \, = \, \, \, \,
\Bigr( {{ (x \, +\, 3)^3} \over { x^2 \, +27) \cdot \, (x\, +9}}  \Bigl)^{1/4}
               \\
  \hspace{-0.98in}&& \quad \quad \quad \quad \quad \quad \quad \quad
                     \times \,
_2F_1\Bigl([{{1} \over {12}}, \, {{5} \over {12}}], \, [1], \,\, 
 \,{\frac {1728 \, \, {x}^{3} \cdot \, ({x}^{2}+9\,x+27)^{3}}{
 (9+x)^{3} \cdot \, (3+x)^{3} \cdot \, ({x}^{2}+27)^{3}}}
\Bigr).  \nonumber            
\end{eqnarray}
We have a {\em modular correspondence}  between these two Hauptmoduls appearing in (\ref{solainfty}) and
(\ref{solainfty2}). The algebraic relation between these two Hauptmoduls corresponds
to the {\em mmodular equation}:
\begin{eqnarray}
\label{modcorresp}
\hspace{-0.98in}&&
262144000000000\,{A}^{3}{B}^{3} \cdot \, (A+B) \,  \,
+4096000000 \, {A}^{2}{B}^{2}  \cdot \, (27\,{A}^{2}-45946\,AB+27\,{B}^{2})
\nonumber \\
\hspace{-0.98in}&& \quad  \quad \,
\,\, +15552000 \cdot \, AB  \cdot \, (A+B) \cdot \,  ({A}^{2}+241433\,AB+{B}^{2})
 \nonumber \\
\hspace{-0.98in}&& \quad  \quad \quad  \,\,\,
+729 \cdot \,({A}^{4}\, +B^4) \,  \,  -779997924 \cdot \, AB \cdot \, ({A}^{2}+{B}^{2})
\, \,  +1886592284694 \cdot  \,{A}^{2}{B}^{2}
\nonumber \\
\hspace{-0.98in}&&  \quad  \quad  \quad  \quad 
\, \,  +2811677184 \cdot \,A B \cdot \, (A+B) \, 
\, \, \, -2176782336 \cdot \,AB  \,  \, \, \, = \, \, \, \, 0,             
\end{eqnarray}
which is a representation of $\, \tau \rightarrow \, 3 \, \tau \, \, $ where $\, \tau$ is the ratio of the periods.

In fact, we have been able to find the two pullbacked hypergeometric solutions of $\, L_2$ and $\, M_2$.
One can actually discover that the pullbacked hypergeometric solutions of $\, L_2$ have the form
\begin{eqnarray}
\label{solL2}
\hspace{-0.98in}&& \quad  \quad  \quad  \,  \,  
{{ (x\, +3) } \over { P_2^{1/2} \cdot \, P_4^{1/4} }}  \, \cdot \, \,
_2F_1\Bigl([{{1} \over {12}}, \, {{5} \over {12}}], \, [1], \,\,\,
{\frac { 108 \cdot \, {x}^{3} \cdot \, ({x}^{2}+9\,x+27)^{3} \cdot \, P_2^{2}}{ P_4^{3}}}\Bigr),               
\end{eqnarray}
where the pullback in (\ref{solL2})  is {\em exactly the same Hauptmodul} $\, {\cal H}_x$ as  (\ref{Hauptxapp})
corresponding to the $\, j$-invariant of the elliptic curve of the split Jacobian of the genus-two curve !!
The pullbacked hypergeometric solutions of $\, M_2$ reads:
\begin{eqnarray}
\label{solM2}
  \hspace{-0.98in}&&
 \quad  \quad  \quad  \,  \,   \,  
{{ 1 } \over { P_2^{1/2} \cdot \, (a ^2 \, -9 \, -3 \, x)^{1/4} }}
\\
 \hspace{-0.98in}&& \quad  \quad  \,   \quad \quad  \quad \quad 
\, \times \, \,
_2F_1\Bigl([{{1} \over {12}}, \, {{5} \over {12}}], \, [1], \,\,\,
                    1 \, \,
- {{ \Bigl((a-3) \cdot \, (2\,a^2 \, +6 \, a \, -9) \, -9 \, a \, x\Bigr)^2 } \over { 4 \cdot \, ( a ^2 \, -9 \, -3 \, x )^3 }} \Bigr).
\nonumber 
\end{eqnarray}
Note that this Hauptmodul in (\ref{solM2})
\begin{eqnarray}
\label{Haupt2x}
  \hspace{-0.98in}&&
\quad  \quad  \quad  \quad \quad  \quad
 {\cal H}_2  \, \, \,\, = \, \, \,\,
                     1 \, \,  \,  \,
 - {{ \Bigl((a-3) \cdot \, (2\,a^2 \, +6 \, a \, -9) \, -9 \, a \, x\Bigr)^2 } \over { 4 \cdot \, ( a ^2 \, -9 \, -3 \, x )^3 }},          
\end{eqnarray}
is not expandable at $\, x \, = \, 0$.
That way we finally find the second  Hauptmodul corresponding to the second elliptic curve
in the split Jacobian $\, {\cal E}_1 \times {\cal E}_2$. Note that, in terms of $\, a$
and $\, b$ this second $\, j$-invariant has a quite simple form:
\begin{eqnarray}
\label{solM2j2}
  \hspace{-0.98in}&&
 \quad  \quad  \quad  \quad \quad  \quad \quad \quad \, \,
j_2 \, \, = \, \, \, \,
 {\frac { 256 \cdot \,  ({a}^{2} \, -3\,b)^{3}}{ {a}^{2}{b}^{2} \,\,  -4\,{a}^{3} -4\,{b}^{3} \,  \, +18\,ab \, -27}}.
\end{eqnarray}

Let us denote 
$\, A$ the first Hauptmodul  (\ref{Hauptxapp}) and $\, B$ this last  Hauptmodul  in (\ref{solM2}). For arbitrary values
of the parameter $\, a$ {\em they are not related by a modular correspondence}. The corresponding $\, j$-invariants
must be seen {\em as two independent}\footnote[1]{Emerging from the Igusa-Shiode invariants
  of the Jacobian.} {\em $\, j$-invariants}. Of course, eliminating $\, x$ one can find, for arbitrary 
values of the parameter $\, a$,
some quite involved (non-symmetric) polynomial relation $\, P(A, \, B) \, = \, \, 0\, $ between these two Hauptmoduls. 
We have, however,  rather simple relation for selected values of
the parameter $\, a$, namely:
\begin{eqnarray}
\label{RelAB}
  \hspace{-0.98in}&& \quad  \quad \quad \quad  \quad 
A \, \, = \, \, \,{\frac {27 \cdot \, {B}^{2} \cdot \,  (16\,{B}^{2}+4\,B+7)^{3}}{
 (4\,B-1)^{3} \cdot \, (4\,{B}^{2}+19\,B+4)^{3}}}, 
\quad \quad \quad \quad \hbox{for $\, a \, = \, \, 3$,}
\end{eqnarray} 
\begin{eqnarray}
\label{RelAB}
  \hspace{-0.98in}&& \quad \quad \quad \quad \quad \quad  \quad 
 A \, \, = \, \, -{\frac {{B}^{2} \cdot \, (4\,B-5)^{3}}{ (5\,B-4)^{3}}}, 
\quad\quad \quad  \quad \quad \hbox{for $\, a \, = \, \, 0$.}
\end{eqnarray}

Note that, in the   $\, a \, = \, \, \infty \, $ limit,
the second  Hauptmodul (\ref{Haupt2x})
trivialises and becomes $\, {\cal H}_2  \, \, \, = \, \, \,0$.
The   $\, a \, = \, \, \infty$ limit is, in fact, a bit tricky. The genus-two
curve (\ref{genusTWO}) becomes
\begin{eqnarray}
\label{G2infty}
  \hspace{-0.98in}&& \quad \quad \quad \quad \quad \quad  
 4 \cdot \, y^2 \,  \,   \,
 +\alpha^2 \, \cdot \, (b^2\,x^2 \, +4\,x^3 \,+2\,b\,x \, +1) \cdot \, x^2 \,   \, \, \, = \, \,   \, \, 0,               
\end{eqnarray}
where $\, \alpha \, = \, 1/a$,  $\, a \, \rightarrow  \, \, \infty$ (i.e.  $\, \alpha \rightarrow  \, \, 0$).
Let us introduce $\, Y \, = \, \, y/\alpha \, = \, a \cdot \, y$  (with $\, a \, \rightarrow  \, \, \infty$).
The previous curve (\ref{G2infty}) becomes the elliptic curve 
\begin{eqnarray}
\label{G2inftyY}
  \hspace{-0.98in}&& \quad \quad \quad \quad \quad \quad  \quad 
4 \cdot \, Y^2 \, \,  \,    + (b^2\,x^2 \, +4\,x^3 \,+2\,b\,x \, +1) \cdot \, x^2 \, \,   \, \, = \, \,   \, \, 0,               
\end{eqnarray}
with a $\, j$-invariant giving the Hauptmodul $\, {\cal H} \, = \, \, 1728/j$:
\begin{eqnarray}
\label{HauptG2inftyY}
\hspace{-0.98in}&& \quad \quad \quad \quad \quad \quad \quad \quad 
{\cal H} \, = \, \, \,
{\frac {  1728 \cdot \, (b-3)  \cdot \, ({b}^{2}+3\,b \, +9) }{{b}^{3} \cdot \, ({b}^{3}-24)^{3}}},   
\end{eqnarray}
which, for $\, b \, = \, \, 3 \, +\, x$, gives exactly the Hauptmodul in (\ref{solainfty2}), namely:
\begin{eqnarray}
\label{HauptG2inftyY}
\hspace{-0.98in}&& \quad \quad \quad \quad \quad \quad \quad \quad
{\cal H}_x \, = \, \, \,  \,
{\frac { 1728 \cdot \,  x \cdot \, (x^{2} \, +9\,x \, +27) }{
(x\, +3)^{3} \cdot \, (x^{3} \, +9\,{x}^{2} \, +27\,x \, +3)^{3} }}.
\end{eqnarray}
In that $\, a \, = \, \, \infty$ limit, the genus-two curve (\ref{genusTWO}) degenerates into
a {\em genus-one} curve. As far as creative telescoping is concerned, this amounts to calculating
the telescoper of the rational function 
\begin{eqnarray}
\label{HauptG2inftyrat}
  \hspace{-0.98in}&& \quad \quad \,\quad 
{{x \, y} \over {4 \cdot \, y^2 \, \,  \,  \,  + \, \Bigl((3 \, +x\, y\, z)^2 \cdot \,x^2 \, \,
+4\,x^3 \, +2 \cdot \, (3\, +x\, y\, z) \cdot \,x  \, \, +1\Bigr) \cdot \, x^2 }}.  
\end{eqnarray}
This telescoper  is an order-three linear differential $\, L_3$ which is the {\em direct-sum} 
 $\, L_3 \, = \, \, D_x \oplus \, {\cal L}_2$, where the order-two linear differential $ \, {\cal L}_2$
 is exactly the  $\, a  = \, \, \infty$ limit of the order-two linear differential operator $\, L_2$
 in the order-four  linear differential operator $\, L_4 \, = \, \, LCLM(L_2, \, M_2)\, = \, \, L_2 \oplus \, M_2$,
 with head polynomial (\ref{head}), and thus has (\ref{solainfty2}) as a solution. With these calculations
 we see, quite clearly, how the split Jacobian, which is
{\em isogenous to the product of two elliptic curves}, degenerates,
 in the $\, a \, = \, \, \infty$ limit, into an elliptic curve.

 Curiously, as a by-product of the calculation
 of the (non-symmetric) polynomial relation $\, P(A, \, B) \, = \, \, 0\, $ between these
 two Hauptmoduls, we find in the  $\, a \, = \, \, \infty$ case, a (spurious) genus-zero algebraic
 symmetric artefact relation 
 \begin{eqnarray}
 \label{RelAB}
   \hspace{-0.98in}&&
 1953125 \cdot \,{A}^{3}{B}^{3} \, \, -187500 \cdot \, {A}^{2}{B}^{2} \cdot \, (A \, +B) \, \,
  +375 \cdot \,AB \cdot \, (16\,{A}^{2}-4027\,AB+16\,{B}^{2})
  \nonumber \\
  \hspace{-0.98in}&&  \quad \quad \quad \quad 
 \, \, -64\, \cdot \, (A \, +B)  \cdot \, ({A}^{2} +1487\,AB +{B}^{2})
 \, \,  \,  +110592 \cdot \,AB \, \,  = \, \,  \, 0, 
 \end{eqnarray}
 which turns out to be the fundamental modular equation, parametrised by:
 \begin{eqnarray}
 \label{RelABfundx}
   \hspace{-0.98in}&& \quad  \,   \quad \quad \quad   \quad  \quad   \quad 
 A \, \, = \, \, \,{\frac { 1728\, z}{ (z+16)^{3}}},
\quad  \quad   \quad \, \, \, 
 B \, \, = \, \, \,{\frac { 1728\, {z}^{2}}{ \left( z+256 \right) ^{3}}}.
\end{eqnarray}

\vskip .1cm

\section{Creative telescoping on rational functions of three variables associated with genus-two curves
with split Jacobians: another example}
\label{splitgoodbisapp}

Let us consider the {\em genus-two curve} $\, C_p(x, \, y) \, = \, \, 0 \, $
given in equation (5) of lemma 4 of~\cite{Shaska2}: 
\begin{eqnarray}
\label{splitbis}
\hspace{-0.98in}&&   \quad  \quad  \quad \quad \quad \quad \quad \quad
 C_p(x, \, y) \, \,\, = \,\,  \,\, \, \,  x^6 \, +x^3 \, \, + p \, \, \, \,  -y^2.    
\end{eqnarray}
Let us introduce the rational function  $\, x\, y/D(x, \, y, \, z)$ where the denominator
$\,  D(x, \, y, \, z)\, \,\,$ is given by:
\begin{eqnarray}
\label{ratsplitbis}
\hspace{-0.98in}&&   \quad  \quad \quad \quad \quad
D(x, \, y, \, z) \, \, = \, \, \, \,  C_{p\, = \, xyz}(x, \, y) \, \, \, = \,\, \,\,  \, \,
x^6 \, +x^3 \, \,\, + x\, y\, z \, \,  \, \,   -y^2.       
\end{eqnarray}
The telescoper of this rational function is an order-two linear differential operator
which has the two hypergeometric solutions
\begin{eqnarray}
\label{solratsplit1bis}
  \hspace{-0.98in}&&   \quad \quad  \quad \quad  \quad  \quad \quad \quad \, \, 
 x^{-1/6} \cdot \,
 _2F_1\Bigl([{{1} \over {6}}, \, {{2} \over {3}}], \, [{{5} \over {6}}], \, \, 4\, x\Bigr),  
\end{eqnarray}
which is a Puiseux series at $\, \, x \, = \, \, 0 \, \, $ and: 
\begin{eqnarray}
\label{solratsplit2bis}
  \hspace{-0.98in}&&   \quad \, \, \quad  \quad \quad  \quad \quad \quad \quad
 x^{-1/6} \cdot \,  _2F_1\Bigl([{{1} \over {6}}, \, {{2} \over {3}}], \, [1], \,\, 1 \, -4\, x\Bigr).
\end{eqnarray}
These two  hypergeometric solutions can be rewritten
as\footnote[2]{The fact that $\, _2F_1\Bigl([{{1} \over {6}}, \, {{2} \over {3}}], \, [1], \, z\Bigr)$
can be rewritten as $\, _2F_1\Bigl([{{1} \over {12}}, \, {{5} \over {12}}], \, [1], \, H(z)\Bigr)$
where the Hauptmodul $\, H(z)$ is solution of a quadratic equation is given
in equation (H.16) of Appendix H of~\cite{malala}.} 
\begin{eqnarray}
\label{solratsplit3bis}
  \hspace{-0.98in}&&   \quad \quad \quad \quad \quad \quad \quad \quad 
 {\cal A}(x) \cdot \,  _2F_1\Bigl([{{1} \over {12}}, \, {{5} \over {12}}], \, [1], \, {{1728}  \over {J}} \Bigr).
\end{eqnarray}
where the $\,j$-invariant $\, J$, in the Hauptmodul $\,\, 1728/J \,$ 
{\em corresponds exactly to the degree-two elliptic subfields} of the the split Jacobian of the genus-two curve. 

Of course, if we change $\, p$ into $\, p \, \rightarrow  \, (1\, -p)/4$ in (\ref{splitbis}), the telescoper
of  the rational function  $\, x\, y/D(x, \, y, \, z) \, \, $ where the denominator
$\,  D(x, \, y, \, z)\, \, $ is given by
\begin{eqnarray}
\label{ratsplitbis}
\hspace{-0.98in}&&   \quad  \quad  \quad  
D(x, \, y, \, z) \, \, = \, \, \, \,  C_{p\, = \, xyz}(x, \, y) \, \, \, = \,\, \,  \, \,
x^6 \, +x^3 \, \, \, +\Bigl({{ 1 \, -x\, y\, z} \over {4}}\Bigr) \, \,  \,   -y^2,       
\end{eqnarray}
is the  order-two linear differential operator corresponding
to the $\, x \, \rightarrow  \, (1\, -x)/4$ pullback of the previous one. It has the two
hypergeometric solutions
\begin{eqnarray}
\label{solratsplit2ter}
  \hspace{-0.98in}&&   \quad \quad \quad \quad  \quad \quad  \quad \quad \quad \quad
(1\, -x)^{-1/6} \cdot \,
_2F_1\Bigl([{{1} \over {6}}, \, {{2} \over {3}}], \, [1], \,\, \, x\Bigr),
\end{eqnarray}
and:
\begin{eqnarray}
\label{solratsplit1bis}
 \hspace{-0.98in}&&   \quad \quad \quad  \quad \quad  \quad \quad \quad \quad \quad
 (1\, -x)^{-1/6} \cdot \,
_2F_1\Bigl([{{1} \over {6}}, \, {{2} \over {3}}], \, [{{5} \over {6}}], \, \, 1\, -\, x\Bigr).  
\end{eqnarray}
The pullbacked hypergeometric solution (\ref{solratsplit2ter}) can also be written
\begin{eqnarray}
 \label{forinstance2}
 \hspace{-0.98in}&&  \, \,  \,  \,  \, 
 (1\, -x)^{-1/6} \cdot \, _2F_1\Bigl([{{1} \over {6}}, \,{{2} \over {3}}], \, [1],   \, \, x\Bigr)
\, \,    = \, \, \,
  \, {\cal A}(x)
 \cdot \, _2F_1\Bigl([{{1} \over {12}}, \,{{5} \over {12}}], \, [1],  \, \, {\cal H}(x) 
 \Bigr), 
\end{eqnarray}
where $\, {\cal H}(x)$ reads
\begin{eqnarray}
 \label{forinstance2}
  \hspace{-0.98in}&&
{\cal H}(x) \, \, = \, \, \,
 4 \cdot \, x \cdot \, {\frac { 1458 \, -1215\,x \, + 125\,{x}^{2} }{( 25\,x \, -9)^{3}}}
  \,    \, + \, 8 \cdot \, x \cdot \,{\frac {(27-11\,x)  \cdot \, (27-25\,x) }{\sqrt {1-x} \cdot \, (9-25\,x)^{3}}}, 
\end{eqnarray}    
and where $\, {\cal A}(x) \, $ reads:
\begin{eqnarray}
 \label{forinstance3}
\hspace{-0.98in}&& \quad \quad \,  
{\cal A}(x)   \, \, = \, \, \, \,
(1\, -x)^{-5/24} \cdot \,  \Bigl( {{ 81} \over { 9 \, -25 \, x}}  \Bigr)^{1/8}
 \cdot \,  \Bigl( {{ 5 \cdot \, (1\, -x)^{1/2} \, -4 } \over {  5 \cdot \, (1\, -x)^{1/2} \, +4}}   \Bigr)^{1/8}.                 
\end{eqnarray}     

\vskip .2cm 

\section{Telescopers of rational functions of several variables: some examples}
\label{limitapp}

\vskip .2cm 

\subsection{Telescopers of rational functions of several variables: a second example with four variables}
\label{limitapp2}

Let us now consider the rational function in {\em four} variables $\, x, \, y, \, z, \, u$: 
\begin{eqnarray}
\label{Ratfonc4second}
  \hspace{-0.98in}&&  \,  \quad 
\, \,  
R(x, \, y, \, z, u)  \, \, \,  = \, \,  \quad 
  \\
\hspace{-0.98in}&&  
\quad   \quad  \quad  \, \,     \quad
 {{1} \over { 1 \,\, \, +9\, x \,\,+3\, y \,\,+z \,\,\, +9 \,y \,z \,\,\,  +3 \,u
 \,x \,\, \,\,\,\,+2\,x\,y
\, \,+5 \,x \,z \, \,\,  + 7 \,x^2\,y \,\,\,  + 11 \,z^2\,y}}.
\nonumber
\end{eqnarray}
The telescoper of this rational function of four variables is the {\em same order-two linear differential
operator} $\, L_2$ as for the telescoper of (\ref{Ratfonc4first}). It  has the
same pullbacked hypergeometric solution (\ref{telescRatfonc4u}).

Performing the intersection of the
codimension-one algebraic variety
$$\,\, 1 \,+9\, x \,+3\, y \,+z \, +3 \,u \,x \,+9 \,y \,z \,\,\, \,+2\,x\,y
\,+5 \,x \,z \, + 7 \,x^2\,y \,+ 11 \,z^2\,y\, \, = \, \, 0,$$
corresponding to the denominator of (\ref{Ratfonc4second}), with the hyperbolae
$\, p \, = \, \, x \, y \, z \, u \, $
amounts to eliminating, for instance $\, u$ (writing
$\, u\, = \, \, p/x/y/z$). This
gives $\, P_u \, = \, \, 0\, $ where $\, P_u$ reads:
\begin{eqnarray}
\label{Ratfonc4u}
\hspace{-0.98in}&& \, \,  \quad \quad  \quad \quad \quad 
P_u \,\, \, = \, \,\, \,
7\,{x}^{2}{y}^{2}z \,\,  \,+11\,{y}^{2}{z}^{3} \,\, \, +2\,x{y}^{2}z \,\, \, +5\,xy{z}^{2}
\,\, \, +9\,{y}^{2}{z}^{2}
\nonumber \\
\hspace{-0.98in}&& \quad  \quad \quad \quad \quad \quad \quad  \quad  \quad  \quad 
  \, +9\,xyz \,+3\,{y}^{2}z \, \, +y{z}^{2}  \, \, +yz \,\,   \, +3\,p.
\end{eqnarray}
Assuming $\, x$ to be constant\footnote[1]{If one assumes $\, z$ to be constant, the previous condition $\, P_u(y, \, z) \, =  \, 0 \, $
becomes a genus-zero curve. If one assumes $\, y$  to be constant,  the previous condition $\, P_u(y, \, z) \, =  \, 0 \, $
is again a genus-one curve, but the corresponding Hauptmodul, which depends on $\, y$
is not simply related to (\ref{HauptvRatfonc4uDD}) for any selected value of $\, y$.}
the previous condition $\, P_u(y, \, z) \, =  \, 0 \, $
is an algebraic curve. Calculating its genus, one finds immediately
that it is  {\em genus-one}.
Calculating its $\, j$-invariant, one finds 
\begin{eqnarray}
\label{JinvRatfonc4u}
  \hspace{-0.98in}&& \, \,   \quad \quad \quad
 J \, \, = \, \, \,
-\, {{N^3} \over { 27 \cdot \, p^3 \cdot \, (7x^2 +2x +3)^2 \cdot D}}
 \quad \quad \quad \quad \quad \quad  \hbox{where:}
 \\
\hspace{-0.98in}&& \, \quad \quad  \,  \, \, \, \, \,
N \, = \, \, \,
 81 \cdot \, (280\,p+81) \cdot \, {x}^{4} \, +36\, (376\,p+81)  \cdot  \,{x}^{3}
                   \nonumber \\
\hspace{-0.98in}&& \quad \quad \quad  \quad \quad  \quad \quad \quad \quad \,\,\, 
   -18 \cdot \, (1848\,{p}^{2} +292\,p-27) \cdot \, {x}^{2}
 \nonumber \\
\hspace{-0.98in}&& \quad \quad \quad  \quad \quad \quad \quad  \quad  \quad \quad \quad  \,
-36 \cdot \, (264\,{p}^{2}+20\,p-1)  \cdot \, x
\, \,\,\, -2592\,{p}^{2} \, \,\, +1,
\nonumber               
\end{eqnarray}
and where:
\begin{eqnarray}
  \label{JinvRatfonc4uDD}
  \hspace{-0.98in}&&  \quad  \quad  \quad\quad \, \,
  D \, \, = \, \, \, \,   70875 \cdot \, (35\,p +9)  \cdot \, {x}^{8} \,\,
     +97200 \cdot \, (35\,p+8)  \cdot \, {x}^{7}
\nonumber \\
   \hspace{-0.98in}&&  \quad \, \, \quad \quad  \quad       \quad  \quad             \,
  +270 \cdot  \, (38655\,p+9047) \cdot \, {x}^{6}
 \, \,   +4 \cdot \, (2000295\,p+393278) \cdot \, {x}^{5}
\nonumber \\
 \hspace{-0.98in}&&  \quad \, \,  \quad \quad   \quad  \quad  \quad  \quad \,
 - (21704760\,{p}^{2}+1329219\,p-446680) \cdot \, {x}^{4}
\nonumber \\
\hspace{-0.98in}&& \quad \quad  \quad \quad \quad \quad \quad\quad  \quad
 \, -36 \cdot  \, \left( 332112\,{p}^{2}+28965\,p-1888 \right) \cdot \, {x}^{3}
 \\
  \hspace{-0.98in}&& \quad \quad  \quad \quad \quad \quad \quad\quad  \quad
\, +2 \cdot  \,(8049888\,{p}^{3}-864324\,{p}^{2}-70038\,p+2903) \cdot \, {x}^{2}
\nonumber \\
  \hspace{-0.98in}&& \quad \quad  \quad \quad \quad \quad \quad \quad\quad  \quad
 +24 \cdot  \, (191664\,{p}^{3} \, -16218\,{p}^{2} -246\,p+11) \cdot \, x
 \nonumber \\
 \hspace{-0.98in}&& \quad  \quad\quad  \quad\quad \quad  \quad \quad \quad \quad \quad \quad\quad  \quad
\,\,  +2665872\,{p}^{3}\,-19440\,{p}^{2}\,-12\,p \,\,\, +5.
\nonumber
\end{eqnarray}
In the $\, \, x \, \rightarrow \, \, 0 \,\,  \, $ limit of  the 
Hauptmodul $\, \, H_{p,x} \, = \, \, 1728/J$,  one finds:
\begin{eqnarray}
\label{HauptvRatfonc4uDD}
\hspace{-0.98in}&& \, \quad  \quad \quad  \quad \quad  \, \, 
  H_p \, \, = \, \, \,\,
- \, {{  419904 \cdot \, p^3 \cdot \, (5 \, -12\, p \, -19440\, p^2 \, +2665872\, p^3) )
} \over { (1 \, - \, 2592\, p^2)^3}}. 
\end{eqnarray}
{\em which actually corresponds to the Hauptmodul in} (\ref{telescRatfonc4u}). 

\vskip .2cm

\subsection{Telescopers of rational functions of several variables: a third example with four variables}
\label{limitapp1}

Let us consider the rational function in {\em four} variables $\, x, \, y, \, z, \, u$: 
\begin{eqnarray}
\label{Ratfonc4third}
  \hspace{-0.98in}&&  \,  \quad 
\, \,  
R(x, \, y, \, z, \, u)  \, \, \,  = \, \,  \quad 
  \\
\hspace{-0.98in}&&  
\quad   \quad  \, \,  \, \,  \, \,     \quad
 {{1} \over { 1 \,\, \,\,+3\, y \,\,+z \,\, \,\,+9 \,y \,z \,
\, \, \,  + 11 \,z^2\,y  \,  \, \, +3 \,u \,x \,\, \,\,  \, + x \cdot \, P_1(y, \, z)  \,  \, + x^2 \cdot \, P_2(y, \, z) }},
\nonumber
\end{eqnarray}
where $\, P_1(y, \, z)$ and  $\, P_2(y, \, z)$ are the two simple polynomials 
$\, P_1(y, \, z) \,  = \, \, y^2\, z^2 $
and  $\, P_2(y, \, z) \,  = \, \, y^3$. 
The telescoper of this rational function of {\em four variables} is the
{\em same order-two linear differential
operator} $\, L_2$ as for the telescoper of (\ref{Ratfonc4first}). It  has the
same pullbacked hypergeometric solution (\ref{telescRatfonc4u}). Actually the diagonal
of  the rational function (\ref{Ratfonc4first}) is the expansion (\ref{telescRatfonc4uexp})
of the pullbacked hypergeometric function (\ref{telescRatfonc4u}).
For $\, P_1(y, \, z) \,  = \, \, y^2\, z^2\, $ and  $\, P_2(y, \, z) \,  = \, \, y^3\, $
the elimination of $\, u\, = \, p/x/y/z \, $ in the vanishing condition of the denominator
(\ref{Ratfonc4third}) gives the algebraic curve:
\begin{eqnarray}
\label{genus5}
  \hspace{-0.98in}&&  \, \,  \, \quad
 x^2\, y^4 \, z \,  \, +x\, y^3\, z^3 \,  \, +11\, y^2\, z^3 \, \,  +9 \, y^2\, z^2 \,
 +3\, y^2 \, z \, \,  +y\, z^2 \, \,  \, +y \, z \, \,  + 3 \, p \, \, = \, \, \, 0.                    
\end{eqnarray}
For $\, x$ fixed (and of course $\, p$ fixed) this algebraic curve (\ref{genus5})
is a {\em genus-five} curve, but in the $\, x\, \rightarrow \, 0$ limit it
reduces to the {\em same genus-one} curve as for the first example (\ref{Ratfonc4first}),
namely:
\begin{eqnarray}
\label{genus5one}
  \hspace{-0.98in}&&  \,  \quad \quad \quad \quad \quad\quad
  11\, y^2\, z^3 \, \, +9 \, y^2\, z^2 \,\, 
 +3\, y^2 \, z \, \, +y\, z^2 \,\,  \, +y \, z \, \,\,  + 3 \, p \, \, \, \, = \,\,  \, \, 0.                    
\end{eqnarray}
which corresponds to the Hauptmodul (\ref{HauptvRatfonc4uDD}).

The generalisation of this result is straightforward. Let us consider the rational
function in {\em four variables} $\, x,\, y, \, z$ and $\, u$
\begin{eqnarray}
\label{Ratfonc4thirdgeneral}
  \hspace{-0.98in}&&  \,  \quad  \quad  \quad 
\, \,  
R(x, \, y, \, z, \, u)  \, \, \,  = \, \,  \quad 
  \\
\hspace{-0.98in}&&  
\quad   \quad  \quad  \, \,    \quad   \quad  \quad \quad
 {{1} \over { 1 \,\, \,\,+3\, y \,\,+z \,\, \,\,+9 \,y \,z \,
\, \, \,  + 11 \,z^2\,y  \, \,  +3 \,u \,x \, \, \,  +  x \cdot \, P(x, \, y, \, z) }},
\nonumber
\end{eqnarray}
where $\, P(x, \, y, \, z)$ is an {\em arbitrary polynomial}
of the three variables $\, x$, $\, y$ and $\, z$.
On a large set of examples one verifies that the {\em diagonal}
of (\ref{Ratfonc4thirdgeneral}) is actually the expansion (\ref{telescRatfonc4uexp})
of the pullbacked hypergeometric function (\ref{telescRatfonc4u}):
\begin{eqnarray}
\label{telescRatfonc4uexpbis}
  \hspace{-0.98in}&& \quad \, \, \, \,
 1 \, \, +648\,{x}^{2} \, \, -72900\,{x}^{3} \, \, +1224720\,{x}^{4} \,
  -330674400\,{x}^{5} \, +23370413220\,{x}^{6} \,
  \\
  \hspace{-0.98in}&& \quad \quad \, \,           \,
\, -1276733858400\,{x}^{7} \, +180019474034400\,{x}^{8}
 \, -12013427240614800\,{x}^{9} \, \, \, + \, \, \, \cdots \nonumber 
\end{eqnarray}
However, as far as creative telescoping calculations are concerned
the telescoper corresponding to different polynomials $\, P(x, \, y, \, z)$
becomes quickly a quite large non-minimal linear differential operator.
For instance, even for the simple polynomial $\, P(x, \, y, \, z) \, = \, \, x \, +y$, one obtains a quite large
order-ten telescoper. Of course, since this telescoper has the
pullbacked hypergeometric function (\ref{telescRatfonc4u}) as a solution, it is not minimal, it
is rightdivisible by the order-two linear differential operator having
(\ref{telescRatfonc4u}) as a solution.
It is straightforward to see that the previous elimination of $\, u\, = \, p/x/y/z \, $
in the vanishing condition of the denominator
(\ref{Ratfonc4thirdgeneral}) gives an algebraic curve\footnote[1]{Of arbitrary large genus for increasing degrees of the
 polynomial $\, P(x, \, y, \, z)$.}
\begin{eqnarray}
\label{genusM}
  \hspace{-0.98in}&&  \, \,  \, \quad
 11\, y^2\, z^3 \, \,  +9 \, y^2\, z^2 \,
 +3\, y^2 \, z \, \,  +y\, z^2 \, \,  \, +y \, z \, \,  + 3 \, p \, \,
 + y \, z \cdot \, P(x, \, y, \, z) = \, \, \, 0.                    
\end{eqnarray}
which reduces again, in the $\, x\, \rightarrow \, 0$ limit,
to the {\em same genus-one} curve (\ref{genus5one}).

With that general example (\ref{Ratfonc4thirdgeneral}) we see that there is 
{\em an infinite set of rational functions depending
  on an arbitrary polynomials $\,  P(x, \, y, \, z)$ of three variables} which diagonals are 
 a pullbacked $\, _2F_1$ hypergeometric solution. 

 \vskip .2cm

\subsection{Telescopers of rational functions of several variables: some examples}
\label{limitapp}

More generally we find that the diagonal of the rational function in $\, x, \, y, \, z, \, u$
\begin{eqnarray}
\label{Ratfonc4G}
  \hspace{-0.98in}&&
R(x, \, y, \, z, \, u)  \, \, \,  = \, \,  \quad 
  \\
\hspace{-0.98in}&&  
  {{1} \over { a \,\, \,+ b_1 \, y \,  + \, c_1/y  + b_2 \, z  \, + \, c_2 /z \,  + d_1 \, y \, z  \,
  + \,e_1 \, y/z \,  + \,f_1 \, z/y \,  \, +g_1 \,u \,x \, \,  + x^N \cdot P(y, \, z)}},
\nonumber
\end{eqnarray}
gives, for every integer $\, N \ge \, 1$,
a telescoper {\em independent of the arbitrary polynomial}  $\, P(y, \, z)$, namely 
the {\em same} telescoper that the rational function of $\, x, \, y, \, z, \, u$
\begin{eqnarray}
\label{Ratfonc4G}
  \hspace{-0.98in}&& \quad   \quad   
R(x, \, y, \, z, \, u)  \, \, \,  = \, \,  \quad 
  \\
\hspace{-0.98in}&&  
 \quad   \quad   \quad  \quad   \quad   
{{1} \over { a \,\, \,+ b_1 \, y \,  + \, c_1/y  + b_2 \, z  \, + \, c_2 /z \,  + d_1 \, y \, z  \,
  + \,e_1 \, y/z \,  + \,f_1 \, z/y \,  \,+g_1 \,u \,x }}.
\nonumber
\end{eqnarray}
The telescoper annihilates the pullbacked hypergeometric function:
\begin{eqnarray}
\label{solRatfonc4GG}
  \hspace{-0.98in}&& \quad   \quad  \quad  \quad  \quad    \quad \quad     \quad \quad    
D_H^{-1/12} \cdot \,
 _2F_1\Bigl( [ {{1} \over {12}}, \, {{5} \over {12}}], [1], \, {{1728} \over {J}}\Bigr), 
\end{eqnarray}
where $\, D_H$ denotes the denominator of the Hauptmodul $\, H \, = \, 1728/J$
and where the j-invariant $\, J$ is the  j-invariant of the elliptic curve
corresponding to the $\, x \, = \, 0$ limit of 
\begin{eqnarray}
\label{solRatfonc4GG}
\hspace{-0.98in}&& \quad \quad    \quad \quad   \quad    \quad \quad \quad \quad   \quad   
D\Bigl(x, \, y, \, z, \, {{p} \over {x\, y\, z}}\Bigr) \, \, \, = \, \, \, \,  0,            
\end{eqnarray}
namely the most general {\em nine-parameters} biquadratic $\, B(y, \, z) = \, 0$:
\begin{eqnarray}
\label{biquadr}
  \hspace{-0.98in}&&  
 d_1 \, y^2\, z^2  \, +b_1 \, y^2\, z \, +b_2 \, y\, z^2 \,
 + a \, y\, z \, +e_1\, y^2 \, +f_1\, z^2 \, \, +c_1\, z \, +c_2\, y
   \, + g_1 \, p   \, \, = \, \, \, 0.
\end{eqnarray}

\vskip .2cm

\subsection{A simple $\, u$-extension of the bicubic case.}
\label{extension}

Let us perform a similar ``reversed engineering" with  the (selected) ten-parameters bicubics
like (\ref{alg2bis}) that are elliptic curves. 
Let us consider one of these  (selected) bicubic  $ \, B(x, \, y)$: 
\begin{eqnarray}
\label{alg2bisextselected}
  \hspace{-0.98in}&&
        \quad \quad  \quad \,   \quad \quad \, \,
B(x, \, y)  \,  \, \, = \, \, \,\,
2\,{x}^{3}{y}^{3} \, +5\,{x}^{2}{y}^{3} \, +3\,{x}^{2}{y}^{2} \,
 \nonumber \\
 \hspace{-0.98in}&&  \quad  \quad \quad \quad \quad \quad \quad \quad \quad \quad \quad \, \,\,
 +x{y}^{3} \, +x{y}^{2} \, +3\,{y}^{3} \, +xy \, +3\,{y}^{2} \, +2\,y \, +5
\end{eqnarray}
The bicubic equation  $\, B(x, \, y) \, + \, p = \, \, 0 \, \, $  {\em is an elliptic curve}. One
can calculate its $\, j$-invariant
and the corresponding Hauptmodul $\, 1728/j$:
\begin{eqnarray}
\label{alg2bisextselected}
  \hspace{-0.98in}&&  \quad \quad  \quad \quad \, \,
{\cal H} \, \, = \, \, \,
 -1728\,{\frac {p_4(p) }{ \left( 2424\,p+11305 \right) ^{3}}},
 \quad \quad  \quad \quad\quad \hbox{where:}
  \\
\hspace{-0.98in}&&  \, \, \quad \quad \quad \quad
p_4(p) \,\, = \, \,\,
 99015075\,{p}^{4}\,+1743092117\,{p}^{3}\, +11512110810\,{p}^{2} \\
  \hspace{-0.98in}&& \quad  \quad \quad \quad  \quad \quad \quad  \quad \quad \, \,\,
 +33804556190\,p \,\,+37237506697.
 \nonumber 
\end{eqnarray}
Let us now consider the rational function of four variables
\begin{eqnarray}
  \label{alg2bisextselectedRat}
  \hspace{-0.98in}&&  \quad \quad  \quad \quad \quad \quad   \quad \quad \, \,
R(x, \, y, \, z, \, u) \, \, = \, \, \, \,
{{x \, y^2 } \over { B(x, \, y) \, \, + \, \, u \, x\, y \, z}}. 
\end{eqnarray}
Its telescoper is an order-two linear differential operator $\, L_2$
with pullbacked $\, _2F_1$ hypergeometric solutions
\begin{eqnarray}
  \label{alg2bisextselectedsol}
  \hspace{-0.98in}&&  \quad \,  \, \,
 (2424 \,x \, +11305)^{-1/4} \cdot \,
 _2F_1\Bigl([{{1} \over {12}}, \,  {{5} \over {12}}], \, [1], \,
-1728\,{\frac { p_4(x) }{ (2424\,x+11305)^{3}}}\Bigr),
\nonumber \\
 \hspace{-0.98in}&&  \, \,  \hbox{where:} \quad \quad \quad
p_4(x) \,\, = \, \,\,
 99015075\,{x}^{4} \, +1743092117\,{x}^{3} \, +11512110810\,{x}^{2}
  \\
  \hspace{-0.98in}&& \quad \quad  \quad \quad \quad \quad \quad \quad \quad
 \, \,\, +33804556190\,x \, +37237506697,
\nonumber                
\end{eqnarray}
which actually corresponds to the  Hauptmodul (\ref{alg2bisextselected})
of the bicubic equation  $\, B(x, \, y) \, + \, p = \, \, 0$.

\vskip .2cm

{\bf Comment:} if one considers, instead of (\ref{alg2bisextselectedRat}),  the rational function
\begin{eqnarray}
  \label{alg2bisextselectedRatxy}
  \hspace{-0.98in}&&  \quad \quad  \, \quad \quad \quad \quad \quad \quad \, \,
 R(x, \, y, \, z, \, u) \, \, = \, \, \, \, {{x \, y} \over { B(x, \, y) \, \, + \, \, u \, x\, y \, z}}, 
\end{eqnarray}
one finds an  {\em order-four} telescoper which factorises into {\em two order-two}
linear differential operators $\, M_4 \, = \, \, M_2 \cdot N_2$,
where $\, N_2$ has algebraic functions solutions, 
and where $\, M_2$ is homomorphic to the previous order-two linear differential operator $\, L_2$.

\vskip .2cm
\vskip .1cm

\subsection{Another simple $\, u$-extension of the bicubic case.}
\label{extension2}

Let us consider the rational function in three variables
\begin{eqnarray}
\label{alg2bisextA}
\hspace{-0.98in}&& \quad \quad \quad 
R(x, \, y, \, z) \, = \, \,
 \\
\hspace{-0.98in}&& \quad \quad \quad \quad   \quad \quad 
 {{1} \over {1 +x +2\,y \, +3\,z \,  +2\,y\, z \, +5\, x\, z \,  +7\, x\, z \,  +x^2\, y +y^2\, z +2\, z^2\,x }}.
 \nonumber                 
\end{eqnarray}
If one substitute  $\, z \, = \, p/x/y \, $ in the rational function (\ref{alg2bisextA}), one gets
\begin{eqnarray}
\label{alg2bisextAnewrat}
  \hspace{-0.98in}&&  \quad \quad \quad \quad \quad \, \,
 {\cal R}_p(x, \, y) \, \, = \, \, \,   {{  x\, y^2} \over { {\cal B}_p(x, \, y)}}
\quad  \quad \quad \quad  \quad  \quad \hbox{where:}
\\
\hspace{-0.98in}&&   \quad
{\cal B}_p(x, \, y)   \, \, = \, \,   \, \, \,
2\,{p}^{2} \, \, +y \cdot \, ({y}^{2}+12\,x+2\,y+3) \cdot \, p \, \, \,+{x}^{3}{y}^{3}\,
 +{x}^{2}{y}^{2}\,+2\,x{y}^{3}\,+x{y}^{2}.
\nonumber        
\end{eqnarray}
Let us now consider the rational function in {\em four} variables $\, x, \, y, \, z$ and $\, u$
which is  $\, {\cal R}_p(x, \, y)$ given by (\ref{alg2bisextAnewrat}) where  $\, p \, = \, \, x\, y\, z\, u$:
\begin{eqnarray}
\label{alg2bisextAnewratB}
 \hspace{-0.98in}&& 
 {{  x\, y^2} \over {
 2\, (x\, y\, z\, u)^{2} \, +y \cdot \, ({y}^{2}+12\,x+2\,y+3) \cdot \,  x\, y\, z\, u \, \,
 +{x}^{3}{y}^{3}\, +{x}^{2}{y}^{2}\,+2\,x{y}^{3}\, +x{y}^{2}
 }}.  
\end{eqnarray}
The telescoper of the rational function in three variables (\ref{alg2bisextA}), and
the telescoper of the rational function of four variables (\ref{alg2bisextAnewratB}),
{\em are actually equal}, having the pullbacked $\, _2F_1$ hypergeometric solution
given by  (\ref{2F15HypformAplusplusplus})
in subsection \ref{ten}.  In this particular case the pullbacked $\, _2F_1$
hypergeometric solution reads
\begin{eqnarray}
\label{particular}
  \hspace{-0.98in}&&  \quad \quad   \quad  \, \, \, \,
 (1 -64\,x+7552\,{x}^{2}+3600\,{x}^{3})^{-1/4}
 \nonumber \\
 \hspace{-0.98in}&&  \, \,  \, \quad \quad \quad \quad \quad 
           \quad           \times \,
_2F_1\Bigl([{{1} \over {12}}, \,  {{5} \over {12}}], \, [1], \,
 {{ 1728 \cdot \, x^3 \cdot \, p_9(x)
 } \over {(1-64\,x+7552\,{x}^{2}+3600\,{x}^{3})^3 }}  \Bigr),                   
\end{eqnarray}
where:
\begin{eqnarray}
\label{particularwhere}
  \hspace{-0.98in}&& \quad 
\, p_9(x) \, \, = \, \, \, \,675 \,\, -46908\,x+7579422\,{x}^{2} \, \,
-256103188\,{x}^{3}\,\, +748623104\,{x}^{4}
 \\
 \hspace{-0.98in}&&  \quad \quad  \, \,
-1361870768\,{x}^{5}
 +554260968\,{x}^{6}-1071752256\,{x}^{7}-36904896\,{x}^{8}-314928\,{x}^{9}.
 \nonumber
\end{eqnarray}
The Hauptmodul of the elliptic curve  $\, {\cal B}_p(x, \, y)   \, \, = \, \, 0$
corresponds to the pullback in the pullbacked $\, _2F_1$
hypergeometric solution (\ref{particular}):
\begin{eqnarray}
\label{particularwhereHaupt}
  \hspace{-0.98in}&& \quad \quad \quad \quad \quad \quad \quad \quad \quad \quad
  {\cal H} \, \, = \, \, \,   {{ 1728 \cdot \, p^3 \cdot \, p_9(p)
 } \over {(1-64\,p +7552\,{p}^{2}+3600\,{p}^{3})^3 }}.                 
\end{eqnarray}

\vskip .2cm 

\subsection{One more simple $\, u$-extension of the bicubic case.}
\label{extensionbeurk}

Let us perform a similar ``reversed engineering" with  the (selected) ten-parameters bicubics
like (\ref{alg2bis}) that are elliptic curves. 
Let us recall the fact that the (selected) bicubic $\, B(x, \, y) \, = \, \, 0 \, $
where $ \, B(x, \, y)$ reads: 
\begin{eqnarray}
\label{alg2bisext}
\hspace{-0.98in}&&  \quad \quad \, \,
B(x, \, y)  \,  \, \, = \, \, \,\,
 a \,x \, y^2 \, \,  + \, b_1 \, x^2 \, y^2 \,  + \, b_2 \, x \, y^3 \,   + \, b_3 \,  \, y \,\,\, 
+ \, c_1 \,  \, y^2 \, + \, c_2 \, \, p x \, y \, + \, c_3 \, x^2 \, y^3
\nonumber \\
\hspace{-0.98in}&& \quad \quad \quad \quad \quad \quad \quad \quad 
 \, \, \, + \, \, d_1 \, x^3 \, y^3 \, \,   + \, \, d_2 \, \, y^3
 \,   + \, \, d_3. 
\end{eqnarray}
is an elliptic curve.  Let us consider the denominator of four variables
$\, x$,  $\, y$,  $\, z$ and $\, u$:  
\begin{eqnarray}
\label{alg2bisextD}
\hspace{-0.98in}&&  \quad \quad \, \,\, \quad \quad \,   \quad \quad  \quad
 D(x, y, z, u))  \, \,  \,  =  \, \, \, \,\,
  B(x, \, y) \,\,  \,+ \,  u \cdot x \, y \,  z.               
\end{eqnarray}
It is straightforward to see that, imposing that the
product of the four variables $\, x\, y\, z\, u \, = \, \, p$ is fixed,
the elimination of the fourth variable  $\,u$, by the substitution
$\, u \, = \, p/x/y/z$,  yields a bicubic  $\, B(x, \, y) \,  \,+p \, \,  = \, \, 0$,
which is an elliptic curve. Introducing the rational function of four variables
\begin{eqnarray}
\label{alg2bisextrat}
\hspace{-0.98in}&&  \quad \quad \quad \quad  \quad \quad  \quad \quad\quad \quad
  R(x, y, z, u))  \,  \, \, =  \, \, \, \,
     {{ x \, y^2 } \over { D(x, y, z, u))  }},               
\end{eqnarray}
one finds that the telescoper of this rational function of four variables
(\ref{alg2bisextrat}) is an order-two linear differential operator
\begin{eqnarray}
\label{alg2bisextratAx}
  \hspace{-0.98in}&&  \quad \quad  \quad \quad \quad  \quad\quad \quad\quad \quad
 {\cal A}(x) \cdot \, \, 
_2F_1\Bigl([{{1} \over {12}}, \, {{5} \over {12}} ],  [1], \, \, {\cal H}   \Bigr)                     
\end{eqnarray}
where the Hauptmodul $\,\, {\cal H} \, = \, 1728/j \,$ corresponds
to the $\, j$-invariant of the elliptic curve
 $\, B(x, \, y) \, +p \,  = \, \, 0$.
 Let us just give here a simple example.

The telescoper of (\ref{alg2bisextrat}) with
$$\, \, \, \, B(x, \, y) \,   = \, \, \, 17 x y^2 +x^2 y^2 +5 x y^3 +13 y +y^2 +7 y x
+x^2 y^3 +2 x^3 y^3 +y^3 +3, $$
is an order-two linear differential operator having as solution:
\begin{eqnarray}
&& \hspace{-0.98in} \frac{1}{(15373+1656\,x)^{1/4}} \cdot \, \, _2F_1 
   \Bigl( [ {\frac{1}{12}},{\frac{5}{12}}],[1],
  \, \, {\frac { 64\, q_4(x)  }{ (15373 \, +1656\,x)^{3}}}\Bigr),
   \quad  \quad \hbox{where:} 
\\
  && \hspace{-0.98in}\quad \quad \quad \quad  q_4(x) \, \, = \, \, \,
     1143420561541 \, -475427554218\,x  \,  +67894132770\,{x}^{2}
\nonumber \\
&& \hspace{-0.98in} \quad \quad \quad \quad \quad \quad \quad 
   \, -2161434807\,{x}^{3}\, +22211523\,{x}^{4}.
\nonumber 
\end{eqnarray}
The Hauptmodul appearing in this solution is, as it should, of the form $ \, 1728/j$ where 
$\, j$ is the  $j$-invariant of the bicubic $ \, \, B(x,y) \, + p \, = \, \, 0$. 

\vskip .1cm

\subsection{A simpler $\, u$-extension.}
\label{extension2}

The telescoper of the rational function
\begin{eqnarray}
\label{rat33xzero}
\hspace{-0.98in}&&  \quad \quad  \quad \quad  \quad
 R(x, y, z, u)  \,  \, \, =  \, \, \, \,
{{ 1 } \over {  1\, \, +2\,y+3\,z   +5\,yz  \, +4\,{y}^{2}z    \, \, +3\,x \, u  }},               
\end{eqnarray}
or  the telescoper of the rational function  
\begin{eqnarray} 
\label{rat33}
\hspace{-0.98in}&&   \quad   \,
 R(x, y, z, u)  \, 
 \nonumber \\
 \hspace{-0.98in}&&  \quad \quad  \quad
 \, \, =  \, \, \, \,  {{ 1 } \over {  1\, \, +2\,y+3\,z   +5\,yz  \, +4\,{y}^{2}z    \, \, +3\,xu
   \, \, +(\alpha \, + \, \beta \,z  \, + \gamma \, y) \cdot \, {x}^{n} }},               
\end{eqnarray}
for $\, n \, = \, 1, \, 2, \, \cdots,  5, \, \cdots $, 
are identical (whatever the values of $\, \alpha,\, \beta, \, \gamma$).
This telescoper is an order-two linear differential operator
with the pullbacked hypergeometric solution: 
\begin{eqnarray}
\label{rat33sol}
  \hspace{-0.98in}&&  \quad   \quad
 (1 \, +312\, x \,  -1584 \,{x}^{2})^{-1/4}  \\
 \hspace{-0.98in}&&  \quad \quad  \quad \quad  \quad
                    \times \,
 _2F_1\Bigl([{{1} \over {12}}, \, {{5} \over {12}} ],  [1], \, \,
  \, {\frac { 5038848\,{x}^{3} \cdot \, (13248\,{x}^{3}+2928\,{x}^{2}+368\,x+1) }{
    (1584\,{x}^{2}-312\,x-1)^{3}}}  \Bigr).
\nonumber 
\end{eqnarray}
One can verify that the diagonals of the rational functions (\ref{rat33xzero}) and  (\ref{rat33})
are actually equal and correspond to the expansion of  the pullbacked hypergeometric solution
(\ref{rat33sol}):
\begin{eqnarray}
\label{rat33solexpa}
  \hspace{-0.98in}&&  \quad \quad \quad
                     1 \, \, -78\,x \,\,  +15606\,{x}^{2} \, \, -3888540\,{x}^{3} \,\,
                     +1069866630\,{x}^{4} \,\, -311621002308\,{x}^{5} 
 \nonumber \\
 \hspace{-0.98in}&&  \quad \quad  \quad \quad \quad \quad 
 \,\, +94190901642684\, x^6 \,\,\, \, -29220290149904568 \, x^7  \, \,  + \,  \, \, \cdots 
\end{eqnarray}

The elimination of $\, u$, with   $\, \, u \, = \, p/x/y/z \, $ in the vanishing
condition of the denominator of (\ref{rat33xzero}) gives the elliptic curve
\begin{eqnarray}
\label{ellipticrat33sol}
  \hspace{-0.98in}&&  \quad \quad  \quad \quad  \quad\quad  \quad
 4\,{y}^{3}{z}^{2} \, \, +5\,{y}^{2}{z}^{2} \, \,
 +2\,{y}^{2}z \,  \,+3\,y{z}^{2} \, \, +yz \, \, \, +3\,p \, \, \, = \, \, \, \, 0, 
\end{eqnarray}
which has  a $\, j$-invariant yielding the Hauptmodul
\begin{eqnarray}
\label{Hauptrat33sol}
  \hspace{-0.98in}&&  \quad \quad  \quad \quad \quad \quad \quad
 {\cal H} \, \, = \, \, \,  \,
{\frac { 5038848\,{p}^{3} \cdot \, (13248\,{p}^{3}+2928\,{p}^{2}+368\,p+1) }{
1584\,{p}^{2}-312\,p-1)^{3}}},           
\end{eqnarray}
which is precisely the pullback in (\ref{rat33sol}). 

\vskip .3cm

{\bf Remark:} note that the telescoper of
\begin{eqnarray} 
\label{rat33other}
  \hspace{-0.98in}&&    \, \, \, \, 
 R(x, y, z, u)  \,  \, \, =  \, \, \, \,
 {{ 1 } \over {
  1\, \, +2\,y+3\,z   +5\,yz  \, +4\,{y}^{2}z    \, \, +3\,x u \, \, \,\, \, + \, 11\, z\, x^2  \, + 7\, y\, x }},               
\end{eqnarray}
is a pretty large order-seven linear differential operator, however, this
operator is not the minimal order operator. The {\em minimal order}
linear differential operator
for the diagonal of (\ref{rat33other}) is actually the previous
order-two linear differential operator having the pullbacked hypergeometric solution (\ref{rat33sol}). 
One can verify directly that the {\em diagonal} of (\ref{rat33other}) is actually
the expansion (\ref{rat33solexpa}) for  the pullbacked hypergeometric solution (\ref{rat33sol}).

\vskip .2cm

\subsection{Examples with five, six, ... variables.}
\label{five}

Let us generalise  the four variables rational function (\ref{Ratfonc4G}) introducing
the {\em five} and {\em six}  variables rational functions
\begin{eqnarray}
\label{Ratfonc4G5}
  \hspace{-0.98in}&& \quad  \quad     
R(x, \, y, \, z, \, u, \, v)  \, \, \,  = \, \,  \quad 
  \\
\hspace{-0.98in}&&  
 \quad   \quad   \quad   \quad    \quad   
{{1} \over { a \,\, \,+ b_1 \, y \,  + \, c_1/y  + b_2 \, z  \, + \, c_2 /z \,  + d_1 \, y \, z  \,
  + \,e_1 \, y/z \,  + \,f_1 \, z/y \,  \, +g_1 \,u \, v \,x }},
\nonumber
\end{eqnarray}
and
\begin{eqnarray}
\label{Ratfonc4G6}
  \hspace{-0.98in}&& \quad   \quad    
R(x, \, y, \, z, \, u, \, v, \, w)  \, \, \,  = \, \,  \quad 
  \\
\hspace{-0.98in}&&  
 \quad   \quad   \, \,   \quad     \quad   
{{1} \over { a \,\, \,+ b_1 \, y \,  + \, c_1/y  + b_2 \, z  \, + \, c_2 /z \,  + d_1 \, y \, z  \,
  + \,e_1 \, y/z \,  + \,f_1 \, z/y \,  \,+g_1 \,u \, v \, w \, x }}.
\nonumber
\end{eqnarray}
Their telescopers are the same as the telescoper of (\ref{Ratfonc4G})
which annihilates the pullbacked hypergeometric function (\ref{solRatfonc4GG})
which Hauptmodul is associated with the elliptic (biquadratic) curve (\ref{biquadr}).

\end{document}